\documentclass[11pt]{article}

\usepackage{times}
\usepackage{float}

\usepackage{amsfonts,amstext,amsmath,amssymb,amsthm}

\pdfoutput=1

\usepackage{chngpage}
\usepackage{rotating}
\usepackage{graphicx}
\usepackage{epstopdf} 

\usepackage{caption}
\usepackage{subcaption}
\usepackage{empheq}

\usepackage{paralist}

\long\def\symbolfootnote[#1]#2{\begingroup%
\def\thefootnote{\fnsymbol{footnote}}\footnote[#1]{#2}\endgroup}

\usepackage{hyperxmp}

\usepackage[%
    pdftex,
    bookmarks=true,         
    unicode=false,          
    pdftoolbar=true,        
    pdfmenubar=true,        
    hyperfootnotes=false,   
    pdffitwindow=false,     
    pdfstartview={FitH},    
    pdftitle={Exact Distribution of the Generalized Shiryaev--Roberts Stopping Time Under the Minimax Brownian Motion Setup},    
    pdfauthor={Aleksey S. Polunchenko},     
    pdfsubject={Exact Distribution of the Generalized Shiryaev--Roberts Stopping Time Under the Minimax Brownian Motion Setup},   
    pdfcreator={Aleksey S. Polunchenko},   
    pdfkeywords={First passage times,Generalized Shiryaev--Roberts procedure,Kolmogorov forward equation,Markov diffusion processes,Method of separation of variables,Quickest change-point detection,Parabolic partial differential equations,Sequential analysis,Sturm--Liouville theory}, 
    pdfnewwindow=true,      
    colorlinks=false,       
    linkcolor=red,          
    citecolor=green,        
    filecolor=magenta,      
    urlcolor=cyan           
]{hyperref}

\usepackage[T1]{fontenc}

\usepackage{titlesec} 

\titleformat{\section}{\large\bfseries\uppercase}{\thesection.}{.5em}{}
\titlespacing*{\section}{0pt}{*3}{*2}
\titleformat{\subsection}{\normalfont\bfseries}{\thesubsection.}{.5em}{}
\titlespacing*{\subsection}{0pt}{*3}{*2}
\titleformat{\subsubsection}{\normalfont\bfseries}{\thesubsubsection.}{.5em}{}
\titlespacing*{\subsubsection}{0pt}{*3}{*2}

\numberwithin{equation}{section}

\usepackage[sort&compress]{natbib}

\renewcommand{\Pr}{\mathbb{P}} 
\DeclareMathOperator{\EV}{\mathbb{E}} 

\DeclareMathOperator{\LR}{\Lambda}

\DeclareMathOperator{\ADD}{ADD}
\DeclareMathOperator{\ARL}{ARL}

\DeclareMathOperator{\Ei}{Ei}

\newcommand{\T}{T}

\renewcommand{\le}{\leqslant} 
\renewcommand{\ge}{\geqslant}

\newcommand{\abs}[1]{\left\vert#1\right\vert}

\DeclareMathOperator{\One}{\mathchoice{\rm 1\mskip-4.2mu l}{\rm 1\mskip-4.2mu l}{\rm 1\mskip-4.6mu l}{\rm 1\mskip-5.2mu l}}
\newcommand{\indicator}[1]{{\One_{\left\{#1\right\}}}}

\theoremstyle{plain} 

\theoremstyle{plain}
\newtheorem{remark}{Remark}[section]

\graphicspath{{./gfx/}}

\begin{document}

\title{\textbf{\Large Exact Distribution of the Generalized Shiryaev--Roberts Stopping Time Under the Minimax Brownian Motion Setup}}

\date{}
\author{}
\maketitle

\begin{center}
\null\vskip-2cm\author{
\textbf{\large Aleksey\ S.\ Polunchenko}\\
Department of Mathematical Sciences, State University of New York at Binghamton,\\Binghamton, New York, USA
}
\end{center}
%
%
%
%
\symbolfootnote[0]{\normalsize\hspace{-0.6cm}Address correspondence to A.\ S.\ Polunchenko, Department of Mathematical Sciences, State University of New York (SUNY) at Binghamton, 4400 Vestal Parkway East, Binghamton, NY 13902--6000, USA; Tel: +1 (607) 777--6906; Fax: +1 (607) 777--2450; E-mail:~\href{mailto:aleksey@binghamton.edu}{aleksey@binghamton.edu}.}\\
%
%
{\small\noindent\textbf{Abstract:} We consider the quickest change-point detection problem where the aim is to detect the onset of a pre-specified drift in ``live''-monitored standard Brownian motion; the change-point is assumed unknown (nonrandom). The object of interest is the distribution of the stopping time associated with the Generalized Shryaev--Roberts (GSR) detection procedure set up to ``sense'' the presence of the drift in the Brownian motion under surveillance. Specifically, we seek the GSR stopping time's survival function (the tail probability that no alarm is triggered by the GSR procedure prior to a given point in time), and distinguish two scenarios:\begin{inparaenum}[\itshape(a)] \item when the drift never sets in (pre-change regime) and \item when the drift is in effect {\em ab initio} (post-change regime)\end{inparaenum}. Under each scenario, we obtain a closed-form formula for the respective survival function, with the GSR statistic's (deterministic) nonnegative headstart assumed arbitrarily given. The two formulae are found analytically, through direct solution of the respective Kolmogorov forward equation via the Fourier spectral method to achieve separation of the spacial and temporal variables. We then exploit the obtained formulae numerically and characterize the pre- and post-change distributions of the GSR stopping time depending on three factors:\begin{inparaenum}[\itshape(a)]\item magnitude of the drift, \item detection threshold, and \item the GSR statistic's headstart\end{inparaenum}.
}
\\ \\
%
%
\\ \\
%
%
%
%
%
%
%
%
%

\section{Introduction}
\label{sec:intro}

Sequential (quickest) change-point detection is concerned with the development and evaluation of dependable statistical procedures for early detection of unanticipated changes that may (or may not) occur online in the characteristics of a ``live''-monitored (random) process. Specifically, the process is ``inspected'' continuously so as to keep its characteristics as intended, which is achieved by ``sounding'' an alarm as soon as the process starts to behave otherwise; the challenge is to ``sound'' the alarm as quickly as is possible within an {\it a~priori} set tolerable level of the false positive risk. See, e.g.,~\cite{Shiryaev:Book78}, \cite{Basseville+Nikiforov:Book93}, \cite{Poor+Hadjiliadis:Book09}, \cite{Veeravalli+Banerjee:AP2013}, \cite[Part~II]{Tartakovsky+etal:Book2014} and the references therein.

A change-point detection procedure is identified with a stopping time, $\T$, that is adapted to the filtration, $(\mathcal{F}_t)_{t\ge0}$, generated by the observed process, $(X_t)_{t\ge0}$; the semantics of $\T$ is that it constitutes a rule to stop and declare that the statistical profile of the observed process may have (been) changed. A ``good'' (i.e., optimal or nearly optimal) detection procedure is one that minimizes (or nearly minimizes) the desired detection delay penalty-function, subject to a constraint on the false alarm risk. For an overview of the major optimality criteria see, e.g.,~\cite{Tartakovsky+Moustakides:SA10}, \cite{Polunchenko+Tartakovsky:MCAP2012}, \cite{Polunchenko+etal:JSM2013},~\cite{Veeravalli+Banerjee:AP2013}, and~\cite[Part~II]{Tartakovsky+etal:Book2014}.

\begin{figure}[!htb]
    \centering
    \includegraphics[width=0.9\textwidth]{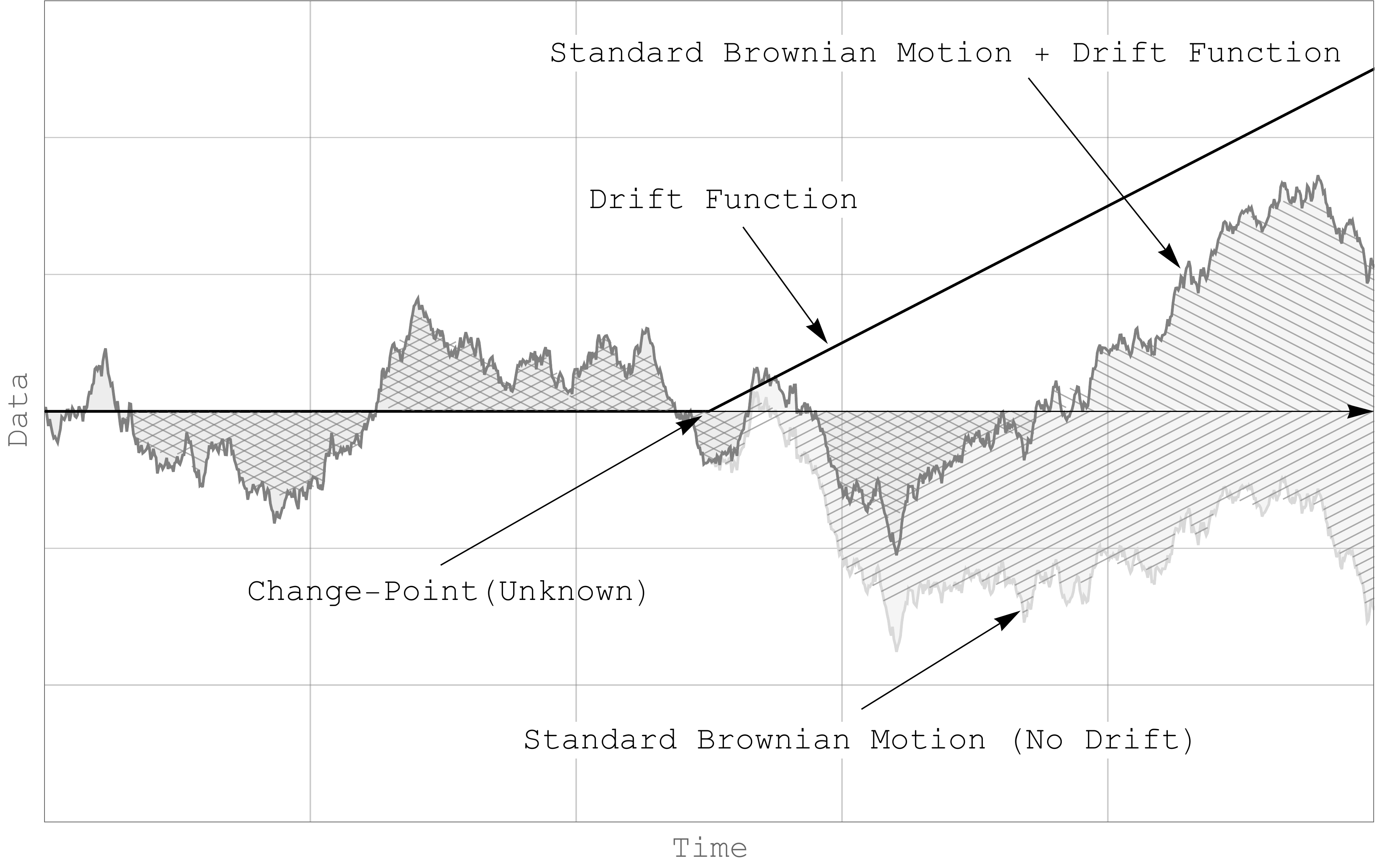}
    \caption{Standard Brownian motion gaining a persistent drift at an unknown time moment.}
    \label{fig:BM-change-point-scenario}
\end{figure}

This work concentrates on the popular minimax setup of the basic change-point detection problem where the observed process, $(X_t)_{t\ge0}$, is standard Brownian motion that at an unknown (nonrandom) time moment $\nu\ge0$---referred to as the change-point---may (or may not) experience an abrupt and permanent change in the drift, from a value of zero initially, i.e., $\EV[X_t]=0$ for $t\in[0,\nu]$, to a known value $\mu\neq0$ following the change-point, i.e., $\EV[X_t]=\mu t$ for $t>\nu$. This is schematically illustrated in Figure~\ref{fig:BM-change-point-scenario}. The goal is to find out---as quickly as is possible within an {\it a~priori} set level of the ``false positive'' risk---whether the drift of the process is no longer zero. See, e.g.,~\cite{Pollak+Siegmund:B85},~\cite{Shiryaev:RMS1996,Shiryaev:Bachelier2002},~\cite{Moustakides:AS2004},~\cite{Shiryaev:MathEvents2006},~\cite{Feinberg+Shiryaev:SD2006},~\cite{Burnaev+etal:TPA2009}, and~\cite[Chapter~5]{Shiryaev:Book2011}.

More formally, under the above Brownian motion change-point scenario, the observed process, $(X_t)_{t\ge0}$, is governed by the stochastic differential equation (SDE):
\begin{align}\label{eq:BM-change-point-model}
dX_t
&=
\mu\indicator{t>\nu}dt+dB_t,\;t\ge0,\;\text{with}\;X_0=0,
\end{align}
where $(B_t)_{t\ge0}$ is standard Brownian motion (i.e., $\EV[dB_t]=0$, $\EV[(dB_t)^2]=dt$, and $B_0=0$), $\mu\neq0$ is the {\em known} post-change drift value, and $\nu\in[0,\infty]$ is the unknown (nonrandom) change-point; here and onward, the notation $\nu=0$ ($\nu=\infty$) is to be understood as the case when the drift is in effect {\it ab initio} (or never, respectively).

The standard way to perform change-point detection under model~\eqref{eq:BM-change-point-model} has been to employ Page's~\citeyearpar{Page:B54} Cumulative Sum (CUSUM) ``inspection scheme''. The choice to use the CUSUM procedure may be justified by the fact (established by~\citealt{Beibel:AS1996}, by~\citealt{Shiryaev:RMS1996}, and by~\citealt{Moustakides:AS2004}) that the CUSUM ``inspection scheme'' is strictly minimax-optimal in the sense of~\cite{Lorden:AMS71}; the discrete-time equivalent of this result was first established by~\cite{Moustakides:AS86}, although an alternative proof was later also offered by~\cite{Ritov:AS90} who exploited a game-theoretic argument.

However, when one is interested in minimax optimality as defined by~\cite{Pollak:AS85}, a sensible alternative to using the CUSUM procedure would be to devise the Generalized Shiryaev--Roberts (GSR) procedure. The latter is due to~\cite{Moustakides+etal:SS11}, and is a headstarted (hence, more general) version of the classical quasi-Bayesian Shiryaev--Roberts (SR) procedure that emerged from the independent work of~\cite{Shiryaev:SMD61,Shiryaev:TPA63} and that of~\cite{Roberts:T66}. With Pollak's~\citeyearpar{Pollak:AS85} definition of minimax optimality in mind, the motivation to prefer the GSR procedure over the CUSUM procedure stems from the results obtained (for the discrete-time analogue of the problem) by~\cite{Tartakovsky+Polunchenko:IWAP10} and by~\cite{Polunchenko+Tartakovsky:AS10}, and then also by~\cite{Tartakovsky+etal:TPA2012} who showed that the GSR procedure with a carefully designed headstart may be faster (in Pollak's~\citeyear{Pollak:AS85} minimax sense) than the CUSUM procedure; as a matter of fact,~\cite{Tartakovsky+Polunchenko:IWAP10} and~\cite{Polunchenko+Tartakovsky:AS10} proved the GSR procedure (with a ``finetuned'' headstart) to be not only faster (in Pollak's~\citeyear{Pollak:AS85} minimax sense) but {\em the} fastest (i.e., the best one can do, again in Pollak's~\citeyear{Pollak:AS85} minimax sense) in two specific (discrete-time) scenarios. For an attempt to extend these results to the (continuous-time) Brownian motion scenario~\eqref{eq:BM-change-point-model}, see, e.g.,~\cite{Burnaev:ARSAIM2009}.

To formally state the problem addressed in this work let us first introduce the GSR procedure. Let $\Pr_{\infty}$ ($\Pr_{0}$) denote the probability measure (distribution law) generated by the observed process, $(X_t)_{t\ge0}$, under the assumption that $\nu=\infty$ ($\nu=0$); note that $\Pr_{\infty}$ is the Wiener measure. Let $\left.\Pr_{\infty}\right|_{\mathcal{F}_{t}}$ ($\left.\Pr_{0}\right|_{\mathcal{F}_t}$) be the restriction of probability measure $\Pr_{\infty}$ ($\Pr_{0}$) to the filtration $\mathcal{F}_{t}$. Define
\begin{align*}
\LR_{t}
&\triangleq
\dfrac{d\left.\Pr_{0}\right|_{\mathcal{F}_{t}}}{d\left.\Pr_{\infty}\right|_{\mathcal{F}_{t}}},\;t\ge0,
\end{align*}
i.e., the Radon--Nikod\'{y}m derivative of $\left.\Pr_{0}\right|_{\mathcal{F}_{t}}$ with respect to $\left.\Pr_{\infty}\right|_{\mathcal{F}_{t}}$. It is well-known that for the Brownian motion scenario under consideration
\begin{align*}
\LR_{t}
&=
\exp\left\{\mu X_{t}-\dfrac{\mu^2}{2} t\right\},\;\text{so that}\; d\LR_{t}=\mu\LR_{t} dX_{t},\;\LR_0=1;
\end{align*}
cf., e.g.,~\cite{Shiryaev:Book1999},~\cite[Formula~(15),~p.~378, and Formula~(16),~p.~379]{Shiryaev:MathEvents2006}, \cite[Formula~(4.31),~p.~49]{Shiryaev:Book2011}, and~\cite{Liptser+Shiryaev:Book2001}.

The process $\{\LR_{t}\}_{t\ge0}$ is the likelihood ratio to test the hypothesis $\mathcal{H}_0\colon\nu=0$ against the alternative $\mathcal{H}_\infty\colon\nu=\infty$, and is the key ingredient of the GSR statistic, conventionally denoted as $(R_t^r)_{t\ge0}$. Specifically, tailored to the Brownian motion scenario at hand, the GSR statistic, $(R_t^r)_{t\ge0}$, is of the form
\begin{align}\label{eq:Rt_r-def}
\begin{split}
R_{t}^r
&\triangleq
r\LR_{t}+\int_0^t\dfrac{\LR_{t}}{\LR_s}\,ds\\
&=
r\exp\left\{\mu X_{t}-\dfrac{\mu^2t}{2}\right\}+\int_0^t\exp\left\{\mu(X_{t}-X_s)-\dfrac{\mu^2(t-s)}{2}\right\}ds,\;t\ge0,
\end{split}
\end{align}
where $R_0^r=r\ge0$ is the headstart (a deterministic point selected so as to optimize the GSR procedure's performance; see, e.g.,~\citealt{Tartakovsky+Polunchenko:IWAP10,Polunchenko+Tartakovsky:AS10,Moustakides+etal:SS11,Tartakovsky+etal:TPA2012,Polunchenko+Sokolov:EnT2014}). When $R_0^r=r=0$, it is said that the GSR statistic has no headstart, in which case it is equivalent to the classical SR statistic. Consequently, the GSR procedure whose statistic has no headstart is equivalent to the classical SR procedure. Hence, the labels ``{\em Generalized} SR statistic'' and ``{\em Generalized} SR procedure'', which appear to have both been coined by~\cite{Tartakovsky+etal:TPA2012}.

The GSR procedure calls for stopping as soon as the GSR statistic $(R_t^r)_{t\ge0}$ either hits or exceeds a certain flat level $A>0$ known as the detection threshold. More formally, the GSR procedure is identified with the stopping time:
\begin{align}\label{eq:T-GSR-def}
\mathcal{S}_A^r
&\triangleq
\inf\{t\ge0\colon R_{t}^r\ge A\}\;\text{such that}\;\inf\{\varnothing\}=\infty,
\end{align}
where $(R_{t}^r)_{t\ge0}$ is the GSR statistic given by~\eqref{eq:Rt_r-def}. The detection threshold $A>0$ is selected in advance so as to control the ``false positive'' risk within acceptable margins. Due to path-continuity of the GSR statistic $(R_{t}^r)_{t\ge0}$, the inequality ``$R_{t}^r\ge A$'' in the above definition of the GSR stopping time $\mathcal{S}_A^r$ may be replaced with equality $R_{t}^r=A$. We note that this is in stark contrast with the discrete-time version of the problem where the GSR statistic is not path-continuous and, as such, is bound to always overshoot the detection threshold; this phenomenon is known as the ``overshoot problem''. We also note that from now on we shall require the headstart $R_0^r=r\ge0$ to come from the interval $[0,A]$ rather than from the interval $[0,+\infty)$, because $\mathcal{S}_A^r=0$ for $R_0^r=r>A\,(>0)$, as can be easily deduced from~\eqref{eq:T-GSR-def} and from~\eqref{eq:Rt_r-def}; in fact, $\mathcal{S}_A^r=0$ for $R_0^r=r=A\,(>0)$ as well, so the detection threshold $A$ is included into the state space of the GSR statistic $(R_{t}^r)_{t\ge0}$ for convenience.

We are now in a position to formulate the specific problem addressed in this paper: to obtain analytically closed-form formulae for the tail probabilities $\Pr_{\infty}(\mathcal{S}_A^r\ge t)$ and $\Pr_{0}(\mathcal{S}_A^r\ge t)$ for any $t\ge0$ and $R_0^r=r\in[0,A]$, with $A>0$. Put otherwise, we are interested in the survival function of the GSR stopping time $\mathcal{S}_A^r$ in two cases:\begin{inparaenum}[\itshape(a)]\item when the observed Brownian motion stays drift-free indefinitely (i.e., $\nu=\infty$) and \item when the Brownian motion is affected by drift {\em ab initio} (i.e., $\nu=0$)\end{inparaenum}. The former scenario corresponds to the pre-change (or pre-drift) regime, and the latter scenario corresponds to the post-change (or post-drift) regime. To the best of our knowledge, neither of the two survival functions has heretofore been obtained explicitly. However, in the discrete-time setup, the problem has been solved by~\cite{Moustakides+etal:SS11} who proposed a general numerical framework to compute a broad range of performance metrics (including the two survival functions) not only for the GSR procedure but also for a whole family of detection procedures with Markovian detection statistics. Moreover, for the GSR procedure specifically, the framework of~\cite{Moustakides+etal:SS11} has been recently improved in terms of accuracy and efficiency by~\cite{Polunchenko+etal:SA2014,Polunchenko+etal:ASMBI2014} and then also by~\cite[Chapter~3]{Du:PhD-Thesis2015}. We also note that, in the discrete-time setup, it is rarely a possibility that the performance of a detection procedure can be found analytically and in a closed-form. The reason is the aforementioned ``overshoot problem''. Hence, the ``solution'' obtained by~\cite{Moustakides+etal:SS11} and then ``refined'' by~\cite{Polunchenko+etal:SA2014,Polunchenko+etal:ASMBI2014} and by~\cite[Chapter~3]{Du:PhD-Thesis2015} is only {\em numerical}, although with a controllably small error. By contrast, the continuous-time model~\eqref{eq:BM-change-point-model} is ``immune'' to the overshoot problem, so the expressions that we obtain in this work for the $\Pr_{\infty}$- and $\Pr_{0}$-survival functions of the GSR stopping time $\mathcal{S}_A^r$ given by~\eqref{eq:T-GSR-def} are {\em exact}.

The remainder of the paper is organized as follows. We begin in Section~\ref{sec:problem+preliminaries} with formally setting up two partial differential equations (PDEs) to then recover the sought survival functions from:\begin{inparaenum}[\itshape(a)]\item one equation corresponding to the pre-drift regime ($\nu=\infty$), and \item one equation corresponding to the post-drift regime ($\nu=0$)\end{inparaenum}. Both PDEs are Kolmogorov forward equations that are each subject to one initial temporal condition and two spacial boundary conditions---one at each of the two end-points of the strip $[0,A]$ to which the GSR statistic $(R_{t}^r)_{t\ge0}$ is confined by virtue of the definition~\eqref{eq:T-GSR-def} of the GSR stopping time $\mathcal{S}_A^r$. While the two Kolmogorov forward equations are different (one assumes that $\nu=\infty$ and the other one assumes that $\nu=0$), they both can be treated {\em simultaneously}, for the two can be combined into one {\em master equation} by introducing an auxiliary ``boolean'' variable (equal to zero when $\nu=\infty$, and to one when $\nu=0$) through which the master equation can be quickly turned into one of the two regime-specific forms. Of course the same boolean ``switch'' also allows to unify the initial and boundary conditions corresponding to different regimes. All this is detailed in Section~\ref{sec:problem+preliminaries} as well. We conclude Section~\ref{sec:problem+preliminaries} with a brief outline the so-called Fourier method to {\em analytically} solve the master equation by means of separating the temporal and spacial variables. The centerpiece of this work is Section~\ref{sec:main-result}, where we devise the Fourier method and solve the master equation explicitly, and obtain exact closed-form formulae for the $\Pr_{\infty}$- and $\Pr_{0}$-survival functions of the GSR stopping time. The obtained formulae are then exploited numerically in Section~\ref{sec:numerical-study} where we offer a numerical study aimed at characterizing the distribution of the GSR stopping time in the pre- and post-drift regimes. To carry out the study, we implemented the obtained formulae in {\em Mathematica}, the popular software package developed by Wolfram Research, Inc. as a programming environment for scientific computing. For more information about Mathematica and Wolfram Research, Inc. see on the Web at \url{www.wolfram.com}. Lastly, Section~\ref{sec:conclusion} summarizes the entire paper.

\section{Preliminaries}
\label{sec:problem+preliminaries}

This section's aim is to briefly outline the principal approach we intend to undertake in the next (main) section to solve the problem we have set out to address in this work, i.e., find $\Pr_{\infty}(\mathcal{S}_A^r\ge t)$ and $\Pr_{0}(\mathcal{S}_A^r\ge t)$ analytically and in a closed-form for any $t\ge0$ and $R_0^r=r\in[0,A]$ with $A>0$ given; recall that the GSR procedure's detection statistic $(R_t^r)_{t\ge0}$ is given by~\eqref{eq:Rt_r-def}, and that the GSR procedure's stopping time $\mathcal{S}_A^r$ is given by~\eqref{eq:T-GSR-def}.

To get started, observe that, as an immediate implication of It\^{o}'s formula applied to the definition~\eqref{eq:Rt_r-def} of the GSR statistic $(R_t^r)_{t\ge0}$, the latter's $\Pr_\infty$-differential is $dR_t^r=dt+\mu R_t^r dB_t$; cf.,~e.g.,~\cite[Formula~(4),~p.~269]{Pollak+Siegmund:B85},~\cite[Formula~(46),~p.~386]{Shiryaev:MathEvents2006} or~\cite[Formula~(1.17),~p.~449]{Feinberg+Shiryaev:SD2006}. Likewise, the respective $\Pr_0$-differential can be seen to be $dR_t^r=(1+\mu^2 R_t^r)dt+\mu R_t^r dB_t$; cf.,~e.g.,~\cite[Formula~(4),~p.~269]{Pollak+Siegmund:B85}. Since either differential is such that the instantaneous drift function and the instantaneous diffusion function both do not depend on time, one can conclude that the GSR statistic $(R_t^r)_{t\ge0}$ is a time-homogeneous Markov diffusion, whether in the pre-drift regime or in the post-drift regime. More importantly, the form of the $\Pr_\infty$-differential and that of the $\Pr_0$-differential are similar enough to be conveniently combined into one:
\begin{align}\label{eq:Rt_r-combined-SDE}
dR_t^r
&=
(1+\theta \mu^2 R_t^r)dt+\mu R_t^r dB_t,
\end{align}
where $\theta$ is either $0$ or $1$ so that $\theta^2=\theta$. Note now that, on the one hand, setting $\theta=0$ in the foregoing differential turns it into the $\Pr_\infty$-differential of $(R_t^r)_{t\ge0}$, and, on the other hand, if $\theta=1$, then the differential~\eqref{eq:Rt_r-combined-SDE} becomes the $\Pr_0$-differential of $(R_t^r)_{t\ge0}$. Let $b(x)\triangleq 1+\theta\mu^2 x$ and $\sqrt{a(x)}\triangleq\mu x$ denote the corresponding instantaneous drift function and diffusion coefficient, respectively.

Next, define
\begin{align*}
p_{\theta}(y,t|x,s)
&\triangleq
\dfrac{\partial}{\partial y}
\begin{cases}
\Pr_{\infty}(R_t^r\le y,\mathcal{S}_A^r\ge t|R_s^r=x),&\text{if $\theta=0$};\\
\Pr_{0}(R_t^r\le y,\mathcal{S}_A^r\ge t|R_s^r=x),&\text{if $\theta=1$},
\end{cases}
\end{align*}
where $0\le s\le t$, i.e., $p_{\theta}(y,t|x,s)$ is the transition probability density of the time-homogeneous Markov diffusion $(R_t^r)_{t\ge0}$ joint with the event that the respective GSR stopping time $\mathcal{S}_A^r$ does not terminate the diffusion $(R_t^r)_{t\ge0}$ prior to a given time point $t\ge0$. Since $(R_t^r)_{t\ge0}$ is time-homogeneous, $p_{\theta}(y,t|x,s)$ depends on $s$ and $t$ only through the difference $t-s\ge0$. Therefore, it suffices to consider only $p_{\theta}(x,t|r)\triangleq p_{\theta}(x,t|r,0)$, because, by definition, $R_0^r=r$. More concretely,
\begin{align}\label{eq:p-theta-def}
p_{\theta}(x,t|r)
&\triangleq
\dfrac{\partial}{\partial x}
\begin{cases}
\Pr_{\infty}(R_t^r\le x,\mathcal{S}_A^r\ge t),&\text{if $\theta=0$};\\
\Pr_{0}(R_t^r\le x,\mathcal{S}_A^r\ge t),&\text{if $\theta=1$},
\end{cases}
\end{align}
where $t\ge 0$. At this point note that since
\begin{align}\label{eq:p-theta-to-survfun}
\Pr_{\infty}(\mathcal{S}_A^r\ge t)
&=
\int_{0}^{A} p_{0}(x,t|r)\,dx
\;\;\text{and}\;\;
\Pr_{0}(\mathcal{S}_A^r\ge t)
=
\int_{0}^{A} p_{1}(x,t|r)\,dx,
\end{align}
where $t\ge0$ and $r\in[0,A]$ with $A>0$, finding $p_{\theta}(x,t|r)$ explicitly for both $\theta=0$ and $\theta=1$ can be seen to be the main stepping stone toward our goal of getting closed-form expressions for $\Pr_{\infty}(\mathcal{S}_A^r\ge t)$ and $\Pr_{0}(\mathcal{S}_A^r\ge t)$ for all $t\ge0$ and $R_0^r=r\in[0,A]$ with $A>0$ given; we remark parenthetically that $\Pr_{\infty}(\mathcal{S}_A^r\ge 0)\equiv1$ and $\Pr_{0}(\mathcal{S}_A^r\ge 0)\equiv1$ for any $r\ge 0$, which is a trivial consequence the definition~\eqref{eq:T-GSR-def} of the GSR stopping time $\mathcal{S}_A^r$.

Since we have now reduced the problem to that of finding $p_{0}(x,t|r)$ and $p_{1}(x,t|r)$ given by~\eqref{eq:p-theta-def}, let us now briefly explain how we plan to find $p_{0}(x,t|r)$ and $p_{1}(x,t|r)$. To that end, the key is exploit the general framework outlined in~\cite[Chapter~6]{Schuss:Book2010} to treat stopped diffusions. See also, e.g.,~\cite[Part~1,~Chapter~4]{Stratonovich:Book1961},~\cite[Chapter~26]{Tikhonov+Mironov:Book1977}, and~\cite[Chapter~5]{Gardiner:Book1985}. Specifically, consider the general diffusion process $(Y_t)_{t\ge0}$ that follows the SDE:
\begin{align*}
dY_t
&=
b_{Y}(Y_t)\,dt+\sqrt{a_{Y}(Y_t)}\,dB_t,\;t\ge0,\;Y_0=y_0,
\end{align*}
where the instantaneous drift function $b_{Y}(y)$ and the instantaneous diffusion coefficient $\sqrt{a_{Y}(y)}$ are both sufficiently smooth. Define the stopping time $\mathcal{T}_B\triangleq\inf\{t\ge 0\colon Y_t\ge B\}$ such that $\inf\{\varnothing\}=\infty$, where $B>0$ is a given threshold. Then, according to~\cite[Chapter~6]{Schuss:Book2010}, the transition density $p^Y(y,t|y_0)\triangleq d\Pr(Y_t\le y,\mathcal{T}_B\ge t)/dy$ simultaneously satisfies two PDEs. Specifically, on the one hand, the density $p^Y(y,t|y_0)$ satisfies the Kolmogorov forward equation
\begin{align}\label{eq:general-Kolmogorov-fwd-PDE}
\frac{\partial}{\partial t}p^Y(y,t|y_0)
&=
-\frac{\partial}{\partial y}\big[b_Y(y)\,p^Y(y,t|y_0)\big]
+\frac{1}{2}\frac{\partial^2}{\partial y^2}\big[a_Y(y)\,p^Y(y,t|y_0)\big],
\end{align}
which, as a PDE of order one in time $t$ (temporal variable) and order two in $x$ (spacial variable), is to be complemented by one {\em initial} temporal condition and two spacial boundary conditions. On the other hand, the density $p^Y(y,t|y_0)$ also satisfies the Kolmogorov backward equation
\begin{align}\label{eq:general-Kolmogorov-back-PDE}
-\frac{\partial}{\partial t}p^Y(y,t|y_0)
&=
b_Y(y)\frac{\partial}{\partial y}\big[p^Y(y,t|y_0)\big]
+\frac{a_Y(y)}{2}\frac{\partial^2}{\partial y^2}\big[p^Y(y,t|y_0)\big],
\end{align}
which, again as a PDE of order one in $t$ and order two in $x$, is to be complemented by one {\em terminal} temporal condition and two spacial boundary conditions. The two equations~\eqref{eq:general-Kolmogorov-fwd-PDE}--\eqref{eq:general-Kolmogorov-back-PDE} are adjoint to each other, and both stem from the seminal work of~\cite{Kolmogorov:MA1931}; incidentally, the forward equation~\eqref{eq:general-Kolmogorov-fwd-PDE} is also important in physics (viz. in quantum mechanics), where it is known as the Fokker--Plank equation, after~\cite{Fokker:APh1914} and~\cite{Planck:PAW1917}, who arrived at the equation before~\cite{Kolmogorov:MA1931}, although using different techniques and motivated by different considerations.

The aforementioned mutual ``adjointness'' of the two Kolmogorov equations~\eqref{eq:general-Kolmogorov-fwd-PDE}--\eqref{eq:general-Kolmogorov-back-PDE} can be illustrated as follows. Introduce the differential operator
\begin{align*}
\mathcal{G}^{*}
&\triangleq
\dfrac{1}{2} a_Y(y)\dfrac{\partial^2}{\partial y^2}+b_Y(y)\dfrac{\partial}{\partial y}
\end{align*}
with
\begin{align*}
\mathcal{G}
&\triangleq
\dfrac{1}{2}\dfrac{\partial^2}{\partial y^2}a_Y(y)-\dfrac{\partial}{\partial y}b_Y(y)
\end{align*}
being the corresponding adjoint operator. Then in terms of the operators $\mathcal{G}$ and $\mathcal{G}^*$, the forward equation~\eqref{eq:general-Kolmogorov-fwd-PDE} can be compactly written as $[\mathcal{G}\circ{p}^{Y}](y,t)=\partial p^{Y}(y,t)/\partial t$, and the operator form of the backward equation~\eqref{eq:general-Kolmogorov-back-PDE} is $[\mathcal{G}^{*}\circ {p}^{Y}](y,t)=-\partial p^{Y}(y,t)/\partial t$. One of the fundamental properties of the operators $\mathcal{G}$ and $\mathcal{G}^*$ is that they can be parameterized as follows
\begin{align}\label{eq:operator-G-param}
\mathcal{G}^{*}
&=
\dfrac{1}{\mathfrak{m}(y)}\dfrac{d}{dy}\dfrac{1}{\mathfrak{s}(y)}\dfrac{d}{dy}
\;\;\text{and}\;\;
\mathcal{G}
=
\dfrac{d}{dy}\dfrac{1}{\mathfrak{s}(y)}\dfrac{d}{dy}\dfrac{1}{\mathfrak{m}(y)},
\end{align}
where
\begin{align}\label{eq:speed+scale-def}
\mathfrak{s}(x)
&\triangleq
\exp\left\{-\int\dfrac{2 b_Y(x)}{a_Y(x)}\,dx\right\}
\;\;\text{and}\;\;
\mathfrak{m}(x)\triangleq\dfrac{2}{a_Y(x)\,\mathfrak{s}(x)}.
\end{align}
i.e., $\mathfrak{s}(x)$ is the solution of the ODE $[\mathcal{G}^{*}\circ\mathfrak{s}](x)=0$ while $\mathfrak{m}(x)$ satisfies the ODE $[\mathcal{G}\circ\mathfrak{m}](x)=0$. It is now direct to see from~\eqref{eq:operator-G-param} that $\mathcal{G}^{*}$ and $\mathcal{G}$ are self-adjoint with respect to $\mathfrak{s}(x)$ and $\mathfrak{m}(x)$, respectively; cf., e.g.,~\cite{Borodin+Salminen:Book2002}. The former function is known as the scale measure, while the function $\mathfrak{m}(x)$ is referred to as the speed measure.

Since the two Kolmogorov equations are mutually adjoint, it follows that either one alone is sufficient to fully characterize the density $p^Y(y,t|y_0)$, provided, however, that the initial (respectively, terminal, if it's the backward equation) temporal condition and the two spacial boundary conditions are properly specified. As a matter of fact, it is the initial (respectively, terminal, if it's the backward equation) condition and the two boundary conditions that not only make the corresponding PDE a complete problem, but also determine the nature of the solution. Since in this work we wish to deal with the forward equation, let us from now on concentrate exclusively on the forward equation~\eqref{eq:general-Kolmogorov-fwd-PDE}.

For the forward equation~\eqref{eq:general-Kolmogorov-fwd-PDE} the initial temporal equation is straightforward: $\lim_{t\to0+}p^Y(y,t|y_0)=\delta(y-y_0)$, where here and onward $\delta(x)$ denotes the Dirac delta function, so that ``$\lim_{t\to0+}p^Y(y,t|y_0)=\delta(y-y_0)$'' is to be understood as equality of distributions. This initial condition merely states that at time zero the process $(Y_t)_{t\ge0}$ is purely deterministic with the entire ``probability mass'' concentrated at one given point $Y_0=y_0$. The two spacial boundary conditions are not as straightforward, because they depend on the particular type of boundaries involved: absorbing, reflective, ``sticky'', natural, entrance, etc. See, e.g.,~\cite[Part~1,~Chapter~4]{Stratonovich:Book1961},~\cite[Chapters~12~\&~26]{Tikhonov+Mironov:Book1977},~\cite[Chapter~5]{Gardiner:Book1985}, and~\cite[Chapter~II,~pp.~14--15]{Borodin+Salminen:Book2002}. Since we are interested in the case when the process $(Y_t)_{t\ge0}$ is restricted to the strip $[0,B]$, we have two boundaries to consider: one at zero and one at $B>0$. In our case, the latter is an absorbing (``killing'') boundary, so that according to~\cite[Chapter~6]{Schuss:Book2010} the corresponding boundary condition is $p^Y(B,t|y_0)=0$ for all $y_0$. For the left end-point of the interval $[0,B]$, we are interested in the case when it is an entrance boundary, which means the process may enter its state space $[0,B]$ through zero but then will never return to it. This is precisely the type of boundary that zero is for the GSR diffusion $(R_t^r)_{t\ge0}$. According to~\cite[Chapter~6]{Schuss:Book2010}, for such boundaries the boundary condition is of the form
\begin{align}\label{eq:Ybnd-cond-0}
\lim_{y\to0+}\left[\frac{1}{\mathfrak{s}(y)}\frac{\partial}{\partial y} \frac{p^{Y}(y,t|y_0)}{\mathfrak{m}(y)}\right]&=0,\;\; y_0\in[0,B),
\end{align}
where $\mathfrak{s}(x)$ and $\mathfrak{m}(x)$ are, respectively, the scale and speed measures given by~\eqref{eq:speed+scale-def}.

It is straightforward to tailor the above brief account of the results presented in~\cite[Chapter~6]{Schuss:Book2010} to the GSR diffusion $(R_t^r)_{t\ge0}$ and the corresponding stopping time $\mathcal{S}_A^r$. Specifically, the density $p_{\theta}(x,t|r)$ defined in~\eqref{eq:p-theta-def} can be seen to satisfy the following Kolmogorov forward equation
\begin{align}\label{eq:master-eqn}
\frac{\partial}{\partial t}p_{\theta}(x,t|r)
&=
-\frac{\partial}{\partial x}\big[(1+\theta\mu^2 x)p_{\theta}(x,t|r)\big]
+\frac{\mu^2}{2}\frac{\partial^2}{\partial x^2}\big[x^2p_{\theta}(x,t|r)\big],\; t\ge0,\; x,r\in[0,A],
\end{align}
subject to\begin{inparaenum}[\itshape(a)]\item the initial condition $\lim_{t\to0+}p_{\theta}(x,t|r)=\delta(x-r)$ valid for all $x$, and \item two boundary conditions---one at $x=0$ (or as $x\to0+$) and one at the absorbing (or ``cemetery'') boundary $x=A$\end{inparaenum}. The former boundary condition is akin to~\eqref{eq:Ybnd-cond-0}, and is of the form:
\begin{align}\label{eq:bnd-cond-00}
\lim_{x\to0+}\left[\frac{1}{\mathfrak{s}(x)}\frac{\partial}{\partial x} \frac{p_{\theta}(x,t|r)}{\mathfrak{m}(x)}\right]&=0,\;\; r\in[0,A),
\end{align}
while the boundary condition at $x=A$ is as follows:
\begin{align}\label{eq:bnd-cond-A}
p_{\theta}(A,t|r)
&=
0,\; r\in[0,A),
\end{align}
which in ``PDEs--speak'' is a Dirichlet--type boundary condition.

Using~\eqref{eq:speed+scale-def} it is easy to see that for equation~\eqref{eq:master-eqn} the corresponding scale and speed measures are
\begin{align}\label{eq:scale+speed-measure-answer}
\mathfrak{s}(x)
&=
x^{-2\theta} e^{\tfrac{2}{\mu^2 x}},
\;\;\text{and}\;\;
\mathfrak{m}(x)
=
\frac{2}{\mu^2 x^2}x^{2\theta} e^{-\tfrac{2}{\mu^2 x}},
\end{align}
and, therefore, the boundary condition~\eqref{eq:bnd-cond-00} at zero can be rewritten more explicitly as follows:
\begin{align}\label{eq:bnd-cond-0}
\lim_{x\to0+}\left[x^{2\theta} e^{-\tfrac{2}{\mu^2 x}}\frac{\partial}{\partial x}\left(x^{2-2\theta} e^{\tfrac{2}{\mu^2 x}}p_{\theta}(x,t|r)\right)\right]&=0,\;\; r\in[0,A).
\end{align}

We shall refer to equation~\eqref{eq:master-eqn} complemented by the boundary conditions~\eqref{eq:bnd-cond-A}--\eqref{eq:bnd-cond-0} as the {\em master equation}. It is obtaining the solution to this equation that is the main objective of this work. Hence, the obvious question to be considered next is that of how exactly we intend to undertake this task. To that end, we shall now give a heuristic outline of our approach to solve the master equation~\eqref{eq:master-eqn}. Let us temporarily ``lighten'' the notation $p_{\theta}(x,t|r)$ to $p(x,t)$. The main idea of our solution strategy is to separate the spacial variable, $x$, and the temporal variable, $t$. More concretely, the idea is to seek $p(x,t)$ that is of the form $p(x,t)=\mathfrak{m}(x)\,\psi(x)\,\tau(t)$ where $\mathfrak{m}(x)$ is the speed measure given by~\eqref{eq:scale+speed-measure-answer}, and $\psi(x)$ and $\tau(t)$ are two unknown functions to be determined. If it were possible to ``fit'' $p(x,t)$ of the form $p(x,t)=\mathfrak{m}(x)\,\psi(x)\,\tau(t)$ into the equation~\eqref{eq:master-eqn}, then the substitution $p(x,t)=\mathfrak{m}(x)\,\psi(x)\,\tau(t)$ would bring the master equation~\eqref{eq:master-eqn} into the following form:
\begin{align*}
\dfrac{\tau'(t)}{\tau(t)}
&=
\dfrac{1}{\mathfrak{m}(x)\,\psi(x)}\left(-\dfrac{d}{dx}\big[(1+\theta\mu^2x)\,\mathfrak{m}(x)\,\psi(x)\big]+\dfrac{\mu^2}{2}\dfrac{d^2}{d x^2}\big[x^2\,\mathfrak{m}(x)\,\psi(x)\big]\right),
\end{align*}
and since $x$ and $t$ are now on two different sides of the equation, the only way to ensure the equation holds for {\em all} $x\in[0,A]$, $A>0$, and $t\ge0$ is to require each of the two sides of the equation to be equal to the same {\em constant}, say $\lambda$. Therefore, the substitution $p(x,t)=\mathfrak{m}(x)\,\psi(x)\,\tau(t)$ effectively splits the original PDE~\eqref{eq:master-eqn} into two ODEs:
\begin{align}\label{eq:two-ODEs}
\dfrac{\tau'(t)}{\tau(t)}
&=
\lambda
\;\;\text{and}\;\;
\dfrac{1}{\mathfrak{m}(x)\,\psi(x)}\left(-\dfrac{d}{dx}\big[(1+\theta\mu^2x)\,\mathfrak{m}(x)\,\psi(x)\big]+\dfrac{\mu^2}{2}\dfrac{d^2}{d x^2}\big[x^2\,\mathfrak{m}(x)\,\psi(x)\big]\right)=\lambda,
\end{align}
for some $\lambda$; the set of all $\lambda$'s that make the foregoing two ODEs hold and yet allow to satisfy the initial and boundary conditions will be required.

The first of the two ODEs~\eqref{eq:two-ODEs}, namely the one for $\tau(t)$, is straightforward to solve: the corresponding general nontrivial solution is simply a multiple of the exponential function $e^{\lambda t}$ considered on the interval $t\in[0,+\infty)$; note that $\tau(t)\neq0$ for all $t\in[0,+\infty)$.

To treat the second of the two ODEs~\eqref{eq:two-ODEs}, namely the one for $\psi(x)$, observe first that in view of~\eqref{eq:scale+speed-measure-answer} and~\eqref{eq:operator-G-param} it can be rewritten as
\begin{align}\label{eq:eigen-eqn-two}
\dfrac{\mu^2}{2}\dfrac{d}{d x}\big[x^2\,\mathfrak{m}(x)\,\psi'(x)\big]
&=
\lambda\,\mathfrak{m}(x)\,\psi(x),
\end{align}
and in this new form it can be easily recognized as the {\em characteristic} equation for the linear differential operator:
\begin{align}\label{eq:operator-D-def}
\mathcal{D}
&\triangleq
\dfrac{\mu^2}{2\,\mathfrak{m}(x)}\dfrac{d}{d x}\,x^2\,\mathfrak{m}(x)\,\dfrac{d}{dx},
\end{align}
i.e., equation~\eqref{eq:eigen-eqn-two} determines the eigenvalues $\lambda$ and the corresponding eigenfunctions $\psi(x)$ of the operator $\mathcal{D}$ given by~\eqref{eq:operator-D-def}.

By exactly the same argument it can be shown that the two boundary conditions~\eqref{eq:bnd-cond-A}--\eqref{eq:bnd-cond-0} under the substitution $p(x,t)=\mathfrak{m}(x)\,\psi(x)\,\tau(t)$ convert to
\begin{align}\label{eq:bdd-cond}
\lim_{x\to0+}[x^2\,\mathfrak{m}(x)\,\psi'(x)]
&=
0
\;\;\text{and}\;\;
\psi(A)
=
0,
\end{align}
where $\mathfrak{m}(x)$ is as in~\eqref{eq:scale+speed-measure-answer}; cf., e.g.,~\cite[Formula~(9),~p.~343]{Linetsky:IJTAF2004}. We also note that to get rid of $\tau(t)$ we used the fact that $\tau(t)\neq 0$ for all $t$.

Complemented with the two boundary conditions~\eqref{eq:bdd-cond}, equation~\eqref{eq:eigen-eqn-two} is a Sturm--Liouville problem. Therefore, by attempting to separate the $x$ and $t$ variables we reduced the original equation~\eqref{eq:master-eqn} to the Sturm--Liouville problem~\eqref{eq:eigen-eqn-two} subject to two boundary conditions~\eqref{eq:bdd-cond}. To emphasize the dependence of $\psi(x)$ on $\lambda$ let from now on $\psi(x,\lambda)$ denote the solution (eigenfunction) corresponding to the eigenvalue $\lambda$. If {\em all} of the eigenvalue-eigenfunction pairs $\{\lambda_k,\psi(x,\lambda_k)\}_{k}$ of the operator $\mathcal{D}$ given by~\eqref{eq:operator-D-def} were known, the solution $p(x,t)$ to the master equation~\eqref{eq:master-eqn} would be given by the expansion
\begin{align}\label{eq:pdf-spectral-expansion}
p(x,t|r=y)
&=
\mathfrak{m}(x)\sum_{k} C_{k}(y)\,e^{\lambda_{k} t}\,\psi(x,\lambda_{k}),
\end{align}
where $C_{k}(y)$ and $\lambda_{k}$ are selected so as to make the solution $p(x,t|r=y)$ satisfy the initial temporal condition as well as the two boundary conditions.

With regard to the initial temporal condition, observe that the eigenfunctions corresponding to two different eigenvalues are orthogonal relative to the ``weight function'' $\mathfrak{m}(x)$ given by~\eqref{eq:scale+speed-measure-answer}. Specifically, it holds that
\begin{align}\label{eq:eigenfun-orthonorm}
\int_{0}^{A}\mathfrak{m}(x)\,\psi(x,\lambda_i)\,\psi(x,\lambda_j)\,dx
&=
\indicator{i=j},
\end{align}
where it is assumed that the two eigenfunctions are each of unit ``length'', i.e., $\|\psi(\cdot,\lambda_{i})\|=1=\|\psi(\cdot,\lambda_{j})\|$, with the ``length'' defined as
\begin{align}\label{eq:norm-def}
\|\psi(\cdot,\lambda)\|^2
&\triangleq
\int_{0}^{A}\mathfrak{m}(x)\,\psi^2(x,\lambda)\,dx,
\end{align}
i.e., also relative to the ``weight function'' $\mathfrak{m}(x)$ given by~\eqref{eq:scale+speed-measure-answer}. This standard result from the Sturm--Liouville theory allows to make the expansion~\eqref{eq:pdf-spectral-expansion} more concrete by finding $C_{k}(y)$ explicitly through utilizing the initial temporal condition. Specifically, multiplying~\eqref{eq:pdf-spectral-expansion} through by $\psi(x,\lambda_j)$ and then integrating both sides the result with respect to $x$ over the interval $[0,A)$, we obtain
\begin{align*}
\int_{0}^{A}p(x,t|r=y)\,\psi(x,\lambda_{j})\,dx
&=
\sum_{i}C_{i}(y)\, e^{\lambda_{i} t}\left[\int_{0}^{A}\mathfrak{m}(x)\,\psi(x,\lambda_{i})\,\psi(x,\lambda_{j})\,dx\right],
\end{align*}
whence, in view of the orthogonality property~\eqref{eq:eigenfun-orthonorm}, one can conclude that
\begin{align*}
C_{k}(y)\,e^{\lambda_{k} t}
&=
\int_{0}^{A}p(x,t|r=y)\,\psi(x,\lambda_k)\,dx,
\end{align*}
and because $C_{k}(y)$ is to be independent of $t$, evaluating both sides of the foregoing identity at $t\to0+$ and making use of the initial condition $\lim_{t\to0+}p(x,t|r=y)=\delta(x-y)$, we obtain
\begin{align*}
C_{k}(y)
&=
\int_{0}^{A}\delta(x-y)\,\psi(x,\lambda_{k})\,dx=\psi(y,\lambda_{k}).
\end{align*}

As a result, we can finally conclude from~\eqref{eq:pdf-spectral-expansion} that
\begin{align}\label{eq:pdf-spectral-expansion1}
p(x,t|r=y)
&=
\mathfrak{m}(x)\sum_{k} e^{\lambda_{k} t}\,\psi(x,\lambda_{k})\,\psi(y,\lambda_{k}),
\end{align}
where $\mathfrak{m}(x)$ is as in~\eqref{eq:scale+speed-measure-answer} and $\{\lambda_{k},\psi(x,\lambda_{k})\}_{k}$ are the eigenvalue-eigenfunction pairs of the operator $\mathcal{D}$ defined by~\eqref{eq:operator-D-def}. We would like to reiterate that the obtained expansion~\eqref{eq:pdf-spectral-expansion1} assumes that the eigenfunctions $\psi(x,\lambda_{k})$ are of unit length in the sense of definition~\eqref{eq:norm-def}, i.e., $\|\psi(\cdot,\lambda_{k})\|=1$ for all $\lambda_{k}$. Incidentally, observe the symmetry $p(x,t|r=y)/\mathfrak{m}(x)=p(y,t|r=x)/\mathfrak{m}(y)$, which is known as the detailed balance equation.

The obtained expansion~\eqref{eq:pdf-spectral-expansion1} is at the heart of the entire separation of variables approach (or the Fourier method) that we effectively just outlined. For a more detailed exposition of this approach, see, e.g.,~\cite{Stratonovich:Book1961},~\cite[Chapters~12~\&~26]{Tikhonov+Mironov:Book1977},~\cite[Chapter~5]{Gardiner:Book1985},~\cite{Schuss:Book2010}, and~\cite{Linetsky:IJTAF2004,Linetsky:BookCh2007}. In particular, it is noteworthy that from the general Sturm--Liouville theory it is known that the series in the right-hand side of~\eqref{eq:pdf-spectral-expansion1} is {\em absolutely} convergent for all $t\ge0$ and $x,y\in[0,A]\times[0,A]$. See, e.g.,~\cite{Levitan:Book1950} or~\cite{Levitan+Sargsjan:Book1975}.

We have now set ourselves in a position to follow through with the separation of variables approach summarized above and manifested in formulae~\eqref{eq:scale+speed-measure-answer},~\eqref{eq:eigen-eqn-two}, and~\eqref{eq:pdf-spectral-expansion1}, and attack the master equation~\eqref{eq:master-eqn} directly. This is precisely the object of the next section, which is the main section of this work.

\section{The Main Result}
\label{sec:main-result}

This section is the centerpiece of this work. It is intended to provide a solution to the main problem of this paper: to obtain closed-form formulae for the GSR stopping time's survival functions under the pre- and post-change regimes, i.e., for, respectively, $\Pr_{\infty}(\mathcal{S}_A^r\ge t)$ and $\Pr_{0}(\mathcal{S}_A^r\ge t)$ for all $t\ge0$ and $R_0^r=r\in[0,A]$, with $A>0$ given. Recall that the problem effectively is to solve the master equation~\eqref{eq:master-eqn} subject to two boundary conditions~\eqref{eq:bnd-cond-A}--\eqref{eq:bnd-cond-0}. The solution will yield the densities $p_{\theta}(x,t|r)$, $\theta=\{0,1\}$, defined by~\eqref{eq:p-theta-def}, and these densities can then be used to get the survival functions through~\eqref{eq:p-theta-to-survfun}.

To devise the separation of variables approach outlined in the preceding section and attack the master equation~\eqref{eq:master-eqn} directly, recall that the gist of the Fourier method is to find the eigenvalues $\lambda$ as well as the corresponding eigenfunctions  $\psi(x,\lambda)$ of the operator $\mathcal{D}$. To recover the eigenfunctions, first observe that the change-of-variables $x\mapsto u\triangleq g(x)$ together with the substitution $f(x,t)=h(x)\,v(u,t)$ bring the equation $[\mathcal{G}\circ f](x,t)=\partial f(x,t)/\partial t$ with
\begin{align*}
\mathcal{G}
&\triangleq
\frac{a(x)}{2}\frac{\partial^2}{\partial x^2}+b(x)\frac{\partial}{\partial x}
\end{align*}
to the form
\begin{align}\label{eq:general-eqn-transform}
\begin{split}
\frac{\partial}{\partial t} v(u,t)
&=
\frac{a(x)}{2}\,[g'(x)]^2\,\frac{\partial^2}{\partial u^2} v(u,t)+\\
&\qquad\qquad\qquad\qquad
+\left\{[\mathcal{G}\circ g](x)+a(x)\,g'(x)\,\frac{h'(x)}{h(x)}\right\}\frac{\partial}{\partial u} v(u,t)+\\
&\qquad\qquad\qquad\qquad\qquad\qquad\qquad\qquad\qquad\qquad\qquad
+\frac{[\mathcal{G}\circ h](x)}{h(x)}\,v(u,t).
\end{split}
\end{align}

Next, note that if
\begin{align*}
g(x)
&=
-\int\frac{2}{a(x)}\,dx,
\end{align*}
so that
\begin{align*}
g'(x)
&=
-\frac{2}{a(x)},\;\;
g''(x)
=
\frac{2 a'(x)}{a^2(x)},
\;\;
\text{and}
\;\;
[\mathcal{G}\circ g](x)=\frac{a'(x)-2b(x)}{a(x)},
\end{align*}
then
\begin{align*}
[\mathcal{G}\circ g](x)+a(x)\,g'(x)\,\frac{h'(x)}{h(x)}
&=
\frac{a'(x)-2b(x)}{a(x)}-2\frac{h'(x)}{h(x)},
\end{align*}
whence it is clear that the choice of $h(x)$ such that the equation
\begin{align}\label{eq:h-fun-eqn}
\frac{h'(x)}{h(x)}
&=
\frac{a'(x)-2b(x)}{2a(x)}
\end{align}
is satisfied will cause the term proportional to $v_u(u,t)\triangleq\partial v(u,t)/\partial u$ in the right-hand side of~\eqref{eq:general-eqn-transform} disappear. Moreover, since, by definition~\eqref{eq:speed+scale-def}, the speed measure $\mathfrak{m}(x)$ solves the equation
\begin{align*}
\frac{1}{2}\frac{\partial}{\partial x}\big[a(x)\,\mathfrak{m}(x)\big]-b(x)\,\mathfrak{m}(x)
&=
0,
\;\;\text{so that}\;\;
\frac{\mathfrak{m}'(x)}{\mathfrak{m}(x)}
=
-\frac{a'(x)-2b(x)}{a(x)},
\end{align*}
it is easy to see that equation~\eqref{eq:h-fun-eqn} is solved by $h(x)=1\,/\sqrt{\mathfrak{m}(x)}$. Finally, since by a simple calculation
\begin{align*}
\frac{[\mathcal{G}\circ h](x)}{h(x)}
&=
-
\frac{[a'(x)]^2-2a''(x)\,a(x)-4a'(x)\,b(x)+4b'(x)\,a(x)+4b^2(x)}{8a(x)},
\end{align*}
we have effectively just shown that the change of variables
\begin{align*}
x\mapsto u\triangleq-\int\frac{2}{a(x)}\,dx
\end{align*}
along with the substitution $f(x,t)=v(u,t)\,/\sqrt{\mathfrak{m}(x)}$ convert the equation $f_t(x,t)=[\mathcal{G}\circ f](x,t)$ into the so-called Schr{\"o}dinger form
\begin{align}\label{eq:Schrodinger-form-eqn1}
v_{uu}(u,t)
-
V(u)\, v(u,t)
&=
\frac{a(u)}{2}\, v_t(u,t),
\end{align}
where
\begin{align}\label{eq:Schrodinger-form-eqn2}
V(u)
&\triangleq
\frac{1}{16}\Bigl\{[a'(u)]^2-2a''(u)\,a(u)-4a'(u)\,b(u)+4b'(u)\,a(u)+4b^2(u)\Bigr\},
\end{align}
and let us also point out that any constant (independent of $u$ and $t$) factor that may be present in the substitution $f(x,t)=v(u,t)\,/\sqrt{\mathfrak{m}(x)}$ can be safely dropped without affecting the equation.

All this can be readily applied our equation $[\mathcal{G}\circ\psi](x)=\lambda\,\psi(x)$ on the eigenvalues and eigenfunctions of the operator $\mathcal{D}$. To that end, since in our case $b(x)=1+\theta\mu^2 x$ and $a(x)=\mu^2x^2$, so that
\begin{align}\label{eq:ux-var-def}
x&\mapsto u=u(x)\triangleq-\int\dfrac{2}{a(x)}\,dx=\dfrac{2}{\mu^2 x},\;\text{whence}\;u\mapsto x=x(u)=\dfrac{2}{\mu^2 u}\;\text{and}\;\dfrac{dx}{x}=-\dfrac{du}{u},
\end{align}
and
\begin{align*}
\psi(x)&\mapsto\psi(u)\triangleq\dfrac{v(u)}{\sqrt{\mathfrak{m}(u)}}=\left(\dfrac{\mu^2}{2}\right)^{\theta+\tfrac{1}{2}} u^{\theta-1}\,e^{\tfrac{u}{2}}\,v(u)\propto u^{\theta-1}\,e^{\tfrac{u}{2}}\,v(u),
\end{align*}
then, in view of~\eqref{eq:Schrodinger-form-eqn1}--\eqref{eq:Schrodinger-form-eqn2} and the fact that $\theta^2=\theta$, our equation $[\mathcal{G}\circ\psi](x)=\lambda\,\psi(x)$ becomes
\begin{align}\label{eq:eigenfcn-eqn-whit-form}
v_{uu}(u)+\left\{-\frac{1}{4}+\frac{1-\theta}{u}+\frac{1/4-\xi^2/4}{u^2}\right\}v(u)
&=
0,
\end{align}
where
\begin{align}\label{eq:xi-def}
\xi
&\equiv
\xi(\lambda)
\triangleq
\sqrt{1+\lambda\frac{8}{\mu^2}}
\;\;\text{so that}\;\;
\lambda
\equiv
\lambda(\xi)
=
\frac{\mu^2}{8}(\xi^2-1),
\end{align}
and we note that $\xi$ is, in general, complex-valued. As a matter of fact, as we shall show shortly, the spectrum $\lambda$ of the operator $\mathcal{D}$ given by~\eqref{eq:operator-D-def} is purely real and lies on the nonnegative part of the real line, which, in view of~\eqref{eq:xi-def}, translates to only two possibilities for $\xi\equiv \xi(\lambda)$---to be either purely real (if $\lambda$ is between $-\mu^2/8$ and $0$; note also that in this case $0\le\xi\le 1$) or purely imaginary (if $\lambda$ is below $-\mu^2/8$). This circumstance will become important below, when we get to recovering the spectrum $\lambda$ of the operator $\mathcal{D}$.

\begin{remark}\label{rem:xi-def-alt}
It is noteworthy that equation~\eqref{eq:eigenfcn-eqn-whit-form} is indifferent with respect to the sign of $\xi\triangleq\xi(\lambda)$, i.e., using
\begin{align}\label{eq:xi-def2}
\xi
&\equiv
\xi(\lambda)
\triangleq
-\sqrt{1+\lambda\frac{8}{\mu^2}}
\end{align}
instead of~\eqref{eq:xi-def} does not affect the equation~\eqref{eq:eigenfcn-eqn-whit-form}. As we will see below, this ambiguity in the definition of $\xi$ is ``harmless'' in that it does not alter the solution in any way.
\end{remark}

The obtained equation~\eqref{eq:eigenfcn-eqn-whit-form} is a particular version of the classical Whittaker~\citeyearpar{Whittaker:BAMS1904} equation
\begin{align}\label{eq:Whittaker-eqn}
w_{zz}(z)+\left\{-\dfrac{1}{4}+\dfrac{a}{z}+\dfrac{1/4-b^2}{z^2}\right\}w(z)
&=
0,
\end{align}
where $w(z)$ is the unknown function of $z\in\mathbb{C}$ and $a,b\in\mathbb{C}$ are two given parameters; see, e.g.,~\cite[Chapter~I]{Buchholz:Book1969}. A self-adjoint homogeneous second-order ODE, Whittaker's~\citeyearpar{Whittaker:BAMS1904} equation~\eqref{eq:Whittaker-eqn} is used to define the well-known two Whittaker functions as the equation's two independent (fundamental) solutions. The two Whittaker functions are special functions conventionally denoted as $W_{a,b}(z)$ and $M_{a,b}(z)$, where the indices $a$ and $b$ are the parameters of the equation. Both functions are, in general, complex-valued, even if the two indices---$a$ and $b$---are both purely real. Yet, even when at least one of the two indices---$a$ or $b$---is complex, the Whittaker functions may still be purely real-valued. Since our equation~\eqref{eq:eigenfcn-eqn-whit-form} is a special case of the Whittaker equation~\eqref{eq:Whittaker-eqn}, the eigenfunctions $\psi(x,\lambda)$ of the operator $\mathcal{D}$ are expressible through the Whittaker $W$ and $M$ functions with appropriately chosen indices and argument. In view of this circumstance it makes sense to briefly pause our solution process and summarize certain essential properties of the two Whittaker functions. For a more thorough treatment of the Whittaker equation~\eqref{eq:Whittaker-eqn} and its fundamental solutions $W_{a,b}(z)$ and $M_{a,b}(z)$, see, e.g.,~\cite{Slater:Book1960} or~\cite{Buchholz:Book1969}.

The Whittaker $M_{a,b}(z)$ function is defined only for $2b\neq -1,-2,-3,\ldots$, and, when defined, $M_{a,b}(z)$ is an analytic function for all $a,z\in\mathbb{C}$. Otherwise, if the condition on the second index $b$ is violated, then $M_{a,b}(z)$ experiences a simple pole, but can be regularized through a division by $\Gamma(1+2b)$. Here and onward $\Gamma(z)$ denotes the well-known Gamma function; see, e.g.,~\cite[Chapter~6]{Abramowitz+Stegun:Handbook1964}.

The Whittaker $W_{a,b}(z)$ function is defined through the $M_{a,b}(z)$ function as follows:
\begin{align}\label{eq:Whittaker-fncs-identity}
W_{a,b}(z)
&=
\dfrac{\Gamma(-2b)}{\Gamma(1/2-b-a)}M_{a,b}(z)+\dfrac{\Gamma(2b)}{\Gamma(1/2+b-a)}M_{a,-b}(z);
\end{align}
cf., e.g.,~\cite[Identity~13.1.34,~p.~505]{Abramowitz+Stegun:Handbook1964}. This definition exploits the fact that the Whittaker equation~\eqref{eq:Whittaker-eqn} is even in $b$, so that $M_{a,-b}(z)$ satisfies the Whittaker equation~\eqref{eq:Whittaker-eqn} as well, and, moreover, $M_{a,-b}(z)$ and $M_{a,b}(z)$ are linearly independent. Hence, it is easy to see from~\eqref{eq:Whittaker-fncs-identity} that $W_{a,b}(z)$ is also a solution of the Whittaker equation~\eqref{eq:Whittaker-eqn}. However, unlike $M_{a,-b}(z)$ and $M_{a,b}(z)$, $W_{a,b}(z)$ and $W_{a,-b}(z)$ are not only dependent, they are identical, i.e., $W_{a,b}(z)\equiv W_{a,-b}(z)$, which can be readily deduced from~\eqref{eq:Whittaker-fncs-identity}. This symmetry of the Whittaker $W$ function with respect to the second index $b$ will play an important role in the sequel. With regard to analyticity properties, $W_{a,b}(z)$ is an analytic function of $z$ for all $a,b,z\in\mathbb{C}$. Moreover, as pointed out, e.g., by~\cite{Dikii:Doklady1960}, $W_{a,b}(z)$ is analytic not only as a function of $z\in\mathbb{C}$ but also as a function of $a\in\mathbb{C}$ and as a function of $b\in\mathbb{C}$. This fact will also prove useful below.

Another important and relevant property of the two Whittaker functions is their Wronskian:
\begin{align}\label{eq:Whittaker-funcs-Wronskian}
\mathcal{W}\big\{M_{a,b}(z),W_{a,b}(z)\big\}
&\triangleq
M_{a,b}(z)\,\dfrac{\partial}{\partial z}W_{a,b}(z)-W_{a,b}(z)\,\dfrac{\partial}{\partial z}M_{a,b}(z)
=
-\dfrac{\Gamma(1+2b)}{\Gamma(b-a+1/2)},
\end{align}
cf., e.g.,~\cite[Identity~2.4.27,~p.~26]{Slater:Book1960}. Therefore, $M_{a,b}(z)$ and $W_{a,b}(z)$ are linearly independent whenever $\Gamma(1+2b)/\Gamma(b-a+1/2)\neq0$. In particular, note that if $b-a+1/2=n-1$, $n\in\mathbb{N}$, then the Gamma function in the denominator of the Wronskian~\eqref{eq:Whittaker-funcs-Wronskian} has a simple pole, so that $\mathcal{W}\big\{M_{a,b}(z),W_{a,b}(z)\big\}=0$. As a result, the two Whittaker functions---$M_{a,b}(z)$ and $W_{a,b}(z)$---become linearly dependent. In that case, both degenerate to a type of polynomial known as the Laguerre polynomial; Laguerre polynomials are constructed from the standard monomial basis $\{1,x,x^2,\ldots,x^n,\ldots\}$ by the Gram--Schmidt procedure and form an orthonormal basis on $x\in\mathbb{R}^{+}$ with respect to the measure $e^{-x}dx$.

Going back to the problem, since our equation~\eqref{eq:eigenfcn-eqn-whit-form} is a special case of the Whittaker equation~\eqref{eq:Whittaker-eqn}, and the latter's two fundamental solutions are the Whittaker $M$ and $W$ functions, i.e., $M_{a,b}(z)$ and $W_{a,b}(z)$, it is easy to see that any eigenfunction $\psi(u,\lambda)$ of the operator $\mathcal{D}$ given by~\eqref{eq:operator-D-def} must be of the general form
\begin{align}\label{eq:eigfun-gen-form}
\psi(u,\lambda)
&=
u^{\theta-1}\,e^{\tfrac{u}{2}}\left\{C_1M_{1-\theta,\tfrac{\xi(\lambda)}{2}}(u)+C_2W_{1-\theta,\tfrac{\xi(\lambda)}{2}}(u)\right\},
\end{align}
where $C_1$ and $C_2$ are arbitrary constants. Since these constants affect not only the ``length'' of $\psi(u,\lambda)$, but also whether or not it ``fits'' the boundary conditions~\eqref{eq:bdd-cond}, the obvious question to be considered next is to ``finetune'' $C_1$ and $C_2$ so as to standardize the general eigenfunction $\psi(u,\lambda)$ given by~\eqref{eq:eigfun-gen-form} in accordance with definition~\eqref{eq:norm-def} and make it satisfy both of the boundary conditions~\eqref{eq:bdd-cond}.

Let us first attempt to ``fit'' the general eigenfunction $\psi(u,\lambda)$ given by~\eqref{eq:eigfun-gen-form} into the boundary condition at zero~\eqref{eq:bnd-cond-0}. To that end, since
\begin{align*}
\psi_{x}(x,\lambda)
&=
-\dfrac{\mu^2 u^2}{2}\psi_{u}(u,\lambda),
\end{align*}
it follows that in terms of $u$ given by~\eqref{eq:ux-var-def} the boundary condition at zero~\eqref{eq:bnd-cond-0}, is equivalent to
\begin{align}\label{eq:bdd-cond0-u}
\lim_{u\to+\infty}u^{2-2\theta}\,e^{-u}\,\psi_u(u,\lambda)
&=
0,
\end{align}
and to verify it we are to first find the first derivative of $\psi(u,\lambda)$ with respect to $u$. To find $\psi_u(u,\lambda)$, it is convenient to reexpress $\psi(u,\lambda)$ given by~\eqref{eq:eigfun-gen-form} via two other special functions, viz. the Kummer function usually denoted as $M(a,b,z)$ and the Tricomi function conventionally denoted as $U(a,b,z)$. See, e.g.,~\cite[Chapter~13]{Abramowitz+Stegun:Handbook1964}. These functions form a pair of fundamental solutions to the Kummer equation, a homogeneous second-order ODE which, up to a particular change of variables, is equivalent to the Whittaker equation~\eqref{eq:Whittaker-eqn}. More concretely, this change of variables is as follows
\begin{align*}
M_{a,b}(z)
&=
e^{-\tfrac{z}{2}}\,z^{b+1/2}\,M(1/2+b-a,1+2b,z)
\;\text{and}\;
W_{a,b}(z)
=
e^{-\tfrac{z}{2}}\,z^{b+1/2}\,U(1/2+b-a,1+2b,z);
\end{align*}
cf.,~e.g.,~\cite[Identity~13.1.32,~p.~505]{Abramowitz+Stegun:Handbook1964} and~\cite[Identity~13.1.33,~p.~505]{Abramowitz+Stegun:Handbook1964}, respectively. As a result, we obtain
\begin{align*}
u^{\theta-1}\,e^{\tfrac{u}{2}}\,M_{1-\theta,\tfrac{\xi(\lambda)}{2}}(u)
&=
u^{\xi/2-1/2+\theta}\,M(\xi/2-1/2+\theta,1+\xi,u),\;\text{and}\\
u^{\theta-1}\,e^{\tfrac{u}{2}}\,W_{1-\theta,\tfrac{\xi(\lambda)}{2}}(u)
&=
u^{\xi/2-1/2+\theta}\,U(\xi/2-1/2+\theta,1+\xi,u),
\end{align*}
which, upon substitution back into~\eqref{eq:eigfun-gen-form}, yields
\begin{align*}
\psi(u,\lambda)
&=
u^{\alpha}\big\{C_1M(\alpha,2\alpha+2-2\theta,u)+C_2 U(\alpha,2\alpha+2-2\theta,u)\big\},
\end{align*}
where $\alpha=\alpha(\lambda,\theta)\triangleq\xi(\lambda)/2-1/2+\theta$ so that $1+\xi(\lambda)=2\alpha+2-2\alpha$. In this new form, the eigenfunction $\psi(u,\lambda)$ is simpler to differentiate with respect to $u$. Specifically, we obtain
\begin{align*}
\begin{split}
\psi_u(u,\lambda)
&=
u^{\alpha-1}\Biggl\{C_1\left[\alpha\,M(\alpha,2\alpha+2-2\theta,u)+u\,\frac{\partial}{\partial u}M(\alpha,2\alpha+2-2\theta,u)\right]+\\
&\qquad\qquad\qquad\qquad\qquad
+C_2\left[\alpha\,U(\alpha,2\alpha+2-2\theta,u)+u\,\frac{\partial}{\partial u}U(\alpha,2\alpha+2-2\theta,u)\right]\Biggr\},
\end{split}
\end{align*}
and since
\begin{align*}
a\,M(a,b,z)+z\,\dfrac{\partial}{\partial z}M(a,b,z)
&=
a\,M(a+1,b,z),
\end{align*}
as given by~\cite[Identity~13.4.10,~p.~507]{Abramowitz+Stegun:Handbook1964}, and because
\begin{align*}
a\,U(a,b,z)+z\,\dfrac{\partial}{\partial z}U(a,b,z)
&=
a\,(1+a-b)\,U(a+1,b,z),
\end{align*}
as given by~\cite[Identity~13.4.23,~p.~507]{Abramowitz+Stegun:Handbook1964}, the above expression for $\psi_{u}(u,\lambda)$ reduces further to
\begin{align*}
\psi_u(u,\lambda)
&=
\alpha\,u^{\alpha-1}\,\big\{C_1\,M(\alpha+1,2\alpha+2-2\theta,u)-C_2\,(1+\alpha-2\theta)\,U(\alpha+1,2\alpha+2-2\theta,u)\big\}.
\end{align*}

Next, since $U(a,b,z)=z^{-a}\,\big[1+\mathcal{O}(1/|z|)\big]$, $\Re(z)\to\infty$, as given by~\cite[Formula~13.1.8,~p.~504]{Abramowitz+Stegun:Handbook1964}, and because
\begin{align*}
M(a,b,z)
&=
\dfrac{\Gamma(b)}{\Gamma(a)}\,e^{z}\,z^{a-b}\,\left[1+\mathcal{O}\left(\dfrac{1}{\abs{z}}\right)\right],\;\abs{z}\to+\infty\,,\Re(z)>0,
\end{align*}
as given by~\cite[Formula~13.1.4,~p.~504]{Abramowitz+Stegun:Handbook1964}, then in view of the fact that, by definition, $u$ is not only purely real but also positive, we obtain
\begin{align*}
\begin{split}
u^{2-2\theta}\,e^{-u}\,\psi_{u}(u,\lambda)
&=
\alpha\,u^{2-2\theta}\Biggl\{C_1\,\dfrac{\Gamma(2\alpha+2-2\theta)}{\Gamma(\alpha+1)}\,\big[1+\mathcal{O}(1/u)\big]-\\
&\qquad\qquad\qquad\qquad
-C_2\,(1+\alpha-2\theta)\,e^{-u}\,\big[1+\mathcal{O}(1/u)\big]\Biggr\},
\end{split}
\end{align*}
whence it is apparent that $C_1$ must be taken to be zero in order for $u^{2-2\theta}\,e^{-u}\,\psi_{u}(u,\lambda)$ to tend to 0 as $u$ goes to $+\infty$, i.e., in order to make the eigenfunction $\psi(x,\lambda)$ satisfy the boundary condition at zero given by~\eqref{eq:bdd-cond0-u}.

We are now able to claim that the {\em nonnormalized} eigenfunctions are of the form
\begin{align}\label{eq:eigfun-gen-form2}
\psi(u,\lambda)
&=
C\,e^{\tfrac{u}{2}}\,u^{\theta-1} W_{1-\theta,\tfrac{\xi(\lambda)}{2}}(u),
\end{align}
so that it is clear that the eigenvalues $\lambda$ are determined entirely by the absorbing boundary condition~\eqref{eq:bnd-cond-A}, while the choice of the constant factor $C\neq0$ must be such that $\|\psi(\cdot,\lambda)\|=1$ for each particular eigenvalue $\lambda$. With regard to ensuring that $\|\psi(\cdot,\lambda)\|=1$ for each particular eigenvalue $\lambda$, observe that
\begin{align*}
\|\psi(\cdot,\lambda)\|^2
&\triangleq
\int_{0}^{A}\mathfrak{m}(x)\,\psi^2(x,\lambda)\,dx
=
C^2\left(\dfrac{2}{\mu^2}\right)^{2\theta}\int_{\tfrac{2}{\mu^2 A}}^{+\infty} W_{1-\theta,\tfrac{\xi(\lambda)}{2}}^{2}(u)\,\dfrac{du}{u^2},
\end{align*}
whence it follows that to ``pin down'' $C$ so as to have $\|\psi(\cdot,\lambda)\|=1$ we are to compute the integral
\begin{align}\label{eq:C-imp-int}
\int_{\tfrac{2}{\mu^2 A}}^{+\infty} W_{1-\theta,\tfrac{\xi(\lambda)}{2}}^{2}(u)\,\dfrac{du}{u^2}
\end{align}
for each particular eigenvalue $\lambda$. The foregoing improper integral can be evaluated with the aid of the more general indefinite integral
\begin{align}\label{eq:Whittaker-indef-int}
\int W_{a,b_1}(z)\,W_{a,b_2}(z)\,\dfrac{dz}{z^2}
&=
\dfrac{1}{b_2^2-b_1^2}\left\{W_{a,b_1}(z)\,\dfrac{\partial}{\partial z} W_{a,b_2}(z)-W_{a,b_2}(z)\,\dfrac{\partial}{\partial z} W_{a,b_1}(z)\right\},\;b_1\neq b_2;
\end{align}
cf.,~e.g.,~\cite[Identity~1.13.3.6,~p.~37]{Prudnikov+etal:Book1990}. Specifically, for any $\lambda_{i}\neq\lambda_{j}$, we have
\begin{align}\label{eq:eigfun-dot-prod1}
\begin{split}
\int_{0}^{A}\mathfrak{m}(x)\,&\psi(x,\lambda_{i})\,\psi(x,\lambda_{j})\,dx
\stackrel{\text{(a)}}{=}
C^2\left(\dfrac{2}{\mu^2}\right)^{2\theta}\int_{\tfrac{2}{\mu^2 A}}^{+\infty} W_{1-\theta,\tfrac{\xi(\lambda_{i})}{2}}(u)\,W_{1-\theta,\tfrac{\xi(\lambda_{j})}{2}}(u)\,\dfrac{du}{u^2}\\
&\stackrel{\text{(b)}}{=}
\left(\dfrac{2}{\mu^2}\right)^{2\theta}\dfrac{4C^2}{\xi^2(\lambda_{j})-\xi^2(\lambda_{i})}\Biggl\{W_{1-\theta,\tfrac{\xi(\lambda_{i})}{2}}(u)\,\dfrac{\partial}{\partial u} W_{1-\theta,\tfrac{\xi(\lambda_{j})}{2}}(u)-\\
&\qquad\qquad\qquad\qquad\qquad\qquad
-W_{1-\theta,\tfrac{\xi(\lambda_{j})}{2}}(u)\,\dfrac{\partial}{\partial u}W_{1-\theta,\tfrac{\xi(\lambda_{i})}{2}}(u)\left.\Biggr\}\right|_{u=\tfrac{2}{\mu^2 A}}^{u\to+\infty}\\
&\stackrel{\text{(c)}}{=}
\left(\dfrac{2}{\mu^2}\right)^{2\theta}\dfrac{4C^2}{\xi^2(\lambda_{i})-\xi^2(\lambda_{j})}\Biggl\{W_{1-\theta,\tfrac{\xi(\lambda_{i})}{2}}(u)\,\dfrac{\partial}{\partial u} W_{1-\theta,\tfrac{\xi(\lambda_{j})}{2}}(u)-\\
&\qquad\qquad\qquad\qquad\qquad\qquad
-W_{1-\theta,\tfrac{\xi(\lambda_{j})}{2}}(u)\,\dfrac{\partial}{\partial u}W_{1-\theta,\tfrac{\xi(\lambda_{i})}{2}}(u)\left.\Biggr\}\right|_{u=\tfrac{2}{\mu^2 A}}
\end{split}
\end{align}
where $\text{(a)}$ is due to~\eqref{eq:scale+speed-measure-answer} and~\eqref{eq:eigfun-gen-form2}, the indefinite integral~\eqref{eq:Whittaker-indef-int} is used in $\text{(b)}$ along with the Fundamental Theorem of Calculus, and $\text{(c)}$ is because $1-\theta\in\mathbb{R}$ and
\begin{align}\label{eq:Whit-fnc-asym-uinf}
W_{1-\theta,b}(u)
&=
e^{-\tfrac{u}{2}}\,u^{1-\theta}\,\left[1+\mathcal{O}\left(\frac{1}{u}\right)\right]\;\text{as}\;u\to+\infty,\;\text{for any $b\in\mathbb{C}$},
\end{align}
which is an immediate consequence of the more general asymptotic property of the Whittaker $W$ function
\begin{align*}
W_{a,b}(z)
&=
e^{-\tfrac{z}{2}}\,z^a\,\left[1+\mathcal{O}\left(\frac{1}{z}\right)\right]\;\text{as}\;|z|\to+\infty,\;\text{for any $b\in\mathbb{C}$, provided $|\arg(z)\,|<\pi$},
\end{align*}
established, e.g., in~\cite[Section~16.3]{Whittaker+Watson:Book1927}. Consequently, if $\lambda_{i}\neq\lambda_{j}$, then from the last equality in~\eqref{eq:eigfun-dot-prod1} and the formula~\eqref{eq:xi-def} for $\xi(\lambda)$, the $\mathfrak{m}(x)$--``weighted'' dot-product of $\psi(x,\lambda_{i})$ and $\psi(x,\lambda_{j})$ can be seen to be
\begin{align}\label{eq:eigfun-dot-prod2}
\begin{split}
\int_{0}^{A}\mathfrak{m}(x)\,&\psi(x,\lambda_{i})\,\psi(x,\lambda_{j})\,dx
=
\left(\dfrac{2}{\mu^2}\right)^{2\theta}\dfrac{\mu^2 C^2}{2(\lambda_{i}-\lambda_{j})}\times\\
&\qquad\qquad\times\Biggl\{W_{1-\theta,\tfrac{\xi(\lambda_{i})}{2}}\left(\dfrac{2}{\mu^2 A}\right)\left.\Biggl[\dfrac{\partial}{\partial u} W_{1-\theta,\tfrac{\xi(\lambda_{j})}{2}}(u)\Biggr]\right|_{u=\tfrac{2}{\mu^2 A}}-\\
&\qquad\qquad\qquad\qquad\qquad\qquad
-W_{1-\theta,\tfrac{\xi(\lambda_{j})}{2}}\left(\dfrac{2}{\mu^2 A}\right)\left.\Biggl[\dfrac{\partial}{\partial u} W_{1-\theta,\tfrac{\xi(\lambda_{i})}{2}}(u)\Biggr]\right|_{u=\tfrac{2}{\mu^2 A}}\Biggr\},
\end{split}
\end{align}
and we remark that while the condition $\lambda_{i}\neq\lambda_{j}$ is critical for the validity of~\eqref{eq:eigfun-dot-prod1} and~\eqref{eq:eigfun-dot-prod2}, neither~\eqref{eq:eigfun-dot-prod1} nor~\eqref{eq:eigfun-dot-prod2} actually assumes that $\lambda_{i}$ and $\lambda_{j}$ are eigenvalues of the operator $\mathcal{D}$. That is, both~\eqref{eq:eigfun-dot-prod1} and~\eqref{eq:eigfun-dot-prod2} are valid merely so long as $\lambda_{i}\neq\lambda_{j}$, and regardless of whether $\lambda_{i}$ and $\lambda_{j}$ do belong to the spectrum of $\mathcal{D}$ or not. This is significant for two reasons. On the one hand, if $\lambda_{i}$ and $\lambda_{j}$ are both eigenvalues of $\mathcal{D}$, then $\psi(x,\lambda_{i})$ and $\psi(x,\lambda_{j})$ are both eigenfunctions, and, as such, must satisfy the absorbing boundary condition~\eqref{eq:bnd-cond-A}, in view of which one can immediately conclude from~\eqref{eq:eigfun-dot-prod2} that
\begin{align*}
\int_{0}^{A}\mathfrak{m}(x)\,\psi(x,\lambda_{i})\,\psi(x,\lambda_{j})\,dx
&=
0,\;\text{for $\lambda_{i}\neq\lambda_{j}$},
\end{align*}
which explicitly confirms the validity of the orthogonality property~\eqref{eq:eigenfun-orthonorm} that, as we mentioned earlier, is to hold for the eigenfunctions corresponding to any two different eigenvalues.

On the other hand, the explicit expression~\eqref{eq:eigfun-dot-prod2} we obtained for the $\mathfrak{m}(x)$--``weighted'' dot-product of $\psi(x,\lambda_{i})$ and $\psi(x,\lambda_{j})$ can also be used to bring the eigenfunctions to a unit ``length'', i.e., to have $\|\psi(\cdot,\lambda)\|=1$, or equivalently compute the improper integral~\eqref{eq:C-imp-int}. To that end, the idea is to fix an eigenvalue $\lambda$ and use
\begin{align*}
\|\psi(\cdot,\lambda)\|^2
&=
\lim_{\epsilon\to0}\int_{0}^{A}\mathfrak{m}(x)\,\psi(x,\lambda+\epsilon)\,\psi(x,\lambda)\,dx,
\end{align*}
i.e., effectively pass~\eqref{eq:eigfun-dot-prod2} to the limit as $|\lambda_{i}-\lambda_{j}|\to0$ (assuming, however, that either $\lambda_{i}$ or $\lambda_{j}$ is an eigenvalue of the operator $\mathcal{D}$). Specifically, if $\lambda$ is an eigenvalue, then, as an eigenfunction, $\psi(x,\lambda)$ must satisfy the absorbing boundary condition~\eqref{eq:bnd-cond-A}, so that for any $\epsilon$ such that $\lambda+\epsilon$ is not an eigenvalue and $\psi(x,\lambda+\epsilon)$ is not an eigenfunction, from~\eqref{eq:eigfun-dot-prod2} we obtain
\begin{align}\label{eq:psi-norm-eps}
\begin{split}
\int_{0}^{A}\mathfrak{m}(x)\,&\psi(x,\lambda+\epsilon)\,\psi(x,\lambda)\,dx
=\\
&\qquad\qquad
=
\left(\dfrac{2}{\mu^2}\right)^{2\theta}\dfrac{\mu^2 C^2}{2\epsilon}\,W_{1-\theta,\tfrac{\xi(\lambda+\epsilon)}{2}}\left(\dfrac{2}{\mu^2 A}\right)\left.\Biggl[\dfrac{\partial}{\partial u} W_{1-\theta,\tfrac{\xi(\lambda)}{2}}(u)\Biggr]\right|_{u=\tfrac{2}{\mu^2 A}},
\end{split}
\end{align}
and before we proceed to taking the limit as $\epsilon\to0$ it is worth recalling the aforementioned observation made by~\cite{Dikii:Doklady1960} that $W_{a,b}(z)$ is analytic not only as a function of $z\in\mathbb{C}$ but also as a function of $b\in\mathbb{C}$. As a result, we have the first-order Taylor expansion
\begin{align*}
W_{1-\theta,\tfrac{\xi(\lambda+\epsilon)}{2}}\left(\dfrac{2}{\mu^2 A}\right)
&=
W_{1-\theta,\tfrac{\xi(\lambda)}{2}}\left(\dfrac{2}{\mu^2 A}\right)
+
\dfrac{\epsilon}{2}\left.\left\{\Biggl[\dfrac{\partial}{\partial y}\xi(y)\Biggr]\left.\Biggl[\dfrac{\partial}{\partial b}W_{1-\theta,b}\left(\dfrac{2}{\mu^2 A}\right)\Biggr]\right|_{b=\tfrac{\xi(y)}{2}}\right\}\right|_{y=\lambda_{\epsilon}^{*}},
\end{align*}
where $\lambda_{\epsilon}^{*}$ is within an $|\epsilon|>0$ distance from $\lambda$, i.e., $\lambda_{\epsilon}^{*}\to\lambda$ as $\epsilon\to0$. Since $\lambda$ is an eigenvalue, the absorbing boundary condition~\eqref{eq:bnd-cond-A} enables us to simplify the above Taylor expansion to
\begin{align*}
W_{1-\theta,\tfrac{\xi(\lambda+\epsilon)}{2}}\left(\dfrac{2}{\mu^2 A}\right)
&=
\dfrac{\epsilon}{2}\left.\left\{\Biggl[\dfrac{\partial}{\partial y}\xi(y)\Biggr]\left.\Biggl[\dfrac{\partial}{\partial b}W_{1-\theta,b}\left(\dfrac{2}{\mu^2 A}\right)\Biggr]\right|_{b=\tfrac{\xi(y)}{2}}\right\}\right|_{y=\lambda_{\epsilon}^{*}},
\end{align*}
or
\begin{align}\label{eq:Whit-fnc-Taylor}
W_{1-\theta,\tfrac{\xi(\lambda+\epsilon)}{2}}\left(\dfrac{2}{\mu^2 A}\right)
&=
\dfrac{2\epsilon}{\mu^2\xi(\lambda_{\epsilon}^{*})}\left.\Biggl[\dfrac{\partial}{\partial b}W_{1-\theta,b}\left(\dfrac{2}{\mu^2 A}\right)\Biggr]\right|_{b=\tfrac{\xi(\lambda_{\epsilon}^{*})}{2}},
\end{align}
because
\begin{align*}
\dfrac{\partial}{\partial\lambda}\xi(\lambda)
&=
\dfrac{4}{\mu^2\xi(\lambda)},
\end{align*}
as can obtained at once from~\eqref{eq:xi-def}. Plugging~\eqref{eq:Whit-fnc-Taylor} back over into~\eqref{eq:psi-norm-eps} yields
\begin{align*}
\begin{split}
\int_{0}^{A}\mathfrak{m}(x)\,&\psi(x,\lambda+\epsilon)\,\psi(x,\lambda)\,dx
=\\
&\qquad\qquad
=
\left(\dfrac{2}{\mu^2}\right)^{2\theta}\dfrac{C^2}{\xi(\lambda_{\epsilon}^{*})}\left.\Biggl[\dfrac{\partial}{\partial b}W_{1-\theta,b}\left(\dfrac{2}{\mu^2 A}\right)\Biggr]\right|_{b=\tfrac{\xi(\lambda_{\epsilon}^{*})}{2}}\left.\Biggl[\dfrac{\partial}{\partial u} W_{1-\theta,\tfrac{\xi(\lambda)}{2}}(u)\Biggr]\right|_{u=\tfrac{2}{\mu^2 A}},
\end{split}
\end{align*}
whence the trivial observation that $\xi(\lambda+\epsilon)\to\xi(\lambda)$ as $\epsilon\to0$ combined with the aforementioned continuity of $W_{a,b}(z)$ as a function of $b$ lead further to
\begin{align*}
\|\psi(\cdot,\lambda)\|^2
&=
\left(\dfrac{2}{\mu^2}\right)^{2\theta}\dfrac{C^2}{\xi(\lambda)}\left.\Biggl[\dfrac{\partial}{\partial b}W_{1-\theta,b}\left(\dfrac{2}{\mu^2 A}\right)\Biggr]\right|_{b=\tfrac{\xi(\lambda)}{2}}\left.\Biggl[\dfrac{\partial}{\partial u} W_{1-\theta,\tfrac{\xi(\lambda)}{2}}(u)\Biggr]\right|_{u=\tfrac{2}{\mu^2 A}},
\end{align*}
so that finally it is apparent that the choice
\begin{align}\label{eq:const-C-sqr-answer}
C^2
&\equiv
C_{\lambda,\theta,A}^2
=
\left(\dfrac{\mu^2}{2}\right)^{2\theta}\xi(\lambda)\left/\Biggl\{\left.\Biggl[\dfrac{\partial}{\partial b}W_{1-\theta,b}\left(\dfrac{2}{\mu^2 A}\right)\Biggr]\right|_{b=\tfrac{\xi(\lambda)}{2}}\left.\Biggl[\dfrac{\partial}{\partial u} W_{1-\theta,\tfrac{\xi(\lambda)}{2}}(u)\Biggr]\right|_{u=\tfrac{2}{\mu^2 A}}\Biggr\}\right.
\end{align}
guarantees that $\|\psi(\cdot,\lambda)\|^2=1$ holds for each particular eigenvalue $\lambda$. The obtained result is in agreement with~\cite[Proposition~1]{Linetsky:OR2004} which, in turn, was established using a different technique, viz. one proposed in \cite[Section~5.1]{Linetsky:IJTAF2004}.

It remains to find the actual eigenvalues $\{\lambda\}$ of the operator $\mathcal{D}$. As the first step toward recovering the spectrum $\{\lambda\}$ of the operator $\mathcal{D}$ given by~\eqref{eq:operator-D-def}, let us demonstrate that, under the boundary conditions~\eqref{eq:bnd-cond-A}--\eqref{eq:bnd-cond-0}, the spectrum cannot lie to the left of the origin, i.e., it is impossible to have $\lambda>0$. Indeed, by multiplying~\eqref{eq:eigen-eqn-two} through by $\psi(x,\lambda)$ and then integrating both sides of the result with respect to $x$ over the interval $[0,A)$, we obtain
\begin{align*}
\dfrac{\mu^2}{2}\int_{0}^{A}\psi(x,\lambda)\,\dfrac{d}{dx}\left[x^2\,\mathfrak{m}(x)\,\dfrac{d}{dx}\psi(x,\lambda)\right]dx
-
\lambda\int_{0}^{A}\mathfrak{m}(x)\,\psi^2(x,\lambda)\,dx
&=
0,
\end{align*}
which, after recognizing the second term in the left-hand side as $\|\psi(\cdot,\lambda)\|^2$, i.e., the squared norm~\eqref{eq:norm-def} of $\psi(x,\lambda)$, reduces further to
\begin{align*}
\lambda\,\|\psi(\cdot,\lambda)\|^2
&=
\dfrac{\mu^2}{2}\int_{0}^{A}\psi(x,\lambda)\,\dfrac{d}{dx}\left[x^2\,\mathfrak{m}(x)\,\dfrac{d}{dx}\psi(x,\lambda)\right]dx,
\end{align*}
or
\begin{align}\label{eq:lambda-step2}
\lambda
&=
\dfrac{\mu^2}{2}\int_{0}^{A}\psi(x,\lambda)\,\dfrac{d}{dx}\left[x^2\,\mathfrak{m}(x)\,\dfrac{d}{dx}\psi(x,\lambda)\right]dx,
\end{align}
because without loss of generality $\psi(x,\lambda)$ may be assumed to be of unit length in the sense of~\eqref{eq:norm-def}, i.e., $\|\psi(\cdot,\lambda)\|^2=1$. Next, integration by parts applied to the integral in the right-hand side of~\eqref{eq:lambda-step2} reduces the latter to
\begin{align*}
\begin{split}
\lambda
&=
\dfrac{\mu^2}{2}\left\{\left.\psi(x,\lambda)\left[x^2\,\mathfrak{m}(x)\,\dfrac{d}{dx}\psi(x,\lambda)\right]\right|_{x\to0+}^{x=A}-\int_{0}^{A}x^2\,\mathfrak{m}(x)\left[\dfrac{d}{dx}\psi(x,\lambda)\right]^2dx\right\}\\
&=
-\dfrac{\mu^2}{2}\int_0^{A}x^2\,\mathfrak{m}(x)\left[\dfrac{d}{dx}\psi(x,\lambda)\right]^2dx,
\end{split}
\end{align*}
where we also used the boundary conditions~\eqref{eq:bnd-cond-A}--\eqref{eq:bnd-cond-0} but in the form~\eqref{eq:bdd-cond}. The obtained result implies that $\lambda\le0$, i.e., the spectrum must be concentrated in the nonpositive half of the real line. Consequently, $\xi(\lambda)$ given by~\eqref{eq:xi-def} is either purely real or purely imaginary.

With regard to actually finding the eigenvalues $\lambda$, in view of the remark we made earlier that the eigenvalues $\lambda$ are determined entirely by the absorbing boundary condition~\eqref{eq:bnd-cond-A}, the problem is effectively to solve the equation $\psi(A,\lambda)=0$ where the unknown is $\lambda\le0$, and $0<A<+\infty$ is given. Written explicitly, the equation to be solved to recover the spectrum of the operator $\mathcal{D}$ is
\begin{align*}
e^{\tfrac{1}{\mu^2 A}}\,\left(\dfrac{2}{\mu^2 A}\right)^{\theta-1} W_{1-\theta,\tfrac{\xi(\lambda)}{2}}\left(\dfrac{2}{\mu^2 A}\right)
&=
0,
\end{align*}
which is equivalent to
\begin{align}\label{eq:eigval-eqn2}
W_{1-\theta,\tfrac{\xi(\lambda)}{2}}\left(\dfrac{2}{\mu^2 A}\right)
&=
0,
\end{align}
and it is worth recalling again that $\xi(\lambda)$ is as in~\eqref{eq:xi-def}. For a fixed $0<A<+\infty$, the solutions, $\lambda$, as well as the total number, $N$, thereof depend the two indices $1-\theta$ and $\xi(\lambda)/2$ of the Whittaker $W$ function present in the left-hand side of~\eqref{eq:eigval-eqn2}. With regard to the number of solutions $N$, one of the key factors that determines $N$ is whether or not the two indices $1-\theta$ and $\xi(\lambda)/2$ of the Whittaker $W$ function are purely real or purely imaginary. Since in our case $1-\theta$ is a real number (in fact, it can take only two values: either 0 or 1) and $\xi(\lambda)$, as we argued above, is either purely real or purely imaginary, there are two cases to consider.

The easiest of the two cases is when $\xi(\lambda)$ is purely real. Since, according to~\eqref{eq:xi-def} this occurs only when $-\mu^2/8\le\lambda$, and because we also have the restriction that $\lambda\le0$, it follows that, if $\xi(\lambda)$ is to be purely real, it has to range between $0$ and $1$. In this case, equation~\eqref{eq:eigval-eqn2} can be handled by appealing, e.g., to~\cite[Theorem~4,~p.~944]{Dikii:Doklady1960}, according to which the number $N$ of {\em real} solutions $z$ to the equation $W_{a,b}(z)=0$ when both indices $a$ and $b$ of the Whittaker $W$ function are purely real and such that $a\ge0$ and $\abs{b}\le 1/2$ is $N=\max\{-[\abs{b}-a+1/2],0\}$, where $[x]$ stands for the largest integer not exceeding $x$. See also, e.g.,~\cite[Theorem~9,~p.~11]{Tsvetkoff:Doklady1941},~\cite{Tricomi:MZ1950}, and~\cite[Theorem~2,~p.~156]{Dyson:PhF1960}. Moreover, under the stated assumptions on the two indices of the Whittker $W$ function, the solutions $z$ of the equation $W_{a,b}(z)=0$, should they exist, must be not only real, but also positive. With this mind, let us now turn our equation~\eqref{eq:eigval-eqn2} around and assume instead that $\xi(\lambda)$ is fixed and that the equation is actually for $A>0$. Then, from the aforementioned~\cite[Theorem~4,~p.~944]{Dikii:Doklady1960} and the observation that $\xi(\lambda)/2$ has to be between 0 and $1/2$, it is easy to see that equation~\eqref{eq:eigval-eqn2} viewed as an equation for $A$ such that $A>0$ is actually inconsistent when $\theta=1$. Put another way, if $\theta=1$, then no choice of $\lambda$ such that $-\mu^2/8\le\lambda\le 0$ (so that $\xi(\lambda)$ is purely real and between 0 and 1) can possibly make equation~\eqref{eq:eigval-eqn2} with $A$ being the unknown (restricted to the positive real semiaxis) have even a single solution. This means that in the post-drift regime the spectrum of the operator $\mathcal{D}$ lies entirely to the left of the point $-\mu^2/8$. However, if $\theta=0$, then for any $\lambda$ lying inside the interval $[-\mu^2/8,0]$ equation~\eqref{eq:eigval-eqn2} with $A$ being the unknown (restricted to the positive real semiaxis) may have a solution but no more than one. Flipping this back around, this means that, if $\theta=0$, then for any given $A>0$, there is at most one $\lambda$ located inside the interval $[-\mu^2/8,0]$ (so that $\xi(\lambda)$ is purely real and between 0 and 1) for which equation~\eqref{eq:eigval-eqn2} is satisfied. We therefore arrive at the conclusion that in the pre-drift regime, the spectrum of the operator $\mathcal{D}$ may have at most a single point $\lambda$ lying inside the interval $[-\mu^2/8,0]$. Specifically, if we let $\alpha_{0,A}$ to denote the solution (should it exist) of the equation
\begin{align}\label{eq:eigval-eqn-real}
W_{1,\tfrac{\alpha_{0,A}}{2}}\left(\dfrac{2}{\mu^2 A}\right)
&=
0,
\end{align}
then from~\eqref{eq:xi-def} the corresponding value of $\lambda$ can be seen to be $\lambda=\mu^2(\alpha_{0,A}^2-1)/8$, and for the reasons explained above $0\le \alpha_{0,A}\le 1$, so that $-\mu^2/8\le\lambda\le 0$. Once again, the need to solve equation~\eqref{eq:eigval-eqn-real} arises only in the pre-drift regime, i.e., when $\theta=0$ (or $\nu=\infty$), and should equation~\eqref{eq:eigval-eqn-real} have a solution, it has to be the only solution.

The situation is drastically different when $\xi(\lambda)$ is purely imaginary, which happens when $\lambda\le-\mu^2/8$. In this case, it is convenient to set $\xi(\lambda)=\imath\beta(\lambda)$ where $\beta(\lambda)\in\mathbb{R}$; here and onward $\imath$ denotes the imaginary unit, i.e., $\imath\triangleq\sqrt{-1}$. Moreover, since the Whittaker $W$ is symmetric with respect to the second index, i.e., $W_{a,b}(z)=W_{a,-b}(z)$ for all $a,b,z\in\mathbb{C}$, it is sufficient to assume that $\beta(\lambda)\ge0$. Going back to Remark~\ref{rem:xi-def-alt}, it is due to this symmetry of the Whittaker $W$ function with respect to the second index that the ambiguity in choosing $\xi(\lambda)$ as in~\eqref{eq:xi-def} or as in~\eqref{eq:xi-def2} is nothing to worry about, as it does not cause the solution to change. Moreover, recall the definition~\eqref{eq:Whittaker-fncs-identity} of the Whittaker $W_{a,b}(z)$ function, and note that when $a=1-\theta\in\mathbb{R}$ and $b=\imath\beta/2$, $\beta\in\mathbb{R}$, it takes the form
\begin{align*}
W_{1-\theta,\tfrac{\imath\beta}{2}}(z)
&=
\dfrac{\Gamma(-\imath\beta)}{\Gamma(-\imath\beta/2-1/2+\theta)}\,M_{1-\theta,\tfrac{\imath\beta}{2}}(z)+\dfrac{\Gamma(\imath\beta)}{\Gamma(\imath\beta/2-1/2+\theta)}\,M_{1-\theta,-\tfrac{\imath\beta}{2}}(z),
\end{align*}
whence, because the two terms in the right-hand side are complex conjugates of each other, one may deduce that $W_{1-\theta,\tfrac{\imath\beta}{2}}(z)$ is necessarily real-valued. More specifically,
\begin{align*}
W_{1-\theta,\tfrac{\imath\beta}{2}}(z)
&=
2\Re\left\{\dfrac{\Gamma(-\imath\beta)}{\Gamma(-\imath\beta/2-1/2+\theta)}\,M_{1-\theta,\tfrac{\imath\beta}{2}}(z)\right\},
\end{align*}
where here and onward $\Re(z)$ denotes the real part of a complex number $z\in\mathbb{C}$. More explicitly, the foregoing identity can be written as follows:
\begin{align}\label{eq:Whittaker-fncs-identity-v3}
\begin{split}
W_{1-\theta,\tfrac{\imath\beta}{2}}(z)
&=
2\abs{\dfrac{\Gamma(-\imath\beta)}{\Gamma(-\imath\beta/2-1/2+\theta)}\,M_{1-\theta,\tfrac{\imath\beta}{2}}(z)}\times\\
&\qquad\qquad\qquad\times\cos\left\{\arg\Gamma(-\imath\beta)-\arg\Gamma(-\imath\beta/2-1/2+\theta)+\arg M_{1-\theta,\tfrac{\imath\beta}{2}}(z)\right\},
\end{split}
\end{align}
where $\arg{z}$ means the complex phase (angle between the real and imaginary components) of a complex number $z\in\mathbb{C}$. Formula~\eqref{eq:Whittaker-fncs-identity-v3} is another, more important consequence of the symmetry of the Whittaker $W$ function with respect to the second index. Specifically, it is now clear that, contrary to the case when $-\mu^2/8\le\lambda\le0$ so that $\xi(\lambda)$ is purely real and between 0 and 1, in the case when $\lambda\le-\mu^2/8$ so that $\xi(\lambda)$ is purely imaginary, the number of solutions to the equation~\eqref{eq:eigval-eqn2} is countably many, whether $\theta=0$ or $\theta=1$; cf.~\cite[Theorem~3,~p.~156]{Dyson:PhF1960} and~\cite[Theorem~5,~p.~157]{Dyson:PhF1960}. In fact, a comment made by~\cite[p.~950]{Dikii:Doklady1960} that because the structure of the Whittaker equation~\eqref{eq:Whittaker-eqn} is such that
\begin{align*}
W_{a,b}(z_0)
&=
0
\;\text{necessarily implies that}\;
\left.\left[\frac{\partial^2}{\partial z^2} W_{a,b}(z)\right]\right|_{z=z_0}
=
0,
\end{align*}
it follows from the theory of implicit functions that
\begin{align*}
\left.\left[\frac{\partial}{\partial z} W_{a,b}(z)\right]\right|_{z=z_0}
&\neq
0,
\end{align*}
combined together with the Wronskian~\eqref{eq:Whittaker-funcs-Wronskian} lead to the conclusion that $M_{a,b}(z_0)\neq 0$ if $W_{a,b}(z_0)=0$. Therefore, setting the right-hand side of~\eqref{eq:Whittaker-fncs-identity-v3} equal to zero is equivalent to requiring the argument of the cosine function in the right-hand side of~\eqref{eq:Whittaker-fncs-identity-v3} to be $\pi/2+\pi k$, $k\in\mathbb{Z}$. This ultimately translates to the number of eigenvalues of the operator $\mathcal{D}$ that lie to the left of the point $-\mu^2/8$ being countably many, no matter whether $\theta$ is $0$ or $1$. Moreover, all these eigenvalues are simple (i.e., of algebraic multiplicity one), which is in agreement with the general Sturm--Liouville theory; cf., e.g.,~\cite{Levitan:Book1950} or~\cite{Levitan+Sargsjan:Book1975}.

We are now in a position to put all of the above together and write down the sought-after density, $p_{\theta}(x,t|r)$, in a closed form. Specifically, we obtain:
\begin{align}\label{eq:p_theta-formula}
\begin{split}
p_{\theta}(x,t|r=y)
&=
\dfrac{\mu^2}{2}\,e^{\tfrac{1}{\mu^2 y}-\tfrac{1}{\mu^2 x}}e^{-\tfrac{\mu^2 t}{8}}\left(\dfrac{y}{x}\right)^{1-\theta}\times\\
&\qquad
\times\Biggl\{\,(1-\theta)\,e^{\tfrac{\mu^2 t}{8}\alpha_{0,A}^2}\,\tilde{C}_{0,0,A}^2\,W_{1,\tfrac{\alpha_{0,A}}{2}}\left(\dfrac{2}{\mu^2 x}\right)W_{1,\tfrac{\alpha_{0,A}}{2}}\left(\dfrac{2}{\mu^2 y}\right)+\\
&\qquad\qquad
+
\sum_{n=1}^{\infty} e^{-\tfrac{\mu^2 t}{8}\beta_{n,\theta,A}^2}\,\tilde{C}_{n,\theta,A}^2\,W_{1-\theta,\tfrac{\imath\beta_{n,\theta,A}}{2}}\left(\dfrac{2}{\mu^2 x}\right)W_{1-\theta,\tfrac{\imath\beta_{n,\theta,A}}{2}}\left(\dfrac{2}{\mu^2 y}\right)\Biggr\},
\end{split}
\end{align}
where $x,y\in[0,A]$ and $t\ge0$, and recall that $\theta$ is either 0 ($\nu=\infty$) or 1 ($\nu=0$), the detection threshold $A>0$ is given, the constant $\alpha_{0,A}\in[0,1]$ is the only zero (should it exist) of the equation
\begin{align}\label{eq:alpha-eqn}
W_{1,\tfrac{\alpha_{0,A}}{2}}\left(\dfrac{2}{\mu^2 A}\right)
&=
0,
\end{align}
which is nothing but equation~\eqref{eq:eigval-eqn-real}, the constant $\tilde{C}_{0,0,A}^2$ is
\begin{align}\label{eq:Cinf-eqn}
\tilde{C}_{0,0,A}^2
&=
\alpha_{0,A}\left/\Biggl\{\left.\Biggl[\dfrac{\partial}{\partial b}W_{1,b}\left(\dfrac{2}{\mu^2 A}\right)\Biggr]\right|_{b=\tfrac{\alpha_{0,A}}{2}}\left.\Biggl[\dfrac{\partial}{\partial u} W_{1,\tfrac{\alpha_{0,A}}{2}}(u)\Biggr]\right|_{u=\tfrac{2}{\mu^2 A}}\Biggr\}\right.,
\end{align}
which comes from~\eqref{eq:const-C-sqr-answer} evaluated at $\lambda_{0,0,A}$ such that $\xi(\lambda_{0,0,A})=\alpha_{0,A}\in[0,1]$, and finally
the series $\{\beta_{n,\theta,A}\}_{n\ge1}$ is formed of the (countably many) solutions $\beta_{\theta,A}\ge0$ of the equation
\begin{align}\label{eq:beta-eqn}
W_{1-\theta,\tfrac{\imath\beta_{\theta,A}}{2}}\left(\dfrac{2}{\mu^2 A}\right)
&=
0,
\end{align}
and
\begin{align}\label{eq:Czero-eqn}
\tilde{C}_{n,\theta,A}^2
&=
\imath\beta_{n,\theta,A}\left/\Biggl\{\left.\Biggl[\dfrac{\partial}{\partial b}W_{1-\theta,b}\left(\dfrac{2}{\mu^2 A}\right)\Biggr]\right|_{b=\tfrac{\imath\beta_{n,\theta,A}}{2}}\left.\Biggl[\dfrac{\partial}{\partial u} W_{1-\theta,\tfrac{\imath\beta_{n,\theta,A}}{2}}(u)\Biggr]\right|_{u=\tfrac{2}{\mu^2 A}}\Biggr\}\right.,
\end{align}
which again comes from~\eqref{eq:const-C-sqr-answer} evaluated at $\lambda_{n,\theta,A}$ such that $\xi(\lambda_{n,\theta,A})=\imath\beta_{n,\theta,A}$. We note that because $\theta$ is either 0 or 1, the first term inside the braces in the right-hand side of~\eqref{eq:p_theta-formula} appears only when $\theta=0$, i.e., in the pre-change regime: only in this regime do we have to find $\alpha_{0,A}\in[0,1]$ from equation~\eqref{eq:alpha-eqn}, and then, should equation~\eqref{eq:alpha-eqn} have a solution, compute constant $\tilde{C}_{0,0,A}$ from~\eqref{eq:Cinf-eqn}. Otherwise, in the post-change regime, i.e., when $\theta=1$, the first term inside the braces in the right-hand side of~\eqref{eq:p_theta-formula} is zero (viz. need not be evaluated altogether) because of the factor of $1-\theta=0$ present in front of it. It also important to repeat the comment we made at the end of Section~\ref{sec:problem+preliminaries} that the expansion~\eqref{eq:pdf-spectral-expansion1} is {\em absolutely} convergent for all $t\ge0$ and $x,y\in[0,A]\times[0,A]$. Therefore, the series in the right-hand side of~\eqref{eq:p_theta-formula} is also absolutely convergent for all $t\ge0$ and $x,y\in[0,A]\times[0,A]$---whether $\theta$ is $0$ or $1$.

The survival functions $\Pr_{\infty}(\mathcal{S}_A^r\ge t)$ and $\Pr_{0}(\mathcal{S}_A^r\ge t)$ can be obtained from~\eqref{eq:p_theta-formula} through~\eqref{eq:p-theta-to-survfun}. That is, to get the two survival functions and thus achieve the main objective of this work, the whole problem now is to merely integrate the right-hand side of~\eqref{eq:p_theta-formula} with respect to $x$ over the interval $[0,A]$. This integration can be carried out with the aid of~\cite[Formula~7.623.7,~p.~824]{Gradshteyn+Ryzhik:Book2007} which states that
\begin{align}\label{eq:Whit-int-formula}
\int_{1}^{+\infty} (x-1)^{c-1} x^{a-c-1}\, e^{-\tfrac{qx}{2}}\,W_{a,b}(qx)\,dx
&=
\Gamma(c)\,e^{-\tfrac{q}{2}}\,W_{a-c,b}(q),
\end{align}
provided that $\Re(c)>0$ and $\Re(q)>0$. Specifically, using the foregoing integral identity, we obtain
\begin{align*}
\int_{0}^{A}\mathfrak{m}(x)\,\psi(x,\lambda)\,dx
&\stackrel{\text{(a)}}{=}
C\int_{0}^{A}\left(\dfrac{2 x}{\mu^2}\right)^{\theta}e^{-\tfrac{1}{\mu^2 x}}\,W_{1-\theta,\tfrac{\xi(\lambda)}{2}}\left(\dfrac{2}{\mu^2 x}\right)\dfrac{dx}{x}\\
&\stackrel{\text{(b)}}{=}
C\left(\dfrac{2A}{\mu^2}\right)^{\theta}\int_{1}^{+\infty}y^{-(1+\theta)}\,e^{-\tfrac{y}{\mu^2 A}}\,W_{1-\theta,\tfrac{\xi(\lambda)}{2}}\left(\dfrac{2y}{\mu^2 A}\right)dy\\
&\stackrel{\text{(c)}}{=}
C\left(\dfrac{2A}{\mu^2}\right)^{\theta}e^{-\tfrac{1}{\mu^2 A}}\,W_{-\theta,\tfrac{\xi(\lambda)}{2}}\left(\dfrac{2}{\mu^2 A}\right),
\end{align*}
where $\text{(a)}$ makes use of the expression~\eqref{eq:scale+speed-measure-answer} for the speed measure $\mathfrak{m}(x)$, the expression~\eqref{eq:eigfun-gen-form2} for the eigenfunction $\psi(x,\lambda)$, and~\eqref{eq:ux-var-def}, $\text{(b)}$ is due to the change of variables $x\mapsto y\triangleq y(x)=A/x$, and $\text{(c)}$ is identity~\eqref{eq:Whit-int-formula} with $c=1$, $a=1-\theta$, and $q=2/(\mu^2 A)$ (note that for this choice of $c$ and $q$ the conditions $\Re(c)>0$ and $\Re(q)>0$ required for the integral to hold are fulfilled). Therefore, we obtain:
\begin{align}\label{eq:survfun-inf-formula1}
\begin{split}
\Pr_{\infty}(\mathcal{S}_{A}^{r=y}\ge t)
&=
\dfrac{\mu^2 y}{2}\,e^{\tfrac{1}{\mu^2 y}-\tfrac{1}{\mu^2 A}}e^{-\tfrac{\mu^2 t}{8}}\times\\
&\qquad
\times\Biggl\{e^{\tfrac{\mu^2 t}{8}\alpha_{0,A}^2}\,\tilde{C}_{0,0,A}^2\,W_{0,\tfrac{\alpha_{0,A}}{2}}\left(\dfrac{2}{\mu^2 A}\right)W_{1,\tfrac{\alpha_{0,A}}{2}}\left(\dfrac{2}{\mu^2 y}\right)+\\
&\qquad\qquad
+
\sum_{n=1}^{\infty}e^{-\tfrac{\mu^2 t}{8}\beta_{n,\theta,A}^2}
\tilde{C}_{n,0,A}^2\,W_{0,\tfrac{\imath\beta_{n,0,A}}{2}}\left(\dfrac{2}{\mu^2 A}\right)W_{1,\tfrac{\imath\beta_{n,0,A}}{2}}\left(\dfrac{2}{\mu^2 y}\right)\Biggr\},
\end{split}
\end{align}
and
\begin{align}\label{eq:survfun-zero-formula1}
\begin{split}
\Pr_{0}(\mathcal{S}_{A}^{r=y}\ge t)
&=
\dfrac{2 A}{\mu^2}\,e^{\tfrac{1}{\mu^2 y}-\tfrac{1}{\mu^2 A}}e^{-\tfrac{\mu^2 t}{8}}\times\\
&\qquad\qquad
\times
\sum_{n=1}^{\infty}e^{-\tfrac{\mu^2 t}{8}\beta_{n,1,A}^2}\,\tilde{C}_{n,1,A}^2\,W_{-1,\tfrac{\imath\beta_{n,1,A}}{2}}\left(\dfrac{2}{\mu^2 A}\right)W_{0,\tfrac{\imath\beta_{n,1,A}}{2}}\left(\dfrac{2}{\mu^2 y}\right),
\end{split}
\end{align}
where, as before, $\alpha_{0,A}\in[0,1]$ is the (at most one) root of equation~\eqref{eq:alpha-eqn}, the series $\beta_{n,\theta,A}$ with $\theta$ either $0$ or $1$ are formed of the countably many solutions $\beta_{\theta,A}\ge0$ of equation~\eqref{eq:beta-eqn}, and constants $\tilde{C}_{n,\theta,A}^2$ for $n\ge0$ and $\theta=\{0,1\}$ are as in~\eqref{eq:Cinf-eqn} and~\eqref{eq:Czero-eqn}, respectively. That said, unlike the series in the right-hand side of~\eqref{eq:p_theta-formula}, the series in the right-hand side of~\eqref{eq:survfun-inf-formula1} and that in the right-hand side of~\eqref{eq:survfun-zero-formula1} are convergent for all $y\in[0,A]$ but only for $t>0$. This is not a big problem, because, as we discussed above, at $t=0$, either of the two survival functions is identically equal to $1$, which is a consequence of the definition~\eqref{eq:T-GSR-def} of the GSR stopping time $\mathcal{S}_A^r$.

To conclude this section, we note that the obtained formulae~\eqref{eq:survfun-inf-formula1} and~\eqref{eq:survfun-zero-formula1} can be simplified somewhat with the aid of~\cite[Identity~(2.4.21),~p.~25]{Slater:Book1960} according to which
\begin{align*}
(1/2-a-b)\,(1/2-a+b)\,W_{a-1,b}(z)
&=
(z/2-a)\,W_{a,b}(z)+z\left[\dfrac{\partial}{\partial z}\,W_{a,b}(z)\right].
\end{align*}
Specifically, setting $a=1-\theta$, $b=\imath\beta_{n,\theta,A}/2$, and $z=2/(\mu^2 A)$ in the foregoing identity leads to
\begin{align*}
\left[\left(\theta-\dfrac{1}{2}\right)^2+\dfrac{\beta_{n,\theta,A}^2}{4}\right]&W_{-\theta,\tfrac{\imath\beta_{n,\theta,A}}{2}}\left(\dfrac{2}{\mu^2 A}\right)
=\\
&=
\left(\dfrac{1}{\mu^2 A}-1+\theta\right)\,W_{1-\theta,\tfrac{\imath\beta_{n,\theta,A}}{2}}\left(\dfrac{2}{\mu^2 A}\right)+\dfrac{2}{\mu^2 A}\left.\left[\dfrac{\partial}{\partial u} W_{1-\theta,\tfrac{\imath\beta_{n,\theta,A}}{2}}(u)\right]\right|_{u=\frac{2}{\mu^2 A}},
\end{align*}
whence it follows that
\begin{align*}
\left.\left[\dfrac{\partial}{\partial u} W_{1-\theta,\tfrac{\imath\beta_{n,\theta,A}}{2}}(u)\right]\right|_{u=\frac{2}{\mu^2 A}}
&=
\dfrac{\mu^2 A}{8}\,(1+\beta_{n,\theta,A}^2)\, W_{-\theta,\tfrac{\imath\beta_{n,\theta,A}}{2}}\left(\dfrac{2}{\mu^2 A}\right),
\end{align*}
because
\begin{align*}
W_{1-\theta,\tfrac{\imath\beta_{n,\theta,A}}{2}}\left(\dfrac{2}{\mu^2 A}\right)
&=
0
\;\;\text{and}\;\;
\left(\theta-\dfrac{1}{2}\right)^2
=
\dfrac{1}{4}
\end{align*}
where the second identity is true because $\theta$ is either 0 or 1. Plugging this back into~\eqref{eq:Czero-eqn} we can conclude that
\begin{align*}
\dfrac{\mu^2 A}{2}\,\tilde{C}_{n,\theta,A}^2\,W_{-\theta,\tfrac{\imath\beta_{n,\theta,A}}{2}}\left(\dfrac{2}{\mu^2 A}\right)
&=
4\,\dfrac{\imath\beta_{n,\theta,A}}{1+\beta_{n,\theta,A}^2}\left/\Biggl\{\left.\Biggl[\dfrac{\partial}{\partial b}W_{1-\theta,b}\left(\dfrac{2}{\mu^2 A}\right)\Biggr]\right|_{b=\tfrac{\imath\beta_{n,\theta,A}}{2}}\Biggr\}\right.,
\end{align*}
for all $n\ge1$ and $\theta=\{0,1\}$. Likewise, by exactly the same argument, viz. merely by changing the notation $\imath\beta_{n,\theta,A}$ to $\alpha_{0,A}$, it can be shown that
\begin{align*}
\dfrac{\mu^2 A}{2}\,\tilde{C}_{0,0,A}^2\,W_{0,\tfrac{\alpha_{0,A}}{2}}\left(\dfrac{2}{\mu^2 A}\right)
&=
4\,\dfrac{\alpha_{0,A}}{1-\alpha_{0,A}^2}\left/\Biggl\{\left.\Biggl[\dfrac{\partial}{\partial b}W_{1,b}\left(\dfrac{2}{\mu^2 A}\right)\Biggr]\right|_{b=\tfrac{\alpha_{0,A}}{2}}\Biggr\}\right..
\end{align*}

Direct substitution of the last two identities into~\eqref{eq:survfun-inf-formula1} and into~\eqref{eq:survfun-zero-formula1} yields
\begin{align}\label{eq:survfun-inf-formula2}
\begin{split}
\Pr_{\infty}(\mathcal{S}_A^{r=y}\ge t)
&=
4\,\dfrac{y}{A}\,e^{\tfrac{1}{\mu^2 y}-\tfrac{1}{\mu^2 A}}e^{-\tfrac{\mu^2 t}{8}}\times\\
&\qquad
\times\Biggl\{e^{\tfrac{\mu^2 t}{8}\alpha_{0,A}^2}\,\dfrac{\alpha_{0,A}}{1-\alpha_{0,A}^2}\,W_{1,\tfrac{\alpha_{0,A}}{2}}\left(\dfrac{2}{\mu^2 y}\right)\left/\left.\Biggl[\dfrac{\partial}{\partial b}W_{1,b}\left(\dfrac{2}{\mu^2 A}\right)\Biggr]\right|_{b=\tfrac{\alpha_{0,A}}{2}}\right.+\\
&\qquad\qquad
+
\sum_{n=1}^{\infty}e^{-\tfrac{\mu^2 t}{8}\beta_{n,0,A}^2}
\,\dfrac{\imath\beta_{n,0,A}}{1+\beta_{n,0,A}^2}\times\\
&\qquad\qquad\qquad
\times\,W_{1,\tfrac{\imath\beta_{n,0,A}}{2}}\left(\dfrac{2}{\mu^2 y}\right)\left/\left.\Biggl[\dfrac{\partial}{\partial b}W_{1,b}\left(\dfrac{2}{\mu^2 A}\right)\Biggr]\right|_{b=\tfrac{\imath\beta_{n,0,A}}{2}}\right.\Biggr\},
\end{split}
\end{align}
and
\begin{align}\label{eq:survfun-zero-formula2}
\begin{split}
\Pr_{0}(\mathcal{S}_A^{r=y}\ge t)
&=
4\,e^{\tfrac{1}{\mu^2 y}-\tfrac{1}{\mu^2 A}}e^{-\tfrac{\mu^2 t}{8}}\times\\
&\qquad
\times
\sum_{n=1}^{\infty}e^{-\tfrac{\mu^2 t}{8}\beta_{n,1,A}^2}\,\dfrac{\imath\beta_{n,1,A}}{1+\beta_{n,1,A}^2}\times\\
&\qquad\qquad
\times\,W_{0,\tfrac{\imath\beta_{n,1,A}}{2}}\left(\dfrac{2}{\mu^2 y}\right)\left/\left.\Biggl[\dfrac{\partial}{\partial b}W_{0,b}\left(\dfrac{2}{\mu^2 A}\right)\Biggr]\right|_{b=\tfrac{\imath\beta_{n,1,A}}{2}}\right.,
\end{split}
\end{align}
where again $y\in[0,A]$ and $t>0$, and, by definition, either survival function is unity for $t=0$. We reiterate that $\alpha_{0,A}\in[0,1]$ and $\beta_{n,\theta,A}\ge0$, $n\ge1$, $\theta=\{0,1\}$, are found from equations~\eqref{eq:alpha-eqn} and~\eqref{eq:beta-eqn}, respectively. We also note that getting the survival functions corresponding to the classical SR procedure, i.e., when the GSR procedure has no headstart ($R_0^r=r=0$), is a matter of letting $r=y\to0+$ in both of the above formulae~\eqref{eq:survfun-inf-formula2}--\eqref{eq:survfun-zero-formula2} and making use of the asymptotics~\eqref{eq:Whit-fnc-asym-uinf} of the Whittaker $W$ function.

Despite the seemingly high complexity, the obtained formulae~\eqref{eq:survfun-inf-formula2} and~\eqref{eq:survfun-zero-formula2} are fully amenable to numerical evaluation ``as is'' using {\em Mathematica}. A corresponding numerical study is offered in the next section.

\section{A Numerical Study}
\label{sec:numerical-study}

This section's aim is to exploit numerically the expressions~\eqref{eq:survfun-inf-formula2} and~\eqref{eq:survfun-zero-formula2} obtained in the proceeding section for $\Pr_{\infty}(\mathcal{S}_A^r\ge t)$ and $\Pr_{0}(\mathcal{S}_A^r\ge t)$ to examine the statistical properties of the GSR stopping time $\mathcal{S}_A^r$ in the pre- as well as in the post-drift regimes. To that end, as is apparent from~\eqref{eq:survfun-inf-formula2} and~\eqref{eq:survfun-zero-formula2}, in either regime, the distribution of $\mathcal{S}_A^r$ depends on:\begin{inparaenum}[\itshape(a)]\item the magnitude of the drift $\mu\neq0$, \item the detection threshold $A>0$, and \item the headstart $R_0^r=r\in[0,A]$\end{inparaenum}. To demonstrate the effect of each of these factors on the GSR stopping time's distribution and to see how the $\Pr_{\infty}$-statistical profile of $\mathcal{S}_A^r$ is different from the $\Pr_{0}$-one, we have put together a {\em Mathematica} script that evaluates both survival functions as bivariate functions of $R_0^r=r\in[0,A]$ and $t\ge0$ for any given $\mu\neq0$ and $A>0$. Moreover, the script also evaluates the corresponding densities $-\partial\Pr_{\infty}(\mathcal{S}_A^r\ge t)/\partial t$ and $-\partial\Pr_{0}(\mathcal{S}_A^r\ge t)/\partial t$, also as bivariate functions of $R_0^r=r\in[0,A]$ and $t>0$ for any given $\mu\neq0$ and $A>0$. Since the temporal and spacial variables are separated in both~\eqref{eq:survfun-inf-formula2} and~\eqref{eq:survfun-zero-formula2}, the survival functions' densities are straightforward to find explicitly by direct differentiation of~\eqref{eq:survfun-inf-formula2} and~\eqref{eq:survfun-zero-formula2} with respect to time.

To get a bit more technical, in order to guarantee reasonable accuracy our {\em Mathematica} script truncates the infinite series in the right-hand side of~\eqref{eq:survfun-inf-formula2} and that in the right-hand side of~\eqref{eq:survfun-zero-formula2} to the first 500 (five hundred) terms. Empirically, ``chopping off'' the two infinite series that far proved to be more than sufficient to prevent any significant loss of accuracy, at least for practically important parameter values. To boot, the (at most one) solution $\alpha_{0,A}$ of equation~\eqref{eq:alpha-eqn} and the first 500 solutions $\beta_{n,\theta,A}$, $1\le n\le 500$, of equation~\eqref{eq:beta-eqn} for both $\theta=\{0,1\}$ are all computed to within 400 (four hundred) decimal places. Incidentally, the decision to use {\em Mathematica} (instead of, e.g., MATLAB developed by MathWorks, Inc.) was made because of {\em Mathematica}'s phenomenal ability to handle special functions, especially the Whittaker $W$ function. In particular, {\em Mathematica} turned out to be capable of computing the first derivative of the Whittaker $W$ function with respect to its second index, and doing so not only with high precision but also fairly quickly; recall that the first derivative of the Whittaker $W$ function with respect to its second index is involved in both survival functions' formulae~\eqref{eq:survfun-inf-formula2} and~\eqref{eq:survfun-zero-formula2}. For a given pair of $\mu\neq0$ and $A>0$, and for each particular $\theta=\{0,1\}$, we established experimentally that it takes our {\em Mathematica} script about three hours to complete all the calculations and ``spit out'' a vector containing (approximate) values of $\alpha_{0,A}$ (for $\theta=0$ only) and $\beta_{n,\theta,A}$, $1\le n\le 500$. This is assuming the script is run on an average office PC. To speed it up, we parallelized the calculations using the high throughput computing (HTC) infrastructure available at the Department of Mathematical Sciences at SUNY Binghamton. The use of the HTC infrastructure enabled us to boost the script's productivity up by a factor of about ten, depending on the number of available {\em Mathematica} licenses. Once $\alpha_{0,A}$ (for $\theta=0$ only) and $\beta_{n,\theta,A}$, $1\le n\le 500$, $\theta=\{0,1\}$, are all found, the actual evaluation of the corresponding survival function is merely a matter of plugging the obtained $\alpha_{0,A}$ (for $\theta=0$ only) and $\beta_{n,\theta,A}$, $1\le n\le 500$, back into the appropriate formula, either~\eqref{eq:survfun-inf-formula2} or~\eqref{eq:survfun-zero-formula2}.

However, before we present our numerical results, we would like to point out that it would be desirable to have a way to somehow validate the obtained numbers. To that end, some degree of confidence can be obtained using the observation that
\begin{align}\label{eq:firstmoment-int}
\EV_{\infty}[\mathcal{S}_{A}^{r}]
&=
\int_{0}^{\infty}\Pr_{\infty}(\mathcal{S}_{A}^{r}\ge t)\,dt
\;\;\text{and}\;\;
\EV_{0}[\mathcal{S}_{A}^{r}]
=
\int_{0}^{\infty}\Pr_{0}(\mathcal{S}_{A}^{r}\ge t)\,dt,
\end{align}
where $\EV_{\infty}[\cdot]$ and $\EV_{0}[\cdot]$ denote the expectations under the probability measures $\Pr_{\infty}$ and $\Pr_{0}$, respectively. More concretely, the idea is that both first moments $\EV_{\infty}[\mathcal{S}_{A}^{r}]$ and $\EV_{0}[\mathcal{S}_{A}^{r}]$ have actually been already found exact closed-form expressions for in the literature, although using a different approach. Specifically, it is considered a classical result that $\EV_{\infty}[\mathcal{S}_{A}^{r}]=A-r$, $r\in[0,A]$, and that
\begin{align}\label{eq:ADD_zero_formula}
\EV_{0}[\mathcal{S}_{A}^{r}]
&=
\frac{2}{\mu^2}\left\{e^{\tfrac{2}{\mu^2 A}}\,\left[-\Ei\left(-\frac{2}{\mu^2 A}\right)\right]-e^{\tfrac{2}{\mu^2 r}}\,\left[-\Ei\left(-\frac{2}{\mu^2 r}\right)\right]\right\},\;\; r\in[0,A],
\end{align}
where
\begin{empheq}[%
    left={%
        \Ei(x)\triangleq%
    \empheqlbrace}]{align*}
&-\displaystyle\int_{-x}^\infty\dfrac{e^{-t}}{t}\,dt,\;\text{for $x<0$;}\\
&-\lim_{\varepsilon\to+0}\left[\int_{-x}^{-\varepsilon}\dfrac{e^{-t}}{t}\,dt+\int_{\varepsilon}^{\infty}\dfrac{e^{-t}}{t}\,dt\right],\;\text{for $x>0$},
\end{empheq}
is the exponential integral; for a background on the exponential integral, see, e.g.,~\cite[Chapter~5]{Abramowitz+Stegun:Handbook1964}. These pre- and post-change first moment formulae have been previously obtained, e.g., by~\cite{Pollak+Siegmund:B85},~\cite{Shiryaev:MathEvents2006},~\cite{Feinberg+Shiryaev:SD2006}, and~\cite{Burnaev:ARSAIM2009}, in the context of the quickest change-point detection problem. In particular, formula~\eqref{eq:ADD_zero_formula} is a trivial generalization of~\cite[Lemma~3.3,~p.~459]{Feinberg+Shiryaev:SD2006}. As a matter of fact, for a generic detection procedure given by stopping time $\T$, the pre-change first moment $\EV_{\infty}[\T]$ is known as the Average Run Length (ARL) to false alarm, i.e., $\ARL(\T)\triangleq\EV_{\infty}[\T]$, and the post-change first moment $\EV_{0}[\T]$ is known as the Average Detection Delay (ADD), i.e., $\ADD_{0}(\T)\triangleq\EV_{0}[\T]$. The ARL to false alarm and the ADD are the standard performance metrics commonly used in the minimax quickest change-point detection theory.

That is, if the two survival functions formulae~\eqref{eq:survfun-inf-formula2} and~\eqref{eq:survfun-zero-formula2} are correct, then, according to~\eqref{eq:firstmoment-int}, the integration of each with respect to time over the interval $[0,+\infty)$ must yield the corresponding first moment of the GSR stopping time $\mathcal{S}_A^r$, and the expressions for both of these moments have already been obtained. However, the problem with this idea is that the series involved in~\eqref{eq:survfun-inf-formula2} and~\eqref{eq:survfun-zero-formula2} are not convergent for $t=0$. Hence, the integration with respect to $t$ over the interval $[0,+\infty)$ is not an option. Nevertheless, it is possible (and rather simple) to integrate the series with respect to $t$ over the interval $[t^*,\infty)$ for any $t^*>0$. Then, by picking $t^*>0$ to be sufficiently small, it is reasonable to expect each integrated series to be close to the corresponding first moment. Therefore, should we determine that the $dt$-integral of the survival function over the interval $[t^*,\infty)$ for $t^*\approx 0$ provides an accurate approximation of the corresponding actual first moment (computed exactly via one of the aforementioned formulae), then we can be at least somewhat certain in the validity of the survival function formulae~\eqref{eq:survfun-inf-formula2} and~\eqref{eq:survfun-zero-formula2}, and, consequently, in the validity of our numerical results as well. Since this basic ``sanity check'' is not difficult to perform, we carried it out for each set of parameters we picked for our numerical study.

Specifically, for our study we picked two values of $\mu$: $\mu=0.5$ and $\mu=1.5$. These values correspond to small and contrast changes, respectively. We also note that since the survival functions formulae~\eqref{eq:survfun-inf-formula2} and~\eqref{eq:survfun-zero-formula2} are both symmetric with respect to the sign of $\mu$, it is sufficient to restrict attention to only positive $\mu$. As for the detection threshold $A>0$, we also picked two values: $A=10^2$ and $A=10^3$. When the GSR statistic's headstart is either zero or close to zero, these choices correspond to high and moderate false alarm risk levels, respectively.

We would like to organize the presentation of the numerical results as follows. For each of the selected values of $A$ and each $\theta=\{0,1\}$ we would like to report the results in a set of three rows of figures, where each row is two figures, shown one next to the other: the left one corresponds to $\mu=0.5$, and the right one corresponds to $\mu=1.5$. In each set, the first row of figures presents the corresponding first moment, either $\ARL(\mathcal{S}_A^r)\triangleq\EV_{\infty}[\mathcal{S}_A^r]$, if $\theta=0$, or $\ADD_{0}(\mathcal{S}_A^r)\triangleq\EV_{0}[\mathcal{S}_A^r]$, if $\theta=1$, shown as a function of the headstart $r\in[0,A]$. Specifically, each plot of the first moment shows two curves: the first moment computed exactly, using the aforementioned formulae, and the first moment computed by integrating the corresponding survival function with respect to time over the interval $[t^{*},\infty)$ with $t^{*}=10^{-3}$. The first moment computed exactly is shown as a smooth gray curve, and the first moment computed off the survival function is shown as a sequence of separate solid dark dots. With regard to computing the first moment off the survival function, we note that since the interval of integration $[t^{*},\infty)$ starts at $t^{*}=10^{-3}$, i.e., pretty close to zero, it is reasonable to expect each one of the dark dots to lie on the gray curve, should, of course, the expression---either~\eqref{eq:survfun-inf-formula2} or~\eqref{eq:survfun-zero-formula2}---for the corresponding survival function be actually correct, and should the corresponding numerical error be acceptably small.

The second row of figures is intended to present the corresponding density, either $-\partial\Pr_{\infty}(\mathcal{S}_A^r\ge t)/\partial t$, if $\theta=0$, or $-\partial\Pr_{0}(\mathcal{S}_A^r\ge t)/\partial t$, if $\theta=1$. Specifically, the figures show the density as a function of the headstart $r\in[0,A]$ and time $t$ restricted to the interval $[0,10]$. Finally, the third row of figures reports the corresponding survival function, again either $\Pr_{\infty}(\mathcal{S}_A^r\ge t)$, if $\theta=0$, or $\Pr_{0}(\mathcal{S}_A^r\ge t)$, if $\theta=1$, shown as a function of the headstart $r\in[0,A]$ and time $t$, again restricted to the interval $[0,10]$. We recall that, although the survival function formulae~\eqref{eq:survfun-inf-formula2} and~\eqref{eq:survfun-zero-formula2} do not hold for $t=0$, it follows from the definition~\eqref{eq:T-GSR-def} for the GSR stopping time $\mathcal{S}_A^r$ that at $t=0$ either survival function is unity.

\begin{figure}[t!]
    \centering
    \begin{subfigure}{0.45\textwidth}
        \centering
        \includegraphics[width=\linewidth]{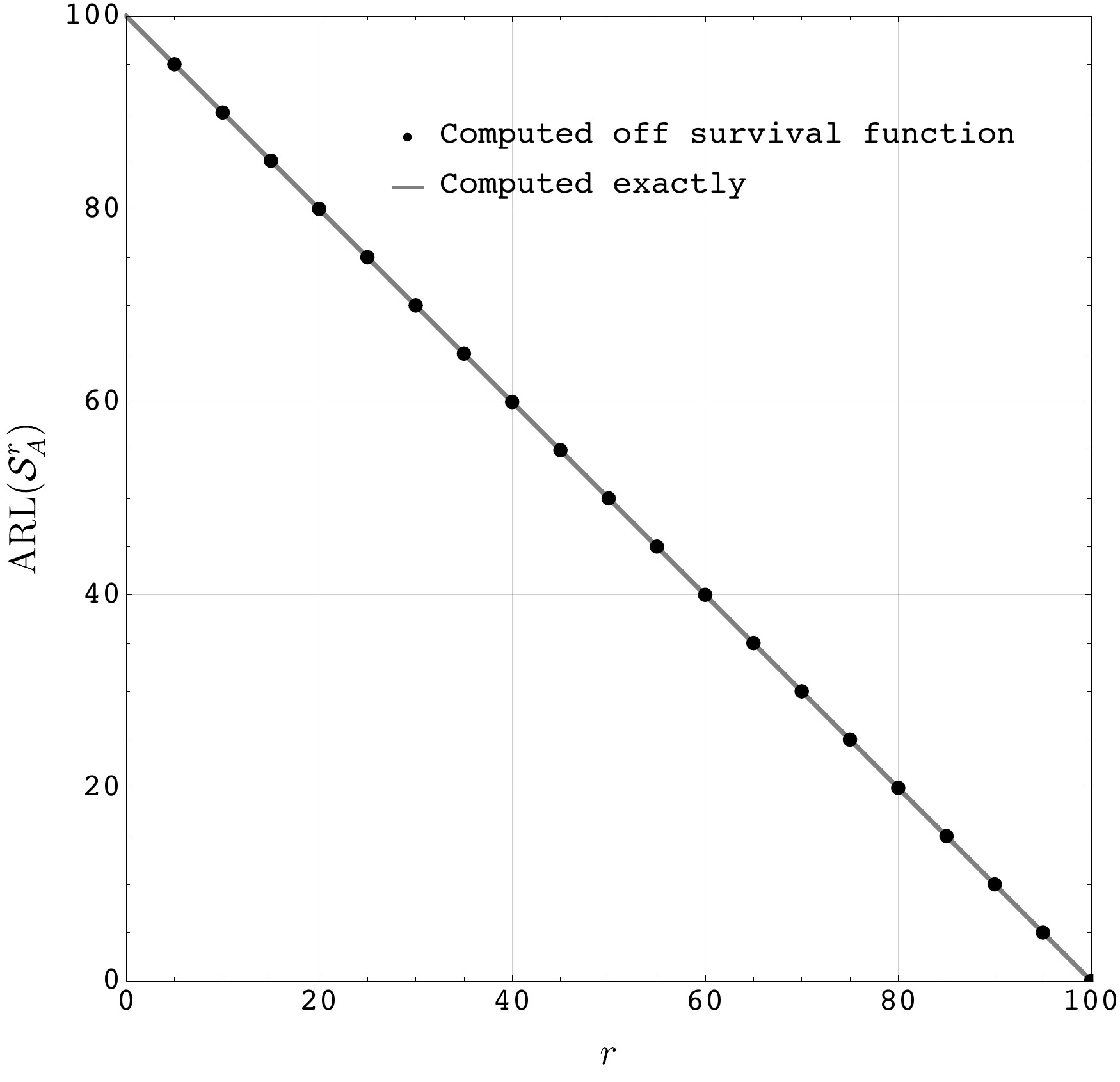}
        \caption{$\mu=0.5$.}
        \label{fig:firstmoment_inf__mu_1over2_A_100}
    \end{subfigure}
    \hspace*{\fill}
    \begin{subfigure}{0.45\textwidth}
        \centering
        \includegraphics[width=\linewidth]{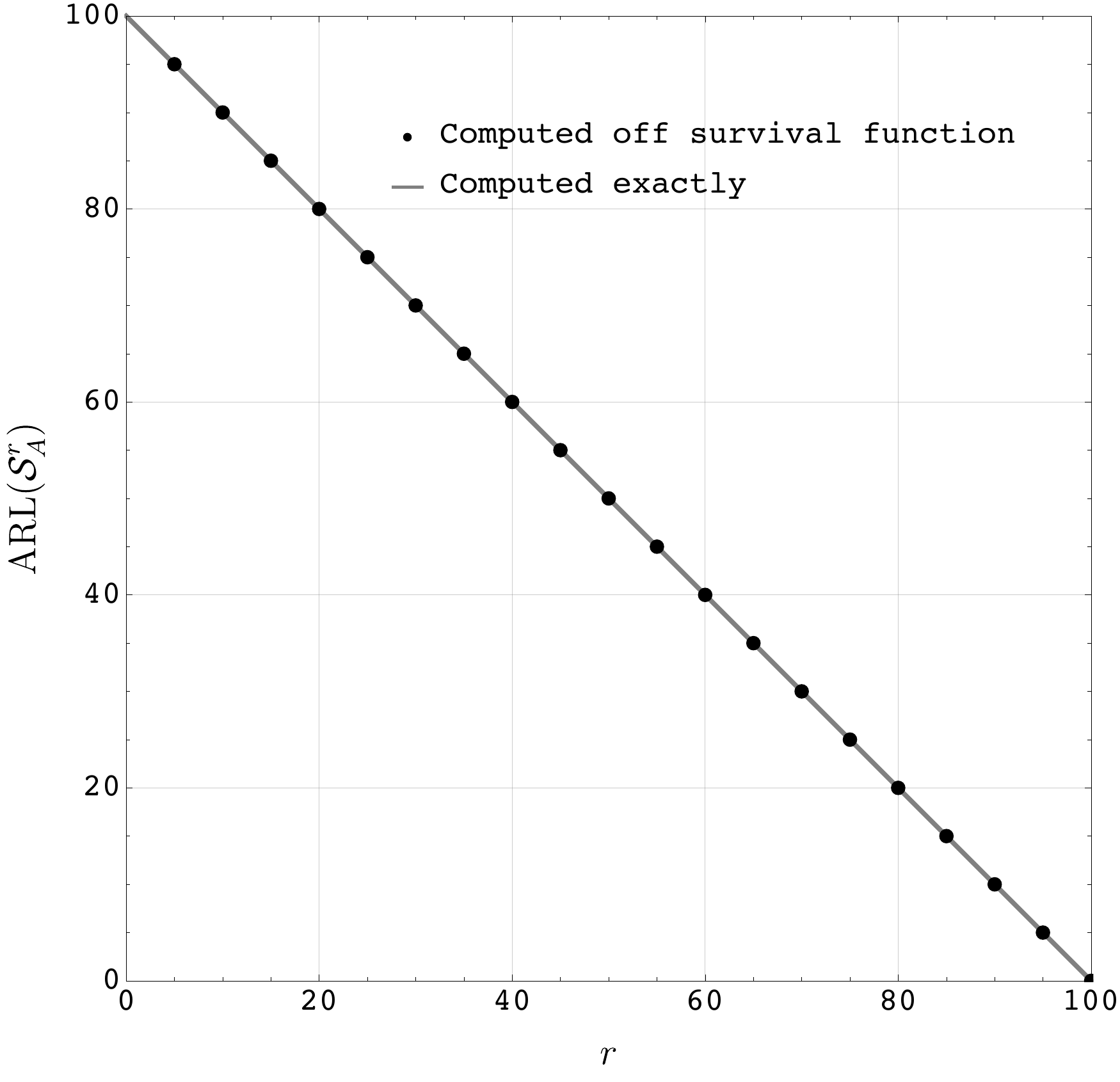}
        \caption{$\mu=1.5$.}
        \label{fig:firstmoment_inf__mu_3over2_A_100}
    \end{subfigure}
    \caption{Pre-change first moment of $\mathcal{S}_A^r$, i.e., $\ARL(\mathcal{S}_A^r)\triangleq\EV_{\infty}[\mathcal{S}_A^r]$, as a function of $r\in[0,A]$ for $A=10^2$ and $\mu=\{0.5,1.5\}$.}
    \label{fig:firstmoment_inf_A_100}
\end{figure}

\begin{figure}[t!]
    \centering
    \begin{subfigure}{0.45\textwidth}
        \centering
        \includegraphics[width=\linewidth]{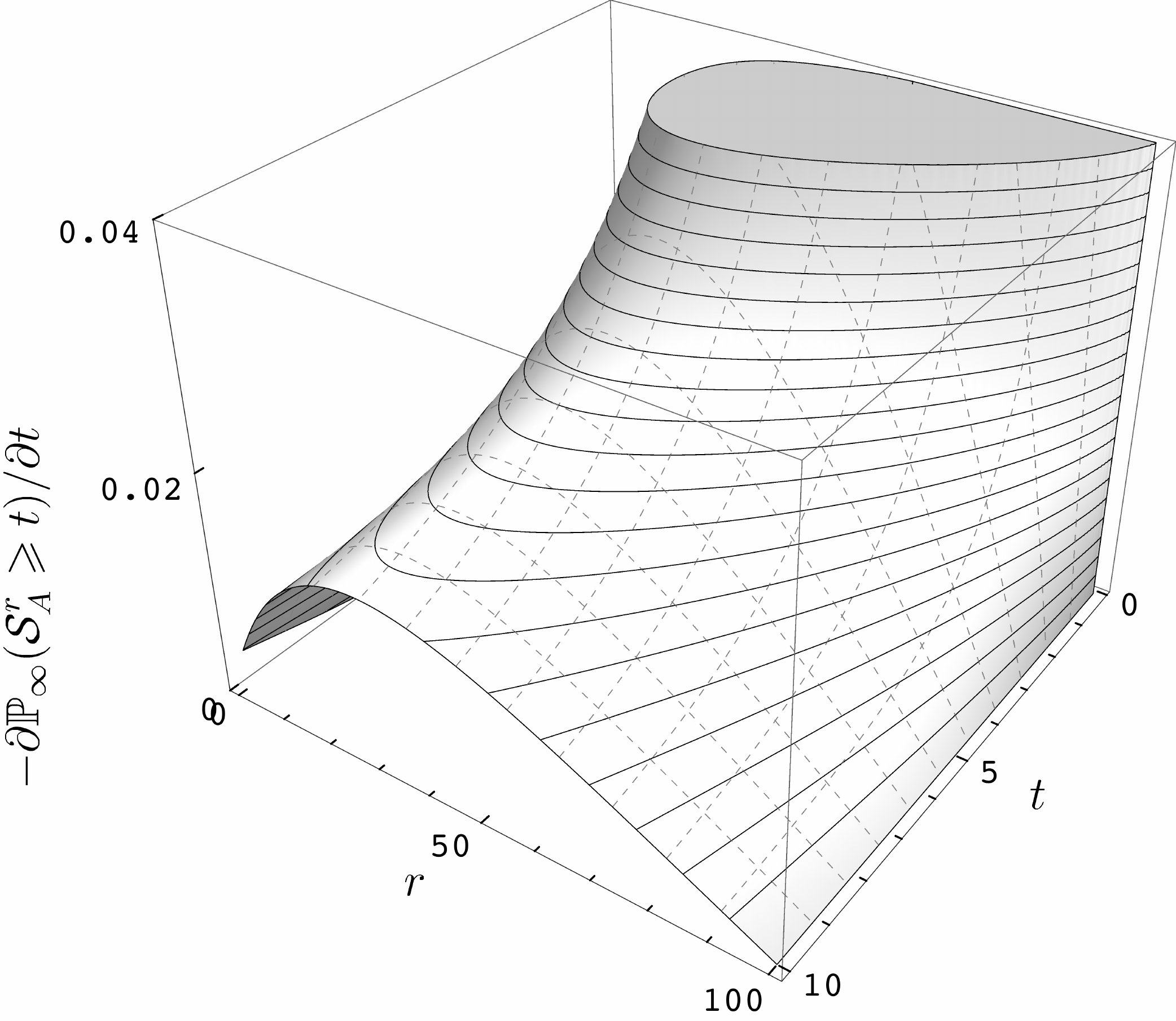}
        \caption{$\mu=0.5$.}
        \label{fig:survfun_pdf_inf__mu_1over2_A_100}
    \end{subfigure}
    \hspace*{\fill}
    \begin{subfigure}{0.45\textwidth}
        \centering
        \includegraphics[width=\linewidth]{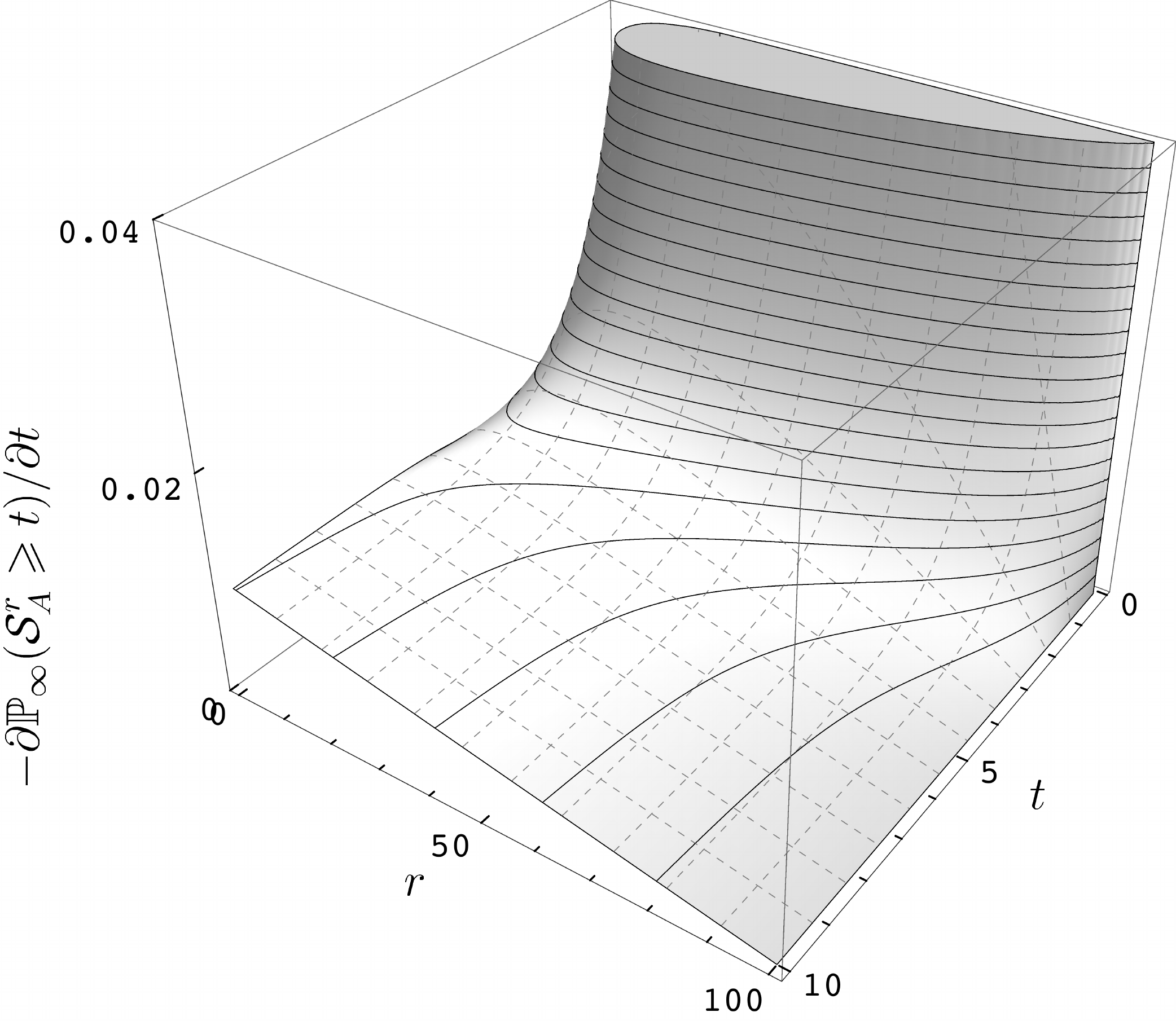}
        \caption{$\mu=1.5$.}
        \label{fig:survfun_pdf_inf__mu_3over2_A_100}
    \end{subfigure}
    \caption{Pre-change survival function density $-\partial\Pr_{\infty}(\mathcal{S}_A^r\ge t)/\partial t$ as a function of $t\in[0,10]$ and $r\in[0,A]$ for $A=10^2$ and $\mu=\{0.5,1.5\}$.}
    \label{fig:survfun_pdf_inf_A_100}
\end{figure}

\begin{figure}[t!]
    \centering
    \begin{subfigure}{0.45\textwidth}
        \centering
        \includegraphics[width=\linewidth]{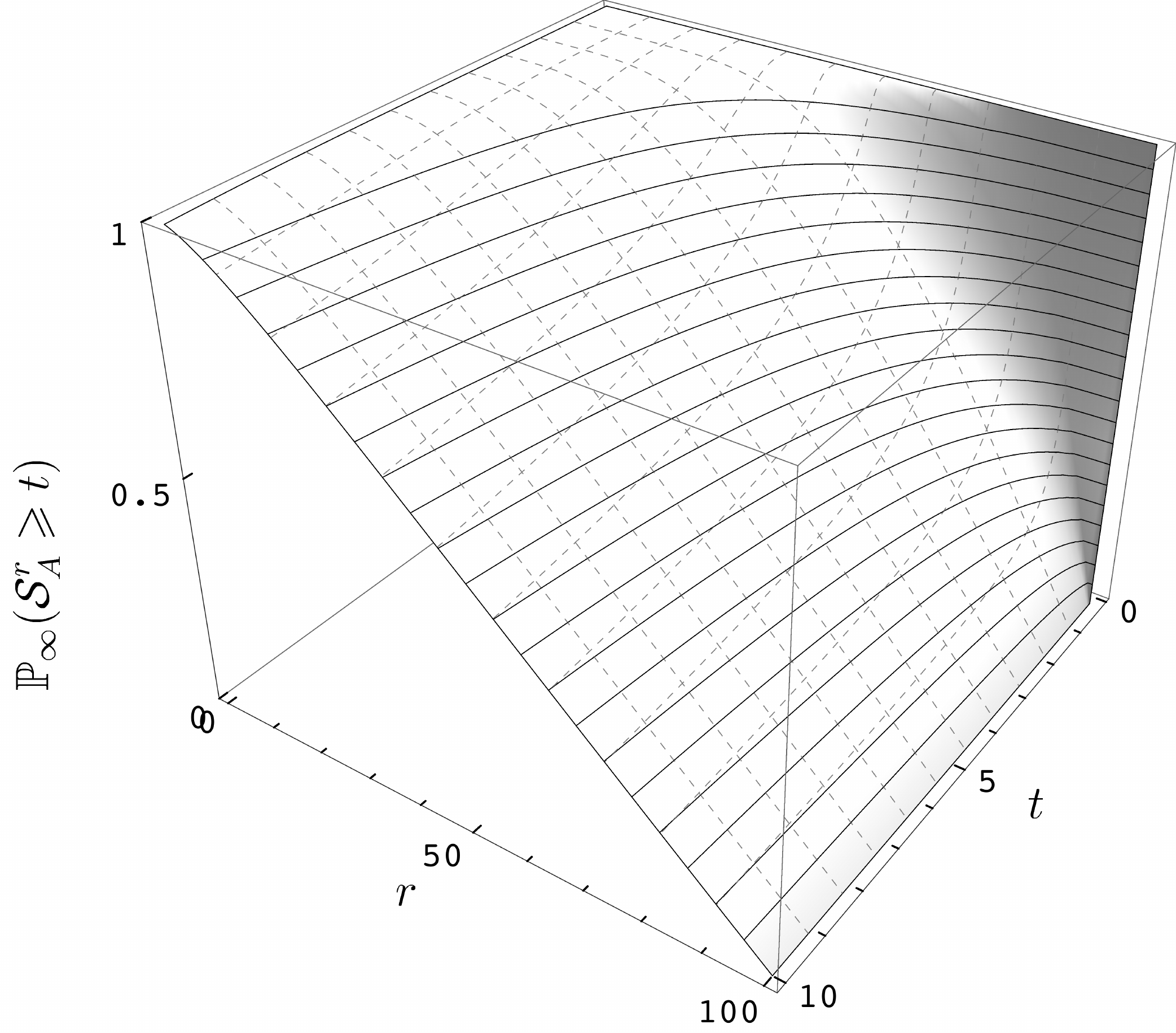}
        \caption{$\mu=0.5$.}
        \label{fig:survfun_inf__mu_1over2_A_100}
    \end{subfigure}
    \hspace*{\fill}
    \begin{subfigure}{0.45\textwidth}
        \centering
        \includegraphics[width=\linewidth]{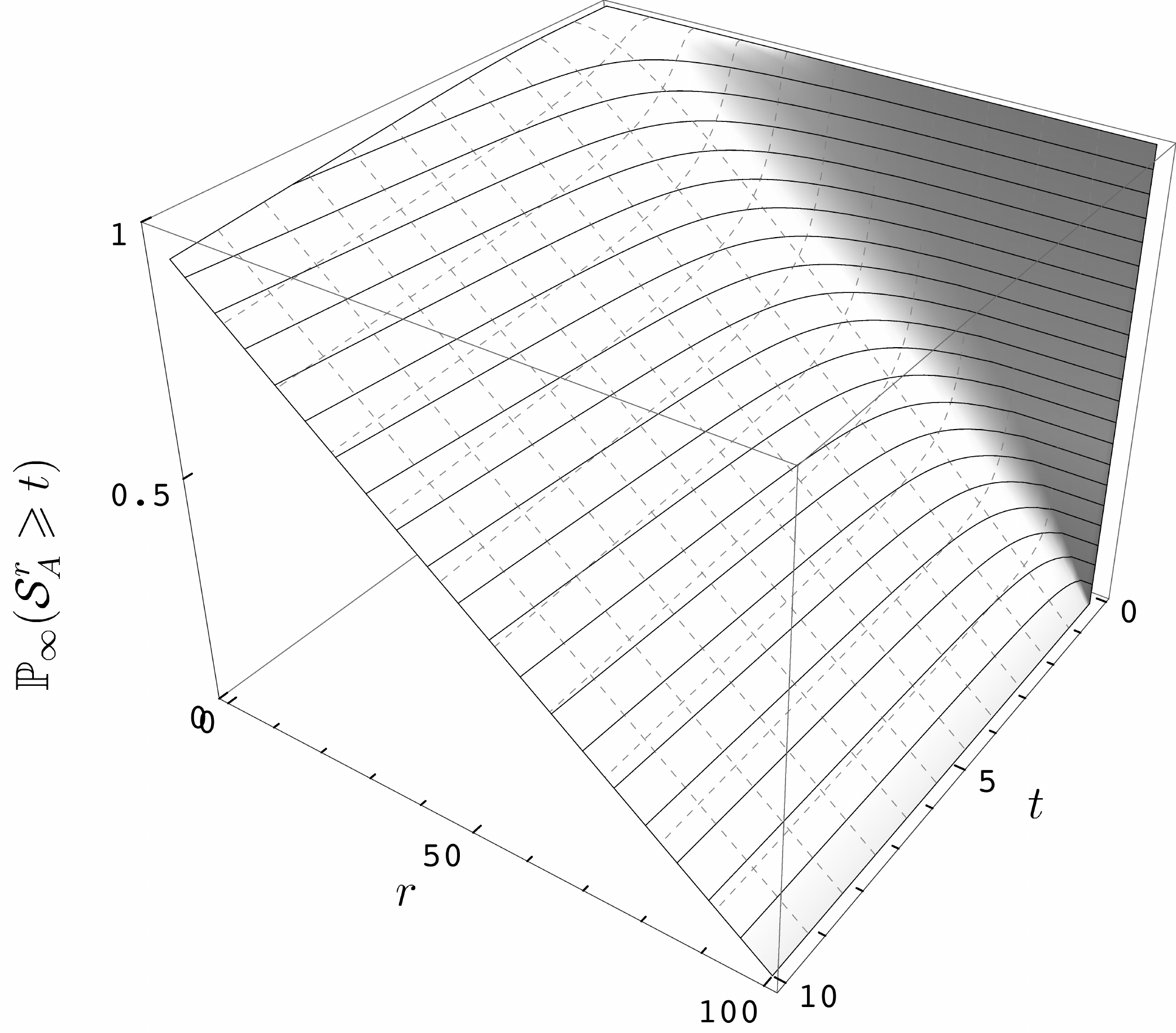}
        \caption{$\mu=1.5$.}
        \label{fig:survfun_inf__mu_3over2_A_100}
    \end{subfigure}
    \caption{Pre-change survival function $\Pr_{\infty}(\mathcal{S}_A^r\ge t)$ as a function of $t\in[0,10]$ and $r\in[0,A]$ for $A=10^2$ and $\mu=\{0.5,1.5\}$.}
    \label{fig:survfun_inf_A_100}
\end{figure}

With all of the above in mind, we now begin our study. The first set of plots is given by Figures~\ref{fig:firstmoment_inf_A_100},~\ref{fig:survfun_pdf_inf_A_100}, and~\ref{fig:survfun_inf_A_100}. These figures all assume that $A=100$ and $\theta=0$, i.e., they all correspond to the pre-change regime with high false alarm risk. We would like to immediately draw attention to Figures~\ref{fig:firstmoment_inf_A_100}. These figures show the first moment, i.e., $\ARL(\mathcal{S}_A^r)\triangleq\EV_{\infty}[\mathcal{S}_A^r]$, as a function of the headstart. The exact values correspond to the solid gray line, and the values computed off the survival function (by means of integration of the survival function formula with respect to $t$) are shown as isolated solid dark dots. The fact that for both values of $\mu$ the dark dots are in perfect agreement with the gray curve provides evidence that the survival function formula~\eqref{eq:survfun-inf-formula2} is likely to be correct. Moreover, it also asserts (at least to some extent) that the accuracy of our numerical results is sufficiently high. As a side comment we note that Figure~\ref{fig:firstmoment_inf__mu_1over2_A_100} which corresponds to $\mu=0.5$ and Figure~\ref{fig:firstmoment_inf__mu_3over2_A_100} which corresponds to $\mu=1.5$ are nearly identical. This is because $\ARL(\mathcal{S}_A^r)=A-r$ for {\em any} $\mu$, and this is a direct consequence the well-known fact that $(R_t^r-r-t)_{t\ge0}$ is a zero-mean $\Pr_{\infty}$-martingale.

Let us next look at Figures~\ref{fig:survfun_pdf_inf_A_100}. These figures show the density $-\partial\Pr_{\infty}(\mathcal{S}_A^r\ge t)/\partial t$ as a function of $r\in[0,A]$ and $t$ between $0$ and $10$. Specifically, Figure~\ref{fig:survfun_pdf_inf__mu_1over2_A_100} corresponds to $\mu=0.5$ and Figure~\ref{fig:survfun_pdf_inf__mu_3over2_A_100} is for $\mu=1.5$. We note that, for either value of $\mu$, the surface has a spike concentrated around the point $r=A$ and $t\approx 0$. This is simple to explain: when $t$ is close to zero, the survival function is close to unity, unless the headstart is close to the detection threshold. Then, as time increases, the surface flattens out, i.e., the spike dissolves, which indicates that the headstart becomes less of a factor. This also makes perfect sense, because obviously $\Pr_{\infty}(\mathcal{S}_A^r\ge t)$ must tend to zero as $t$ increases, whatever by the headstart $r$. This conclusion can be formally reached, e.g., from the Markov inequality. To understand the effect of $\mu$, note that from the $\Pr_\infty$-differential $dR_t^r=dt+\mu R_t^r\,dB_t$, which comes from~\eqref{eq:Rt_r-combined-SDE} with $\theta=0$, it is clear that the variance of $R_t^r$ is directly proportional to $\mu^2$. Hence, in the pre-change regime, the GSR statistic is more volatile for higher values of $\mu$. In terms of the GSR stopping time, this means that the GSR procedure is more likely to stop sooner when $\abs{\mu}$ is large than when $\abs{\mu}$ is small. This is the reason why the density $-\partial\Pr_{\infty}(\mathcal{S}_A^r\ge t)/\partial t$ appears to be more flat (i.e., as though it was ``stretched'' along the $r$-axis) when $\mu=1.5$ than when $\mu=0.5$.

To complete the presentation of the first set of results, Figures~\ref{fig:survfun_inf_A_100} show the corresponding survival function $\Pr_{\infty}(\mathcal{S}_A^r\ge t)$. Specifically, Figure~\ref{fig:survfun_inf__mu_1over2_A_100} shows the survival function for $\mu=0.5$ and Figure~\ref{fig:survfun_inf__mu_3over2_A_100} assumes $\mu=1.5$. As one would expect, the survival function can be seen to be a decreasing function of both the headstart $r\in[0,A]$ and time $t>0$. For the reasons we already explained above, when $\mu=1.5$, the survival function decays down to zero (with respect to both variables) faster than when $\mu=0.5$. However, the rate of decay for $\mu=1.5$ is only slightly higher than that for $\mu=0.5$.

\begin{figure}[t!]
    \centering
    \begin{subfigure}{0.45\textwidth}
        \centering
        \includegraphics[width=\linewidth]{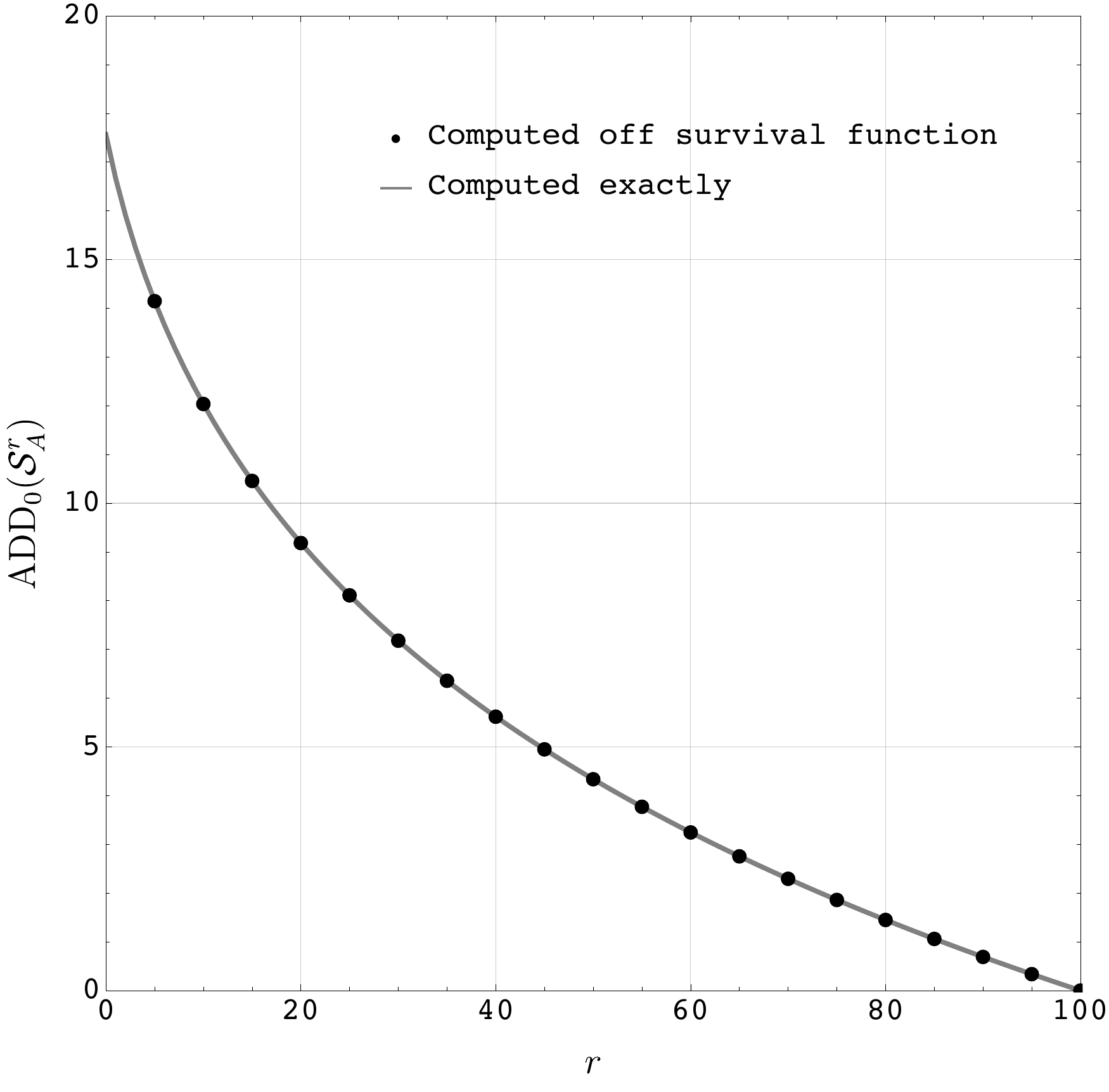}
        \caption{$\mu=0.5$.}
        \label{fig:firstmoment_zero__mu_1over2_A_100}
    \end{subfigure}
    \hspace*{\fill}
    \begin{subfigure}{0.45\textwidth}
        \centering
        \includegraphics[width=\linewidth]{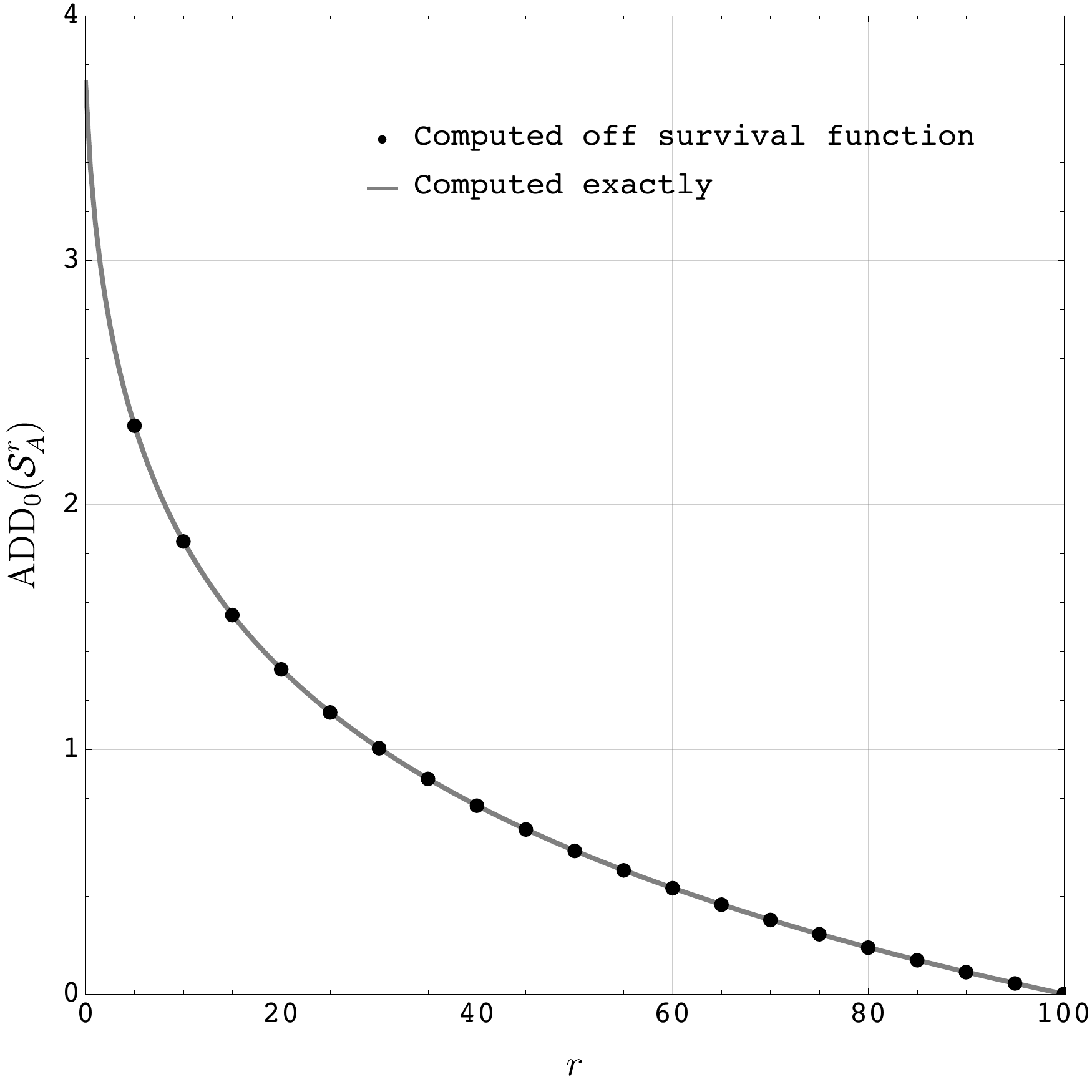}
        \caption{$\mu=1.5$.}
        \label{fig:firstmoment_zero__mu_3over2_A_100}
    \end{subfigure}
    \caption{Post-change first moment of $\mathcal{S}_A^r$, i.e., $\ADD_{0}(\mathcal{S}_A^r)\triangleq\EV_{0}[\mathcal{S}_A^r]$, as a function of $r\in[0,A]$ for $A=10^2$ and $\mu=\{0.5,1.5\}$.}
    \label{fig:firstmoment_zero_A_100}
\end{figure}
\begin{figure}[t!]
    \centering
    \begin{subfigure}{0.45\textwidth}
        \centering
        \includegraphics[width=\linewidth]{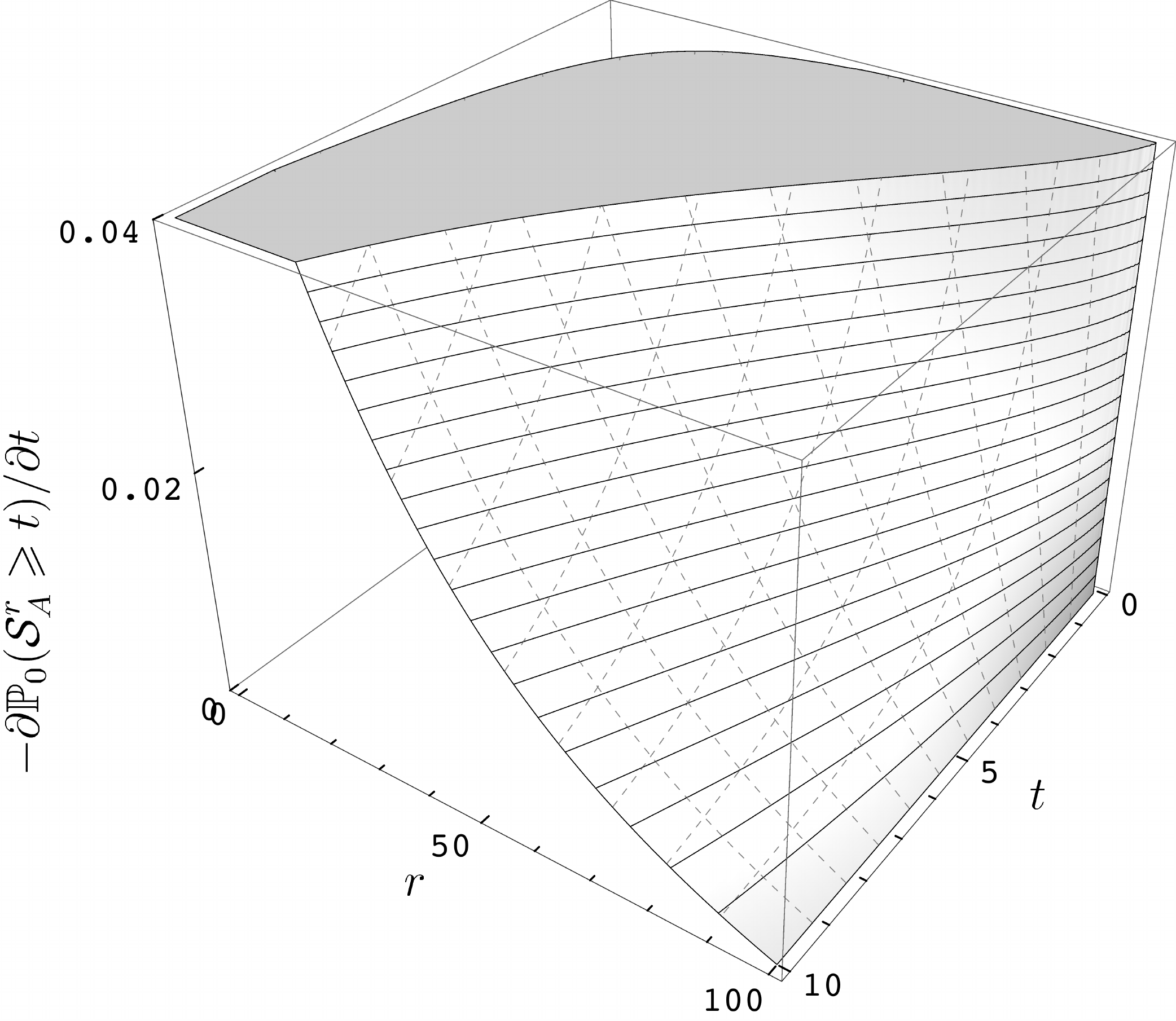}
        \caption{$\mu=0.5$.}
        \label{fig:survfun_pdf_zero__mu_1over2_A_100}
    \end{subfigure}
    \hspace*{\fill}
    \begin{subfigure}{0.45\textwidth}
        \centering
        \includegraphics[width=\linewidth]{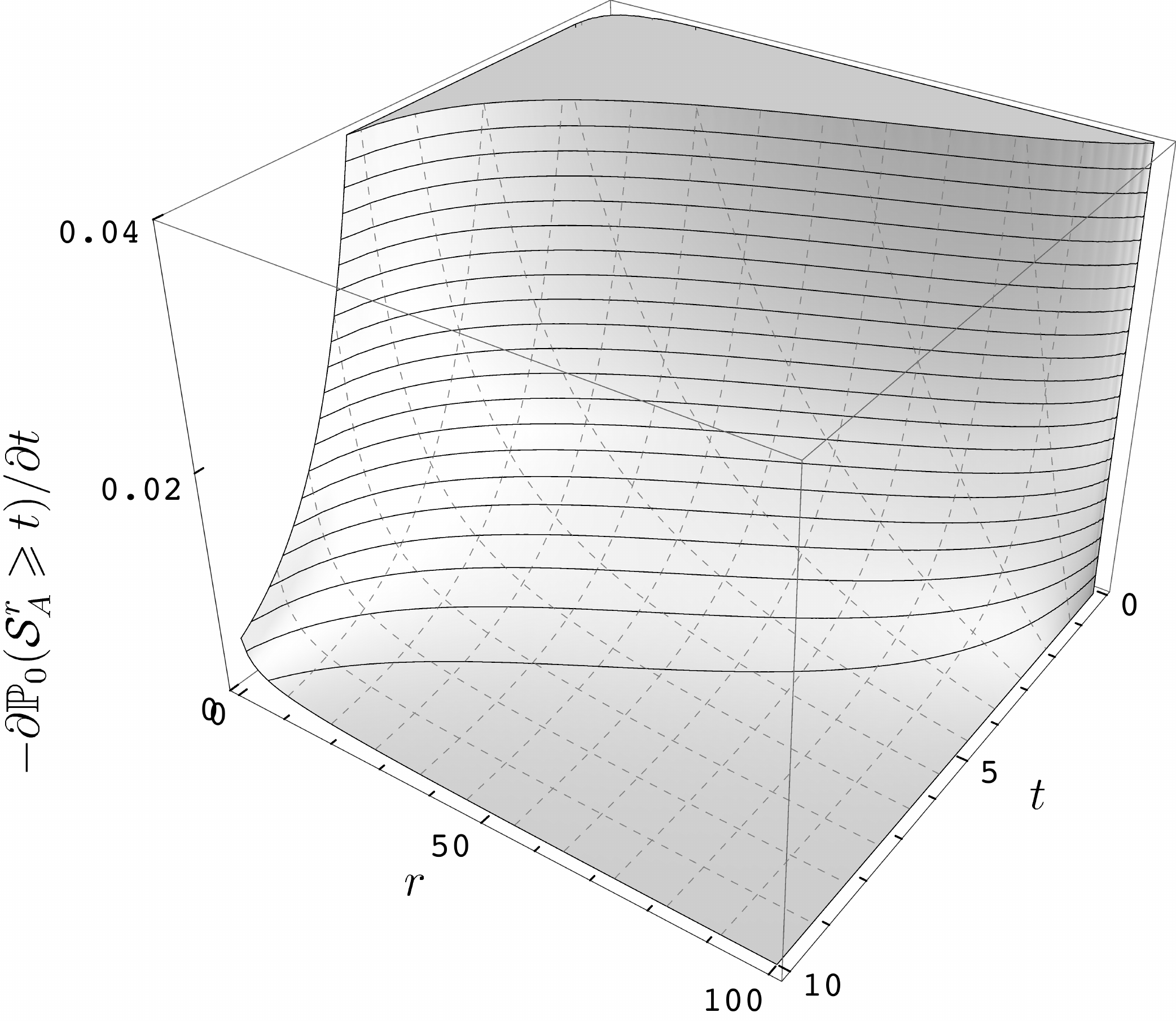}
        \caption{$\mu=1.5$.}
        \label{fig:survfun_pdf_zero__mu_3over2_A_100}
    \end{subfigure}
    \caption{Post-change survival function density $-\partial\Pr_{0}(\mathcal{S}_A^r\ge t)/\partial t$ as a function of $t\in[0,10]$ and $r\in[0,A]$ for $A=10^2$ and $\mu=\{0.5,1.5\}$.}
    \label{fig:survfun_pdf_zero_A_100}
\end{figure}
\begin{figure}[t!]
    \centering
    \begin{subfigure}{0.45\textwidth}
        \centering
        \includegraphics[width=\linewidth]{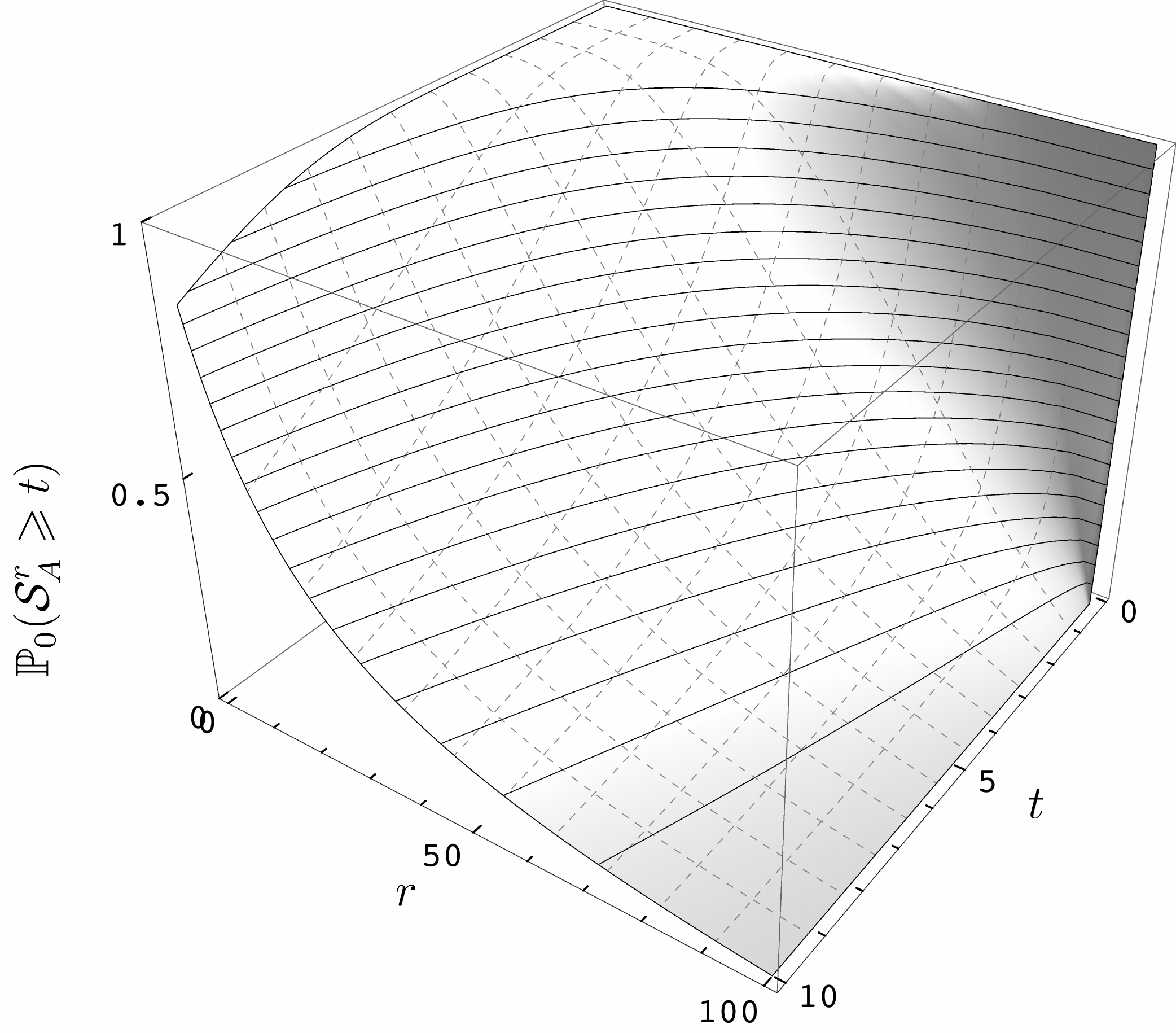}
        \caption{$\mu=0.5$.}
        \label{fig:survfun_zero__mu_1over2_A_100}
    \end{subfigure}
    \hspace*{\fill}
    \begin{subfigure}{0.45\textwidth}
        \centering
        \includegraphics[width=\linewidth]{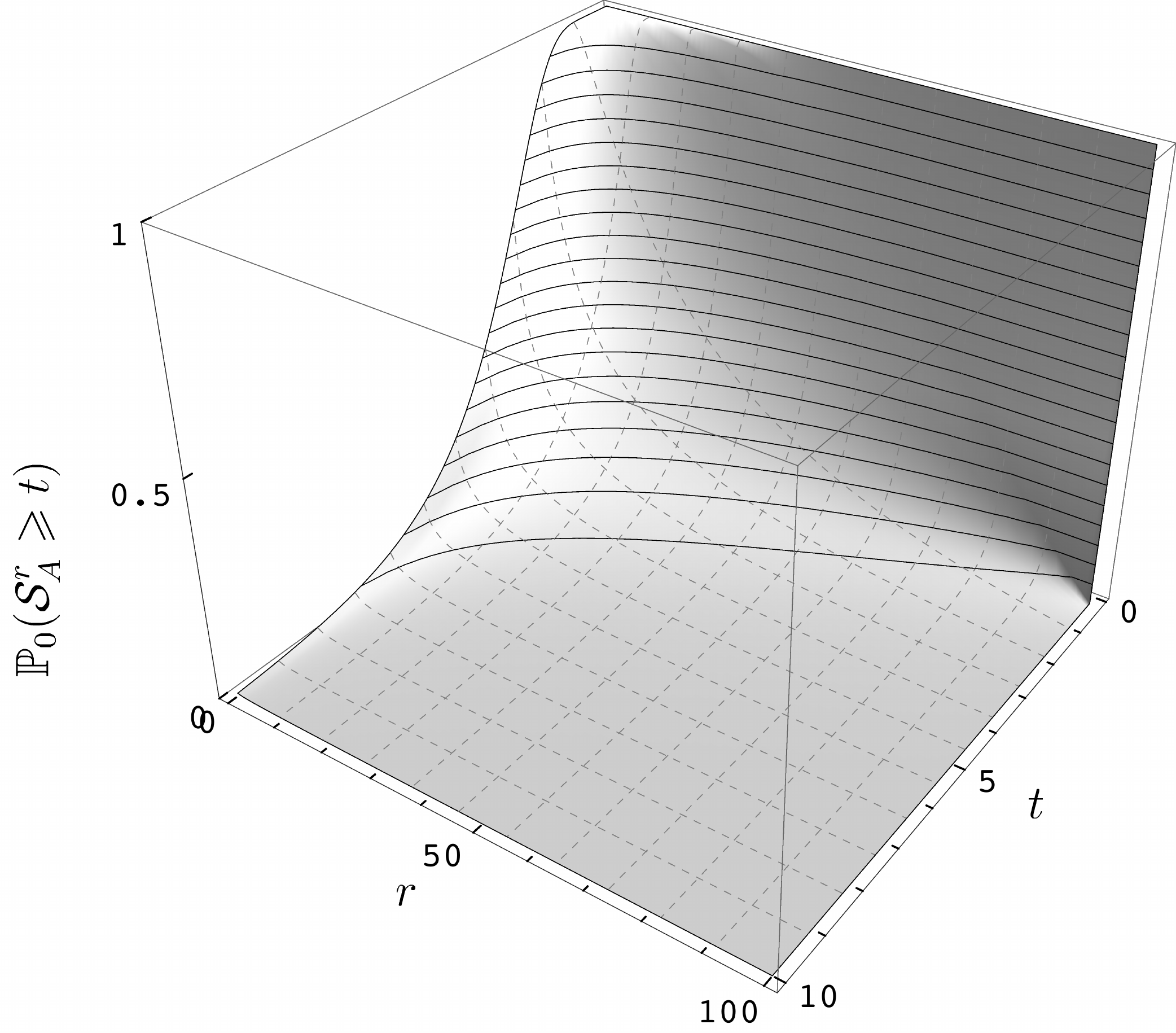}
        \caption{$\mu=1.5$.}
        \label{fig:survfun_zero__mu_3over2_A_100}
    \end{subfigure}
    \caption{Post-change survival function $\Pr_{0}(\mathcal{S}_A^r\ge t)$ as a function of $t\in[0,10]$ and $r\in[0,A]$ for $A=10^2$ and $\mu=\{0.5,1.5\}$.}
    \label{fig:survfun_zero_A_100}
\end{figure}

Let us now see what happens in the post-change regime. To that end, the first set of plots for the case when $A=100$ is formed by Figures~\ref{fig:firstmoment_zero_A_100},~\ref{fig:survfun_pdf_zero_A_100}, and~\ref{fig:survfun_zero_A_100}. As before, we hasten to note the perfect agreement seen in Figures~\ref{fig:firstmoment_zero_A_100} of the values of the corresponding first moment $\ADD_{0}(\mathcal{S}_A^r)\triangleq\EV_{0}[\mathcal{S}_A^r]$ computed exactly and off the survival function. Therefore, once again, the basic ``sanity check'' is successfully passed. However, unlike the pre-change first moment shown Figures~\ref{fig:firstmoment_inf_A_100}, the post-change first moment is not independent of the drift $\mu$, and this is trivial to see from formula~\eqref{eq:ADD_zero_formula}. Therefore, Figure~\ref{fig:firstmoment_zero__mu_1over2_A_100} which shows that post-change first moment for $\mu=0.5$ is actually different from Figure~\ref{fig:firstmoment_zero__mu_3over2_A_100} which shows the post-change first moment for $\mu=1.5$. Specifically, the difference between Figures~\ref{fig:firstmoment_zero__mu_1over2_A_100} and~\ref{fig:firstmoment_zero__mu_3over2_A_100} is in the scale along the vertical axis: for Figure~\ref{fig:firstmoment_zero__mu_1over2_A_100} which corresponds to $\mu=0.5$ the scale along the vertical axis is higher than for Figure~\ref{fig:firstmoment_zero__mu_3over2_A_100} which corresponds to $\mu=1.5$. This is because the post-change first moment $\ADD_{0}(\mathcal{S}_A^r)\triangleq\EV_{0}[\mathcal{S}_A^r]$ represents the average delay to detection (delivered by the GSR procedure), and more contrast changes (higher values of $\abs{\mu}$) are generally detected quicker (i.e., with a lower detection delay) than less contrast changes (lower values of $\abs{\mu}$).

The corresponding density $-\partial\Pr_{0}(\mathcal{S}_A^r\ge t)/\partial t$ is presented in Figures~\ref{fig:survfun_pdf_zero_A_100}. Compared to the pre-change density $-\partial\Pr_{\infty}(\mathcal{S}_A^r\ge t)/\partial t$ shown in Figures~\ref{fig:survfun_pdf_inf_A_100}, the density in the post-change regime is more ``heavy'', i.e., more probability mass is concentrated around the origin. This shouldn't come as a surprise, because, all other things being equal, the post-change first moment $\EV_{0}[\mathcal{S}_A^r]$ is much smaller than the pre-change first moment $\EV_{\infty}[\mathcal{S}_A^r]$. As a matter of fact, it is well-known in quickest change-point detection, that $\EV_{0}[\mathcal{S}_A^r]$ is asymptotically (as $A\to\infty$) on the order of $\log\EV_{\infty}[\mathcal{S}_A^r]$. Alternatively, recall from~\eqref{eq:Rt_r-combined-SDE} that the $\Pr_{0}$-differential of the GSR statistic is $dR_t^r=(1+\mu^2 R_t^r)dt+\mu R_t^r dB_t$. Therefore, since the post-change instantaneous drift function $b(x)=1+\mu^2 x$ dominates its pre-change counterpart $b(x)=1$, it follows that the GSR statistic $R_t^r$ grows faster in the post-change regime than in the pre-change regime. In fact, this is exactly how the GSR statistic ``senses'' the presence of the change (drift) in the observed standard Brownian motion. In addition, note that because the post-change first moment is higher for $\mu=0.5$ than for $\mu=1.5$---see, respectively, Figures~\ref{fig:firstmoment_zero__mu_1over2_A_100} and~\ref{fig:firstmoment_zero__mu_3over2_A_100} above---it follows that the survival function density should tend to zero quicker for $\mu=1.5$ than for $\mu=0.5$. Figures~\ref{fig:survfun_pdf_zero__mu_1over2_A_100} and~\ref{fig:survfun_pdf_zero__mu_3over2_A_100} confirm this.

The ``heavy--tailness'' of the post-change distribution of the GSR stopping time can also be seen in Figures~\ref{fig:survfun_zero_A_100} which show the survival function $\Pr_{0}(\mathcal{S}_A^r\ge t)$. Once again, compared to the survival function in the pre-change regime, in the post-change regime the survival function decays down to zero more rapidly. However, the decay rate for $\mu=1.5$ is substantially higher than that for $\mu=0.5$. See Figures~\ref{fig:survfun_zero__mu_1over2_A_100} and~\ref{fig:survfun_zero__mu_3over2_A_100}, respectively. As before, part of the reason is that when $\mu=1.5$, i.e., for more contrast changes, the delay to detection is much lower than when $\mu=0.5$, i.e., for fainter changes. This causes the survival function appear as though it was ``pushed'' against the ``wall'' given by the vertical plane $t=0$.
\medskip

\begin{figure}[t!]
    \centering
    \begin{subfigure}{0.45\textwidth}
        \centering
        \includegraphics[width=\linewidth]{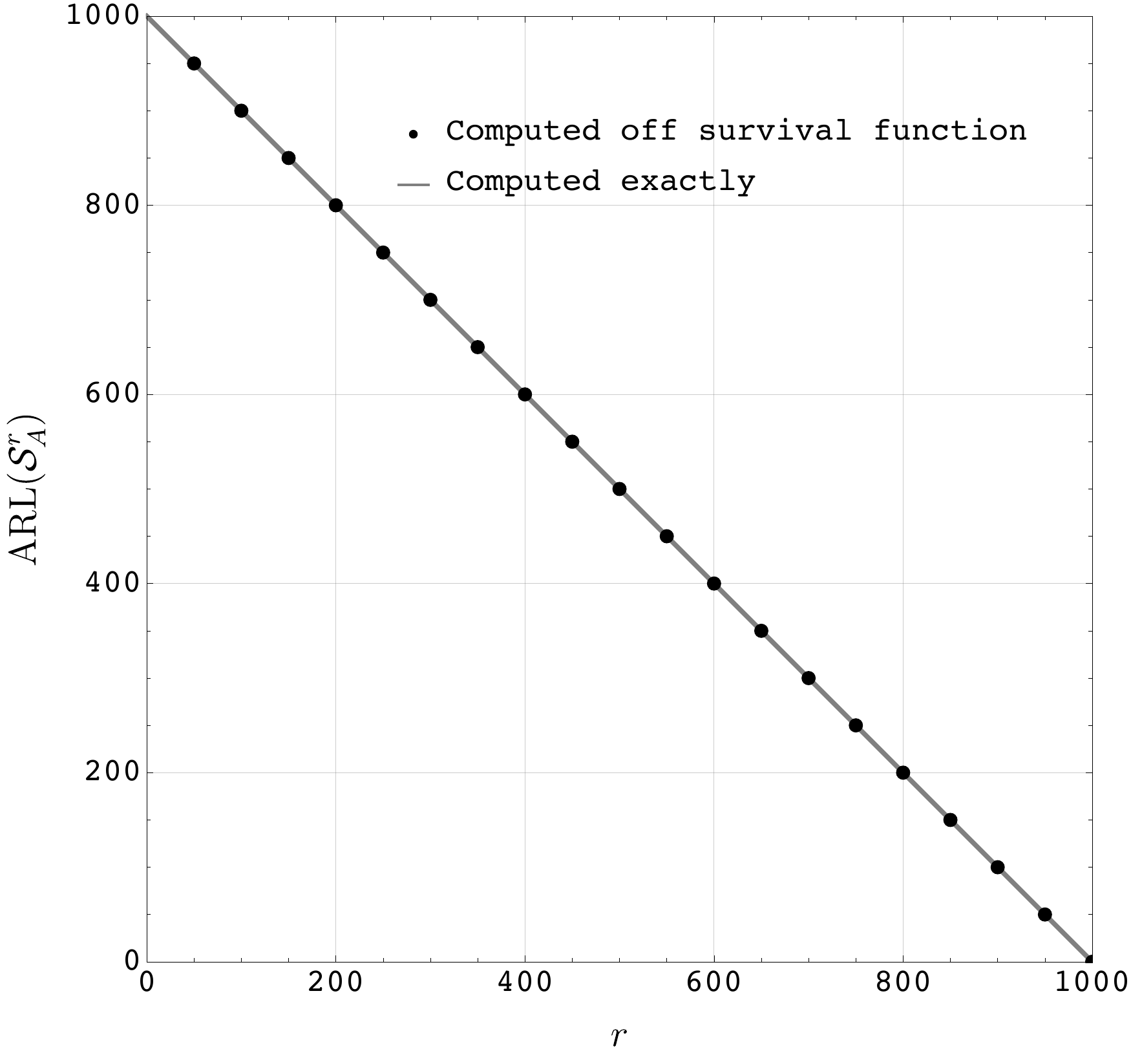}
        \caption{$\mu=0.5$.}
        \label{fig:firstmoment_inf__mu_1over2_A_1000}
    \end{subfigure}
    \hspace*{\fill}
    \begin{subfigure}{0.45\textwidth}
        \centering
        \includegraphics[width=\linewidth]{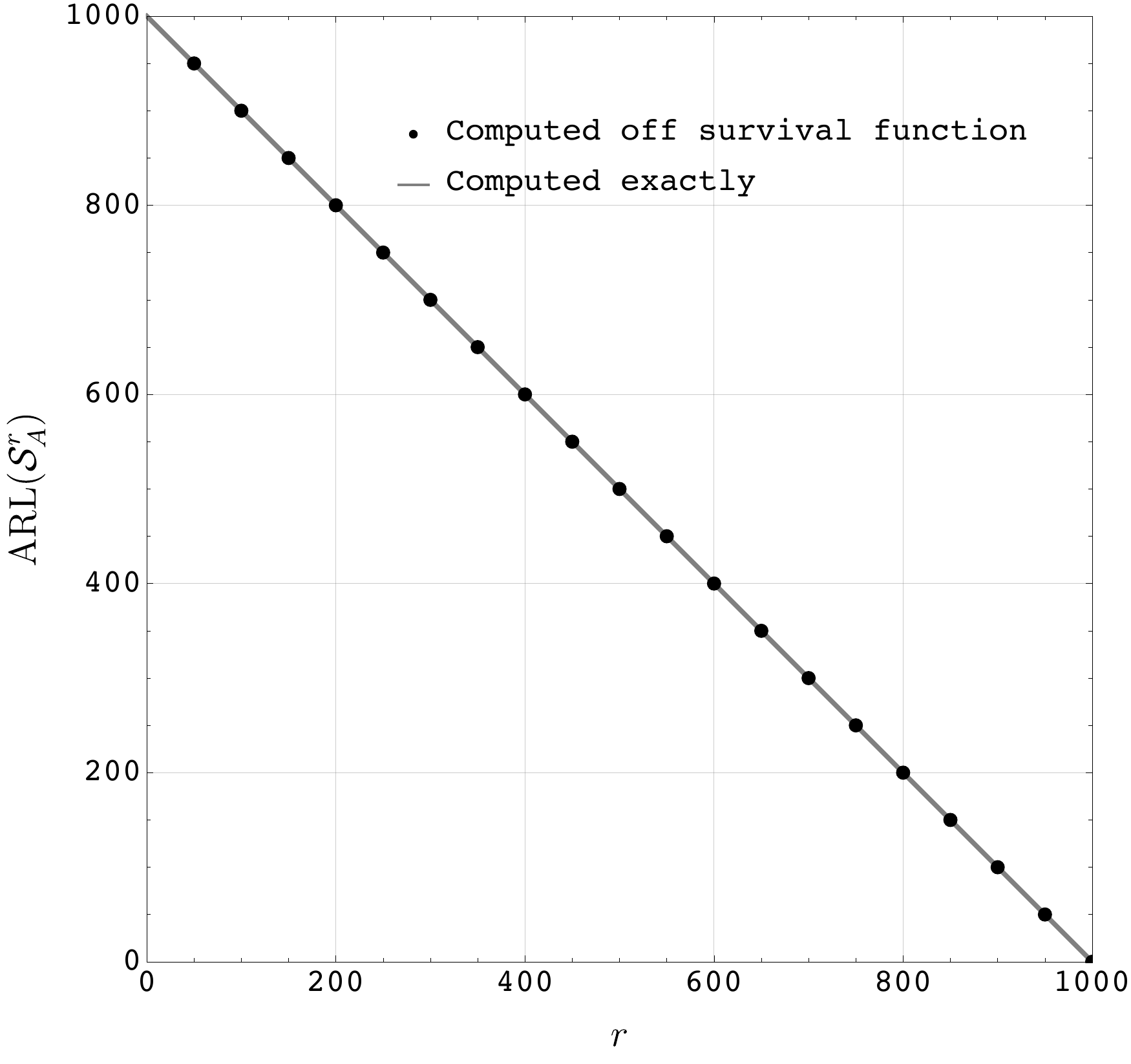}
        \caption{$\mu=1.5$.}
        \label{fig:firstmoment_inf__mu_3over2_A_1000}
    \end{subfigure}
    \caption{Pre-change first moment of $\mathcal{S}_A^r$, i.e., $\ARL(\mathcal{S}_A^r)\triangleq\EV_{\infty}[\mathcal{S}_A^r]$, as a function of $r\in[0,A]$ for $A=10^3$ and $\mu=\{0.5,1.5\}$.}
    \label{fig:firstmoment_inf_A_1000}
\end{figure}

\begin{figure}[t!]
    \centering
    \begin{subfigure}{0.45\textwidth}
        \centering
        \includegraphics[width=\linewidth]{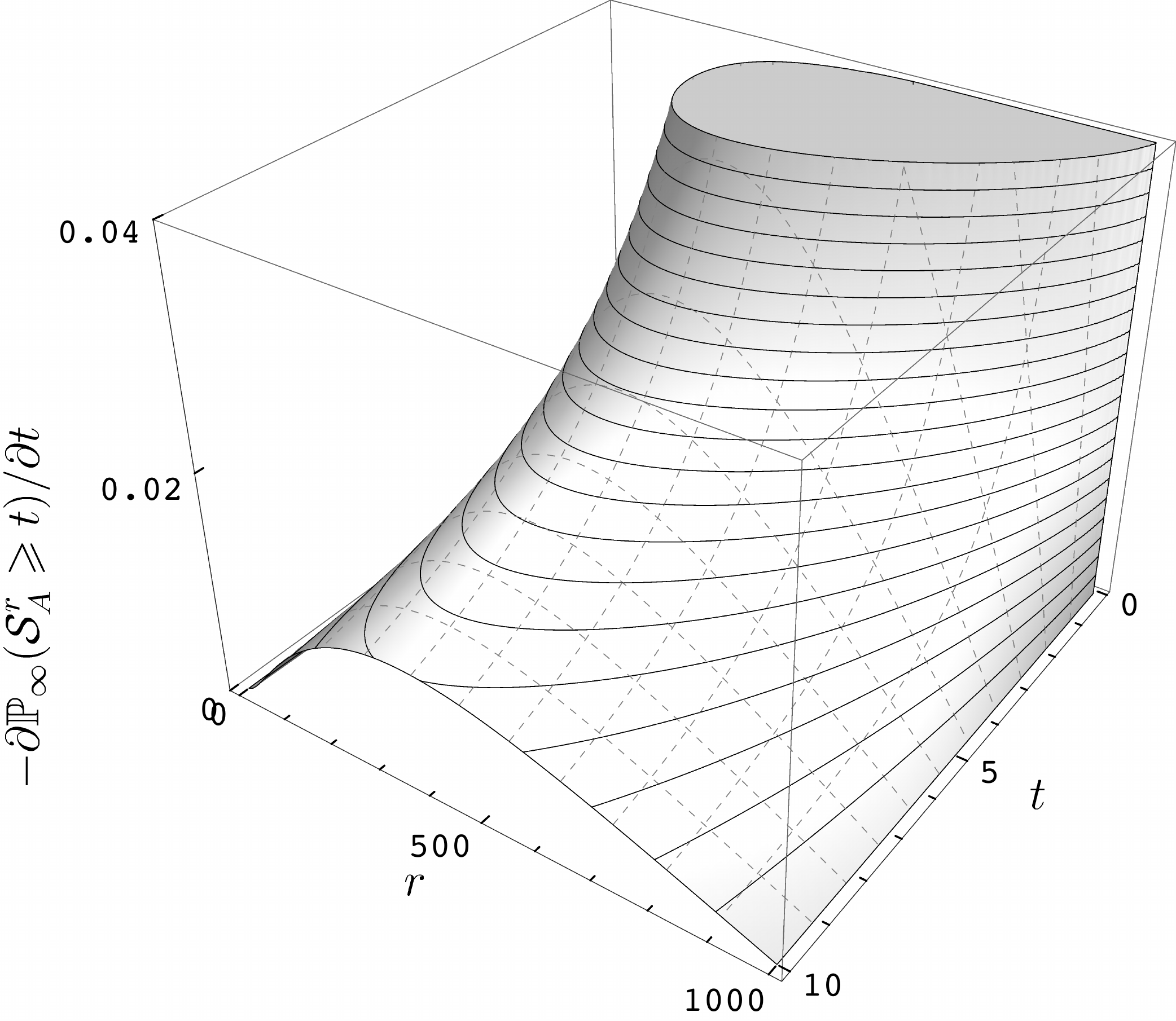}
        \caption{$\mu=0.5$.}
        \label{fig:survfun_pdf_inf__mu_1over2_A_1000}
    \end{subfigure}
    \hspace*{\fill}
    \begin{subfigure}{0.45\textwidth}
        \centering
        \includegraphics[width=\linewidth]{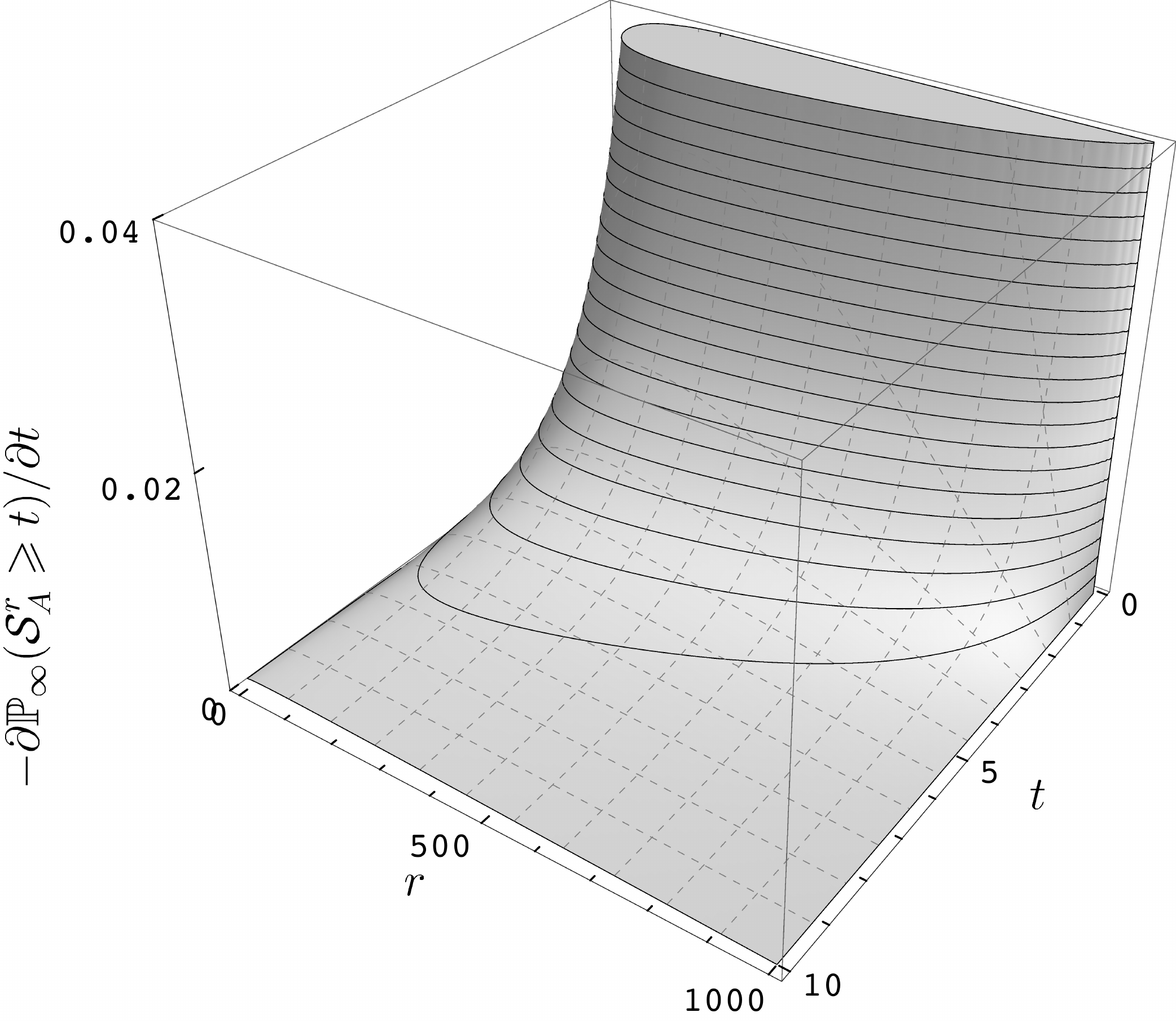}
        \caption{$\mu=1.5$.}
        \label{fig:survfun_pdf_inf__mu_3over2_A_1000}
    \end{subfigure}
    \caption{Pre-change survival function density $-\partial\Pr_{\infty}(\mathcal{S}_A^r\ge t)/\partial t$ as a function of $t\in[0,10]$ and $r\in[0,A]$ for $A=10^3$ and $\mu=\{0.5,1.5\}$.}
    \label{fig:survfun_pdf_inf_A_1000}
\end{figure}

\begin{figure}[t!]
    \centering
    \begin{subfigure}{0.45\textwidth}
        \centering
        \includegraphics[width=\linewidth]{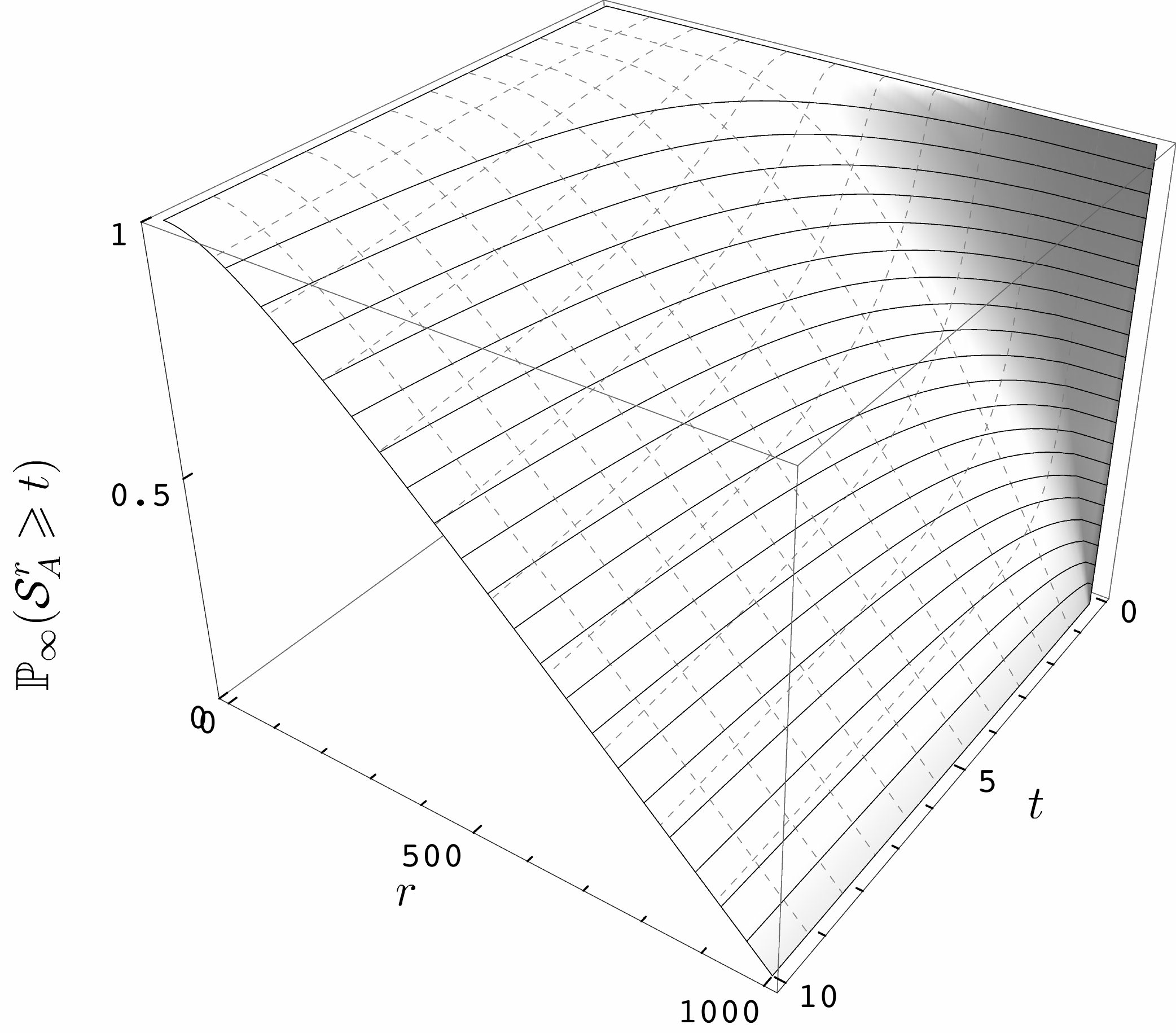}
        \caption{$\mu=0.5$.}
        \label{fig:survfun_inf__mu_1over2_A_1000}
    \end{subfigure}
    \hspace*{\fill}
    \begin{subfigure}{0.45\textwidth}
        \centering
        \includegraphics[width=\linewidth]{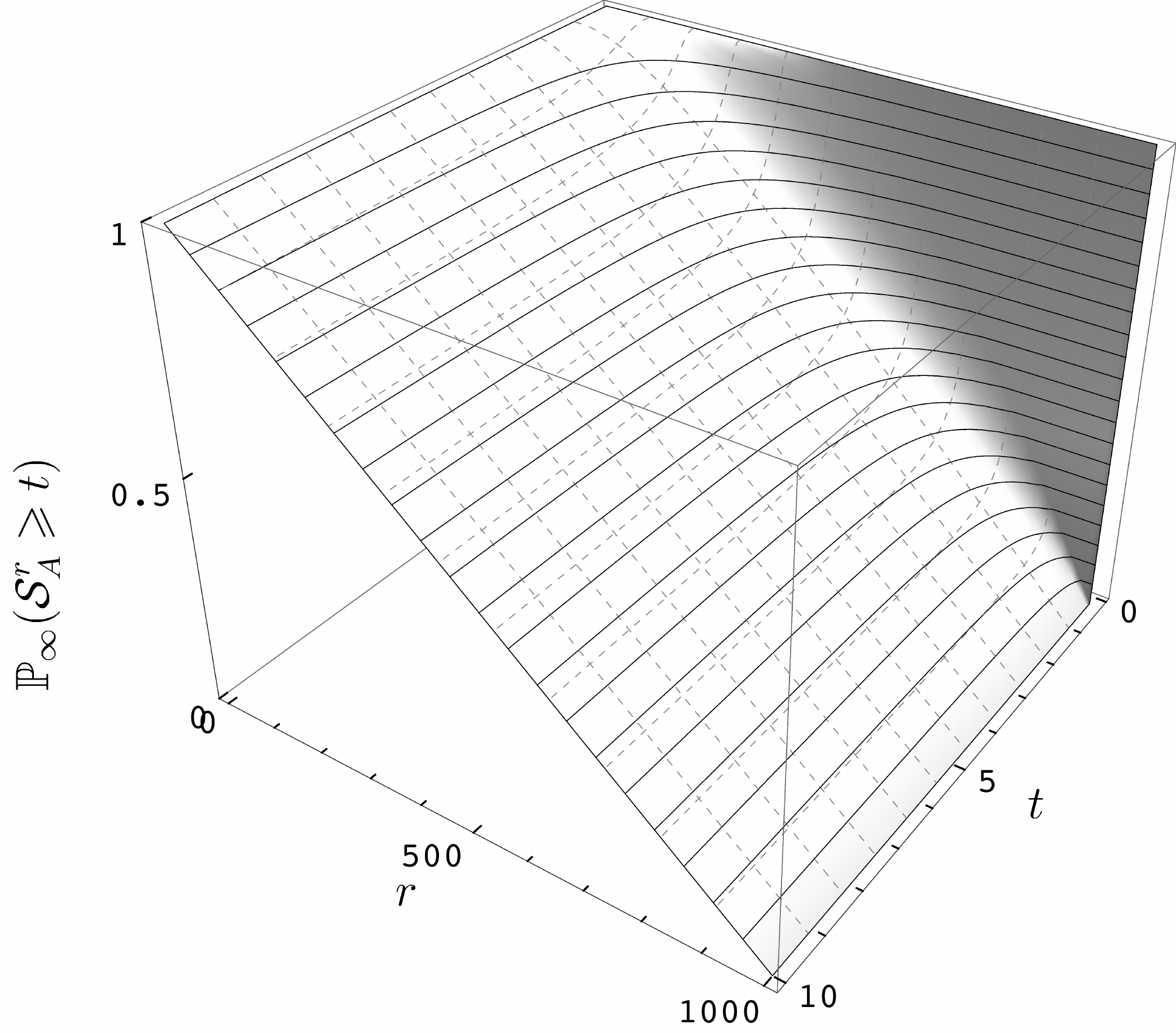}
        \caption{$\mu=1.5$.}
        \label{fig:survfun_inf__mu_3over2_A_1000}
    \end{subfigure}
    \caption{Pre-change survival function $\Pr_{\infty}(\mathcal{S}_A^r\ge t)$ as a function of $t\in[0,10]$ and $r\in[0,A]$ for $A=10^3$ and $\mu=\{0.5,1.5\}$.}
    \label{fig:survfun_inf_A_1000}
\end{figure}

Moving on, let us now consider the case when $A=10^3$. To that end, the first set of results is given in Figures~\ref{fig:firstmoment_inf_A_1000},~\ref{fig:survfun_pdf_inf_A_1000}, and~\ref{fig:survfun_inf_A_1000}, which all correspond to the pre-change regime. All the observations we made for the counterparts of these figures corresponding to the case when $A=10^2$ above immediately extend to these figures as well. However, due to the fact that the detection threshold is now higher, it can be seen in Figures~\ref{fig:survfun_pdf_inf_A_1000} and~\ref{fig:survfun_inf_A_1000} that, as time increases, the headstart ceases to matter quicker than when $A=10^2$. The reason is that the GSR statistic $(R_t^r)_{t\ge0}$ enters its quasi-stationary regime quicker for higher thresholds.
\medskip
\medskip

\begin{figure}[t!]
    \centering
    \begin{subfigure}{0.45\textwidth}
        \centering
        \includegraphics[width=\linewidth]{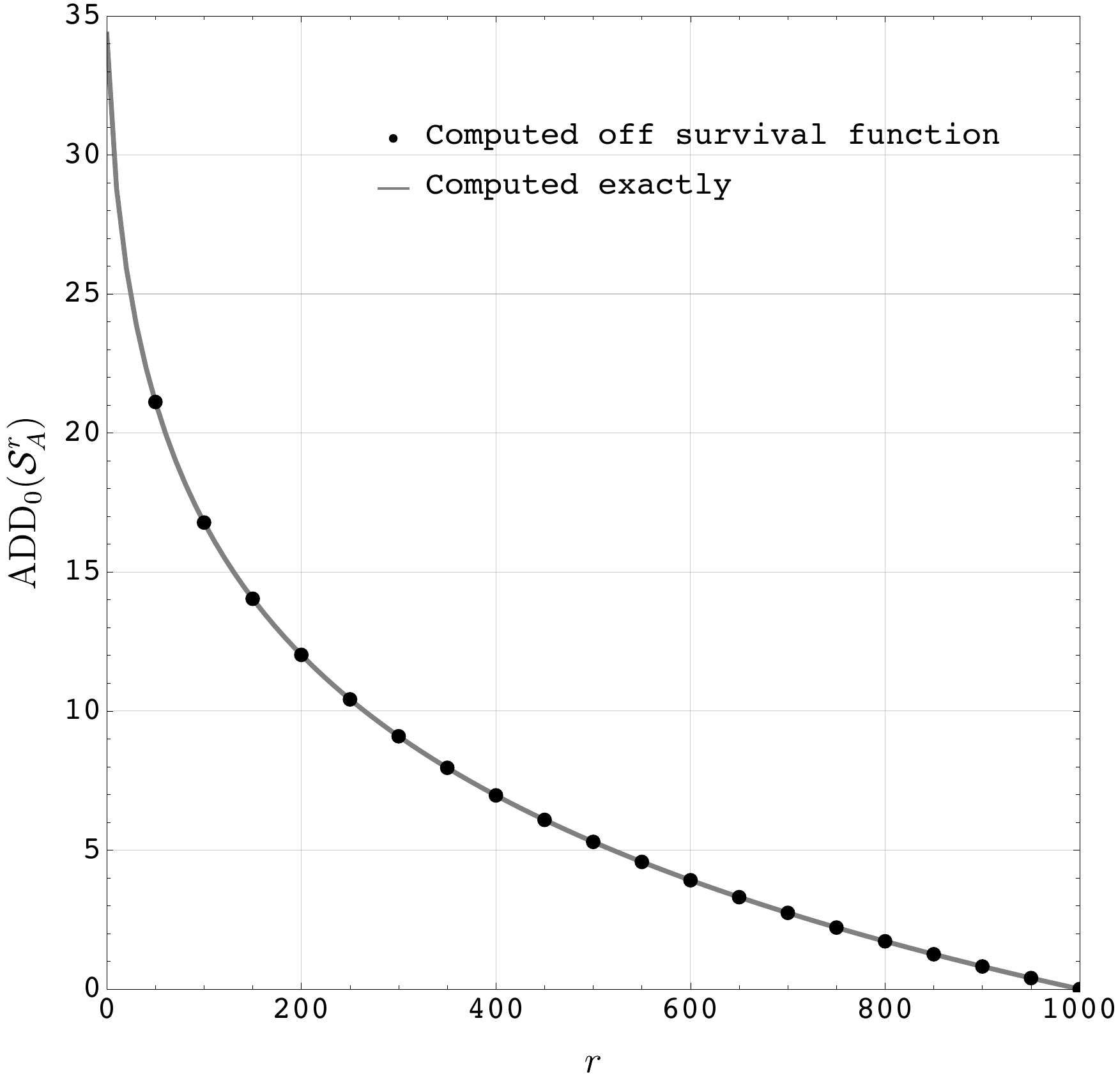}
        \caption{$\mu=0.5$.}
        \label{fig:firstmoment_zero__mu_1over2_A_1000}
    \end{subfigure}
    \hspace*{\fill}
    \begin{subfigure}{0.45\textwidth}
        \centering
        \includegraphics[width=\linewidth]{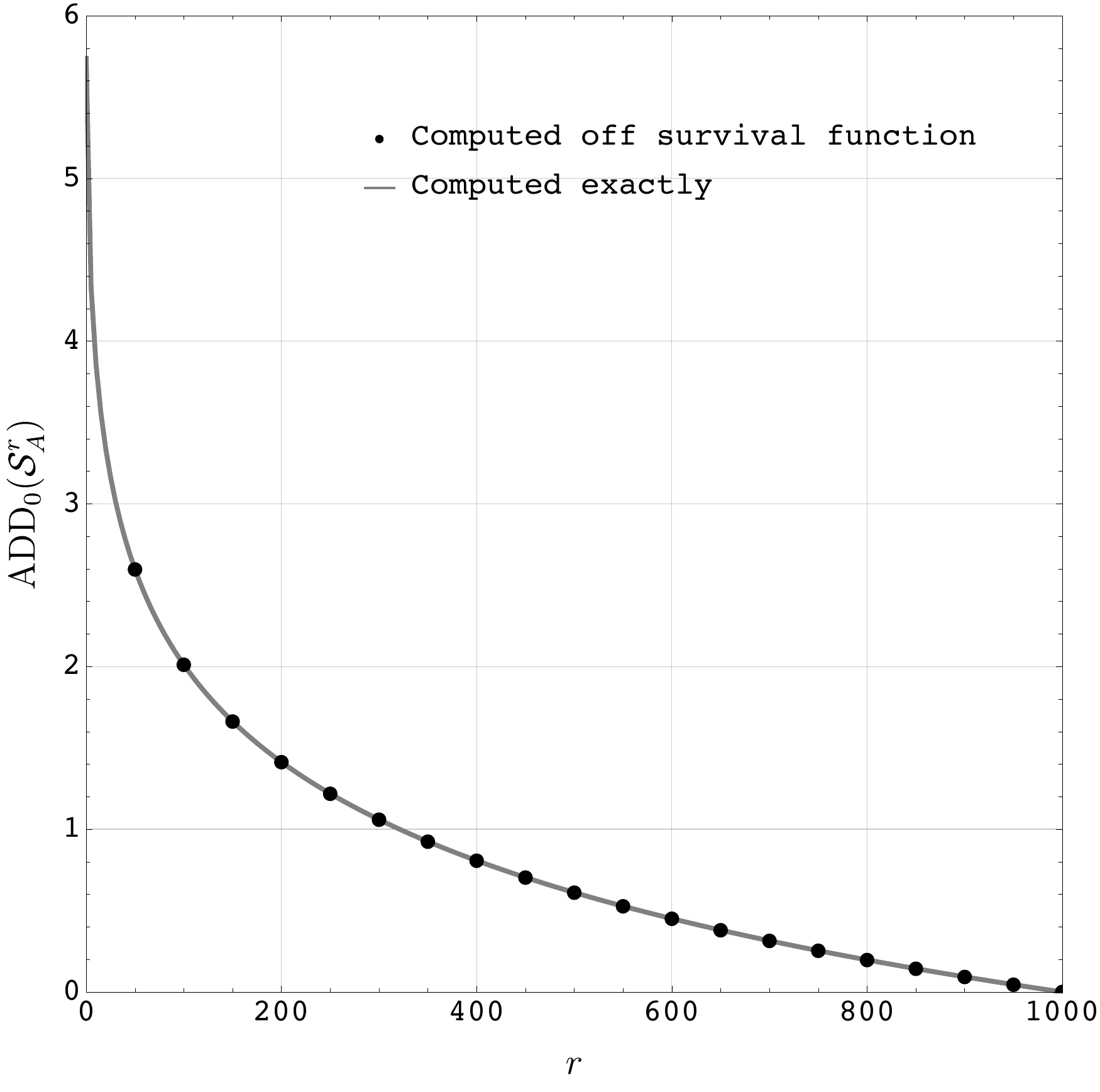}
        \caption{$\mu=1.5$.}
        \label{fig:firstmoment_zero__mu_3over2_A_1000}
    \end{subfigure}
    \caption{Post-change first moment of $\mathcal{S}_A^r$, i.e., $\ADD_{0}(\mathcal{S}_A^r)\triangleq\EV_{0}[\mathcal{S}_A^r]$, as a function of $r\in[0,A]$ for $A=10^3$ and $\mu=\{0.5,1.5\}$.}
    \label{fig:firstmoment_zero_A_1000}
\end{figure}

\begin{figure}[t!]
    \centering
    \begin{subfigure}{0.45\textwidth}
        \centering
        \includegraphics[width=\linewidth]{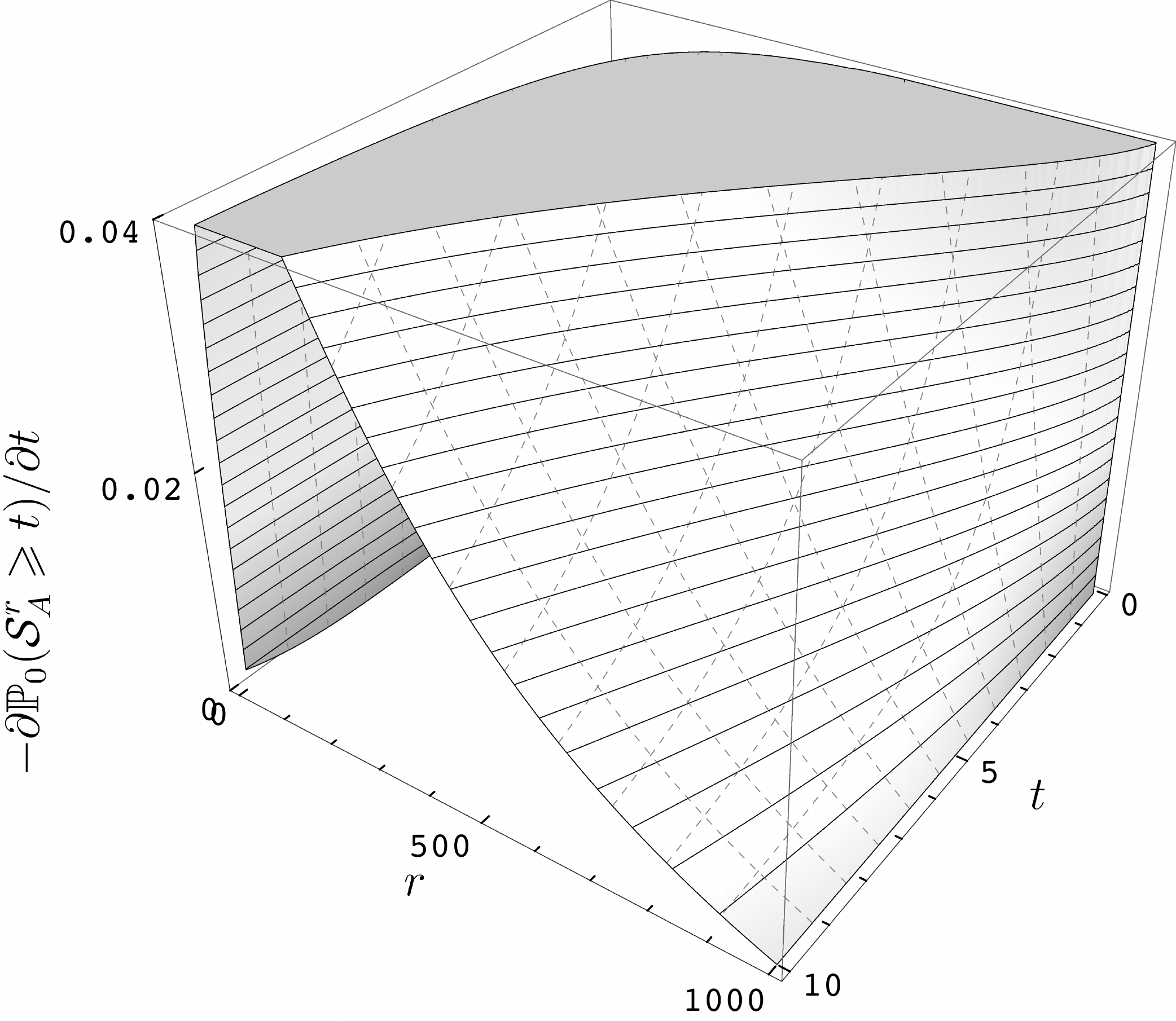}
        \caption{$\mu=0.5$.}
        \label{fig:survfun_pdf_zero__mu_1over2_A_1000}
    \end{subfigure}
    \hspace*{\fill}
    \begin{subfigure}{0.45\textwidth}
        \centering
        \includegraphics[width=\linewidth]{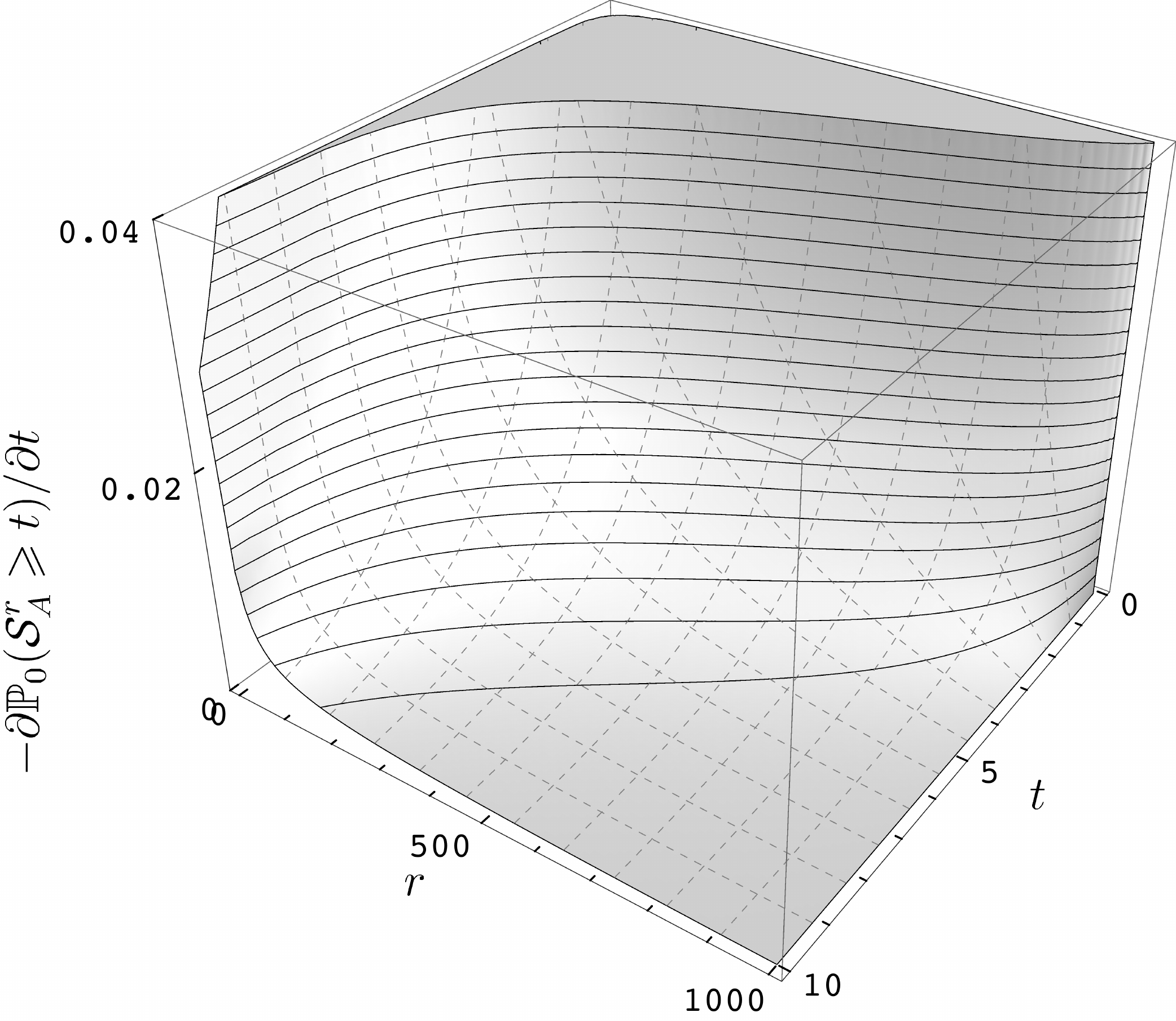}
        \caption{$\mu=1.5$.}
        \label{fig:survfun_pdf_zero__mu_3over2_A_1000}
    \end{subfigure}
    \caption{Post-change survival function density $-\partial\Pr_{0}(\mathcal{S}_A^r\ge t)/\partial t$ as a function of $t\in[0,10]$ and $r\in[0,A]$ for $A=10^3$ and $\mu=\{0.5,1.5\}$.}
    \label{fig:survfun_pdf_zero_A_1000}
\end{figure}

\begin{figure}[t!]
    \centering
    \begin{subfigure}{0.45\textwidth}
        \centering
        \includegraphics[width=\linewidth]{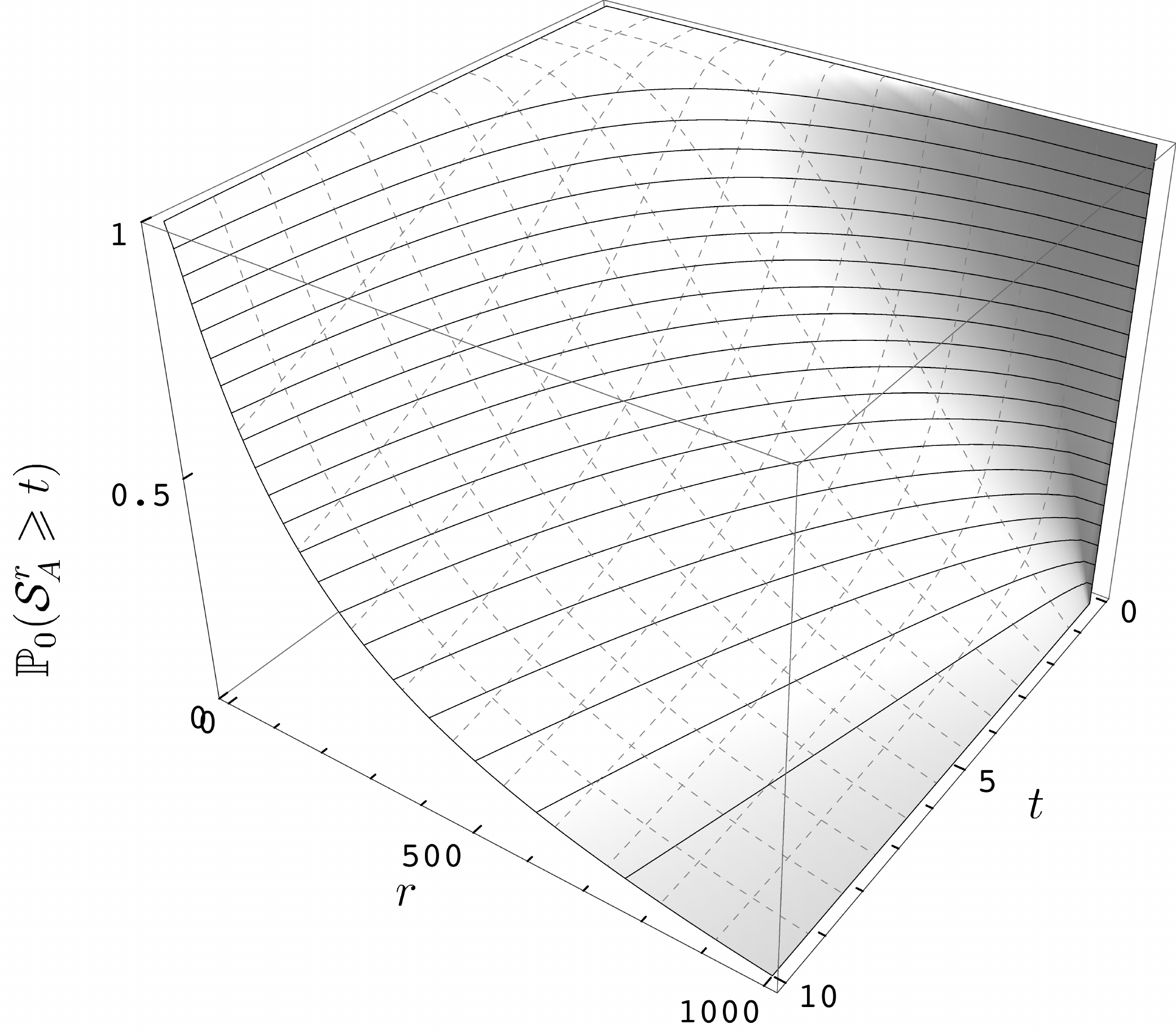}
        \caption{$\mu=0.5$.}
        \label{fig:survfun_zero__mu_1over2_A_1000}
    \end{subfigure}
    \hspace*{\fill}
    \begin{subfigure}{0.45\textwidth}
        \centering
        \includegraphics[width=\linewidth]{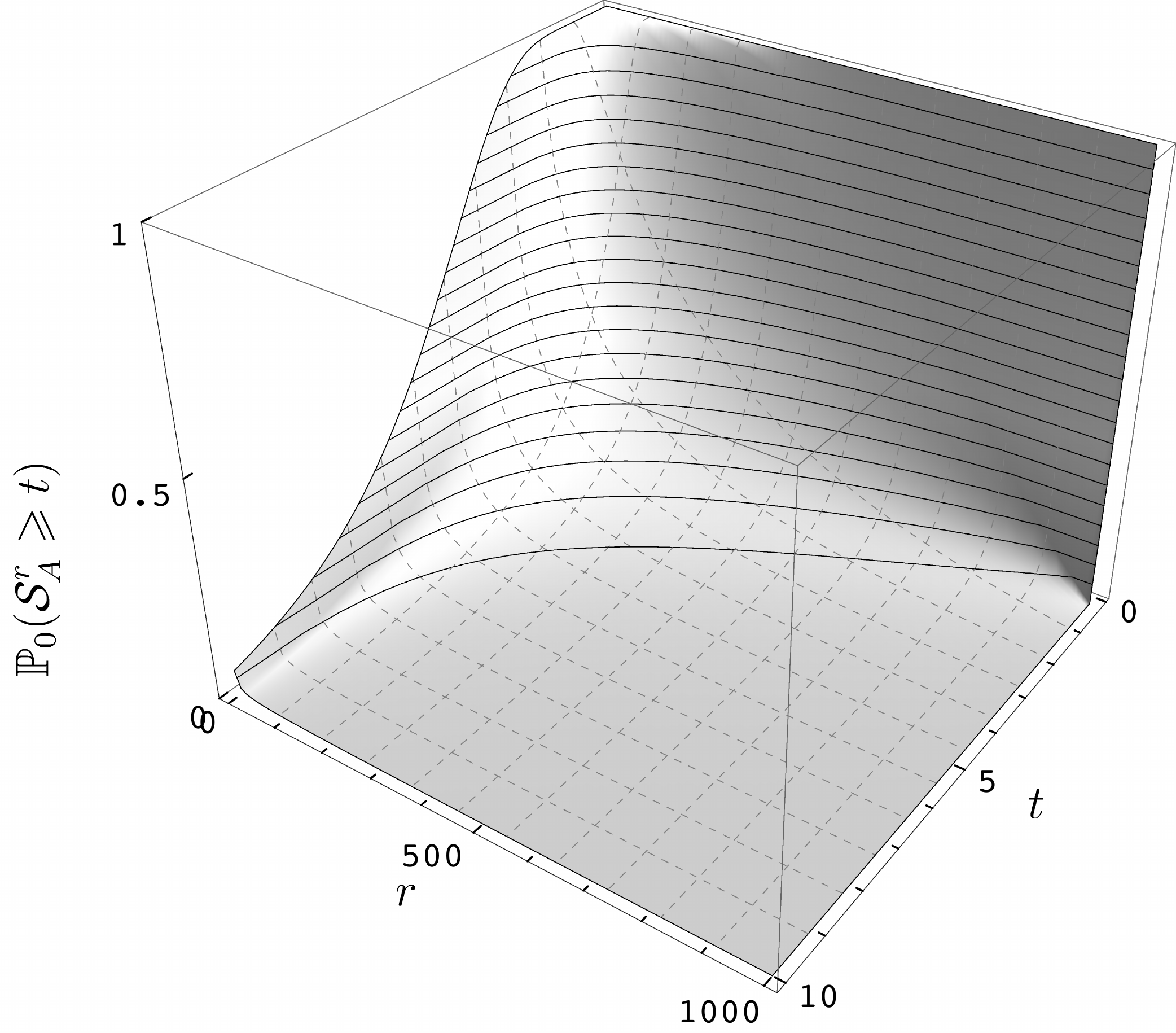}
        \caption{$\mu=1.5$.}
        \label{fig:survfun_zero__mu_3over2_A_1000}
    \end{subfigure}
    \caption{Post-change survival function $\Pr_{0}(\mathcal{S}_A^r\ge t)$ as a function of $t\in[0,10]$ and $r\in[0,A]$ for $A=10^3$ and $\mu=\{0.5,1.5\}$.}
    \label{fig:survfun_zero_A_1000}
\end{figure}

To conclude out numerical study, Figures~\ref{fig:firstmoment_zero_A_1000},~\ref{fig:survfun_pdf_zero_A_1000}, and~\ref{fig:survfun_zero_A_1000} are the post-change counterparts of Figures~\ref{fig:firstmoment_inf_A_1000},~\ref{fig:survfun_pdf_inf_A_1000}, and~\ref{fig:survfun_inf_A_1000}. Again, all the observations we made about the latter three figures can be also made about the former three figures.
\medskip
\medskip

\section{Conclusion}
\label{sec:conclusion}

This work sought to obtain as exhaustive a statistical characterization as possible of the stopping time associated with the Generalized Shiryaev--Roberts (GSR) procedure for quickest change-point detection under the classical minimax Brownian motion drift-shift scenario. Toward that goal, the main contribution of this paper is two exact closed-form formulae for the survival functions of the GSR stopping time, in the pre-drift regime and in the post-drift regime. The two formulae were found analytically, through direct solution of the respective Kolmogorov forward equations, and fully characterize the distribution of the GSR stopping time in the two regimes. On the more applied side, we put the two survival functions' formulae to work in software, and carried out a numerical study of the GSR stopping time's distribution in the two regimes. The study provided, apparently for the first time in the literature, a complete picture of the statistical profile of the GSR stopping time in the pre- and post-drift regimes.

\section*{Acknowledgements}

The author is thankful to the Editor-in-Chief, Nitis Mukhopadhyay (University of Connecticut--Storrs), and to the anonymous referee, whose constructive feedback provided on the first draft of the paper helped improve the quality of the manuscript and shape its current form. The author is also grateful to Dr. Grigory Sokolov (SUNY Binghamton) for the assistance provided with {\em Mathematica}.

\section*{Funding}

The author's effort was supported, in part, by the Simons Foundation via a Collaboration Grant in Mathematics under Award \#\,304574.


\end{document}